\documentclass[12pt]{article}
\usepackage{full page, 
amssymb, amscd, graphicx, 
}
    \title{{\bf Logarithmic
tensor category theory, IV: Constructions of tensor
product bifunctors and the compatibility conditions}}
    \author{Yi-Zhi Huang, James Lepowsky and Lin Zhang}
    \date{}

\begin{document}
    \bibliographystyle{alpha}
    \maketitle

    \newtheorem{rema}{Remark}[section]
    \newtheorem{propo}[rema]{Proposition}
    \newtheorem{theo}[rema]{Theorem}
   \newtheorem{defi}[rema]{Definition}
    \newtheorem{lemma}[rema]{Lemma}
    \newtheorem{corol}[rema]{Corollary}
     \newtheorem{exam}[rema]{Example}
\newtheorem{assum}[rema]{Assumption}
     \newtheorem{nota}[rema]{Notation}
        \newcommand{\ba}{\begin{array}}
        \newcommand{\ea}{\end{array}}
        \newcommand{\be}{\begin{equation}}
        \newcommand{\ee}{\end{equation}}
        \newcommand{\bea}{\begin{eqnarray}}
        \newcommand{\eea}{\end{eqnarray}}
        \newcommand{\nno}{\nonumber}
        \newcommand{\nn}{\nonumber\\}
        \newcommand{\lbar}{\bigg\vert}
        \newcommand{\p}{\partial}
        \newcommand{\dps}{\displaystyle}
        \newcommand{\bra}{\langle}
        \newcommand{\ket}{\rangle}
 \newcommand{\res}{\mbox{\rm Res}}
\newcommand{\wt}{\mbox{\rm wt}\;}
\newcommand{\swt}{\mbox{\scriptsize\rm wt}\;}
 \newcommand{\pf}{{\it Proof}\hspace{2ex}}
 \newcommand{\epf}{\hspace{2em}$\square$}
 \newcommand{\epfv}{\hspace{1em}$\square$\vspace{1em}}
        \newcommand{\ob}{{\rm ob}\,}
        \renewcommand{\hom}{{\rm Hom}}
\newcommand{\C}{\mathbb{C}}
\newcommand{\R}{\mathbb{R}}
\newcommand{\Z}{\mathbb{Z}}
\newcommand{\N}{\mathbb{N}}
\newcommand{\A}{\mathcal{A}}
\newcommand{\Y}{\mathcal{Y}}
\newcommand{\Arg}{\mbox{\rm Arg}\;}
\newcommand{\comp}{\mathrm{COMP}}
\newcommand{\lgr}{\mathrm{LGR}}

\newcommand{\dlt}[3]{#1 ^{-1}\delta \bigg( \frac{#2 #3 }{#1 }\bigg) }

\newcommand{\dlti}[3]{#1 \delta \bigg( \frac{#2 #3 }{#1 ^{-1}}\bigg) }

 \makeatletter
\newlength{\@pxlwd} \newlength{\@rulewd} \newlength{\@pxlht}
\catcode`.=\active \catcode`B=\active \catcode`:=\active \catcode`|=\active
\def\sprite#1(#2,#3)[#4,#5]{
   \edef\@sprbox{\expandafter\@cdr\string#1\@nil @box}
   \expandafter\newsavebox\csname\@sprbox\endcsname
   \edef#1{\expandafter\usebox\csname\@sprbox\endcsname}
   \expandafter\setbox\csname\@sprbox\endcsname =\hbox\bgroup
   \vbox\bgroup
  \catcode`.=\active\catcode`B=\active\catcode`:=\active\catcode`|=\active
      \@pxlwd=#4 \divide\@pxlwd by #3 \@rulewd=\@pxlwd
      \@pxlht=#5 \divide\@pxlht by #2
      \def .{\hskip \@pxlwd \ignorespaces}
      \def B{\@ifnextchar B{\advance\@rulewd by \@pxlwd}{\vrule
         height \@pxlht width \@rulewd depth 0 pt \@rulewd=\@pxlwd}}
      \def :{\hbox\bgroup\vrule height \@pxlht width 0pt depth
0pt\ignorespaces}
      \def |{\vrule height \@pxlht width 0pt depth 0pt\egroup
         \prevdepth= -1000 pt}
   }
\def\endsprite{\egroup\egroup}
\catcode`.=12 \catcode`B=11 \catcode`:=12 \catcode`|=12\relax
\makeatother

\def\hboxtr{\FormOfHboxtr} 
\sprite{\FormOfHboxtr}(25,25)[0.5 em, 1.2 ex] 

:BBBBBBBBBBBBBBBBBBBBBBBBB |
:BB......................B |
:B.B.....................B |
:B..B....................B |
:B...B...................B |
:B....B..................B |
:B.....B.................B |
:B......B................B |
:B.......B...............B |
:B........B..............B |
:B.........B.............B |
:B..........B............B |
:B...........B...........B |
:B............B..........B |
:B.............B.........B |
:B..............B........B |
:B...............B.......B |
:B................B......B |
:B.................B.....B |
:B..................B....B |
:B...................B...B |
:B....................B..B |
:B.....................B.B |
:B......................BB |
:BBBBBBBBBBBBBBBBBBBBBBBBB |

\endsprite

\def\shboxtr{\FormOfShboxtr} 
\sprite{\FormOfShboxtr}(25,25)[0.3 em, 0.72 ex] 

:BBBBBBBBBBBBBBBBBBBBBBBBB |
:BB......................B |
:B.B.....................B |
:B..B....................B |
:B...B...................B |
:B....B..................B |
:B.....B.................B |
:B......B................B |
:B.......B...............B |
:B........B..............B |
:B.........B.............B |
:B..........B............B |
:B...........B...........B |
:B............B..........B |
:B.............B.........B |
:B..............B........B |
:B...............B.......B |
:B................B......B |
:B.................B.....B |
:B..................B....B |
:B...................B...B |
:B....................B..B |
:B.....................B.B |
:B......................BB |
:BBBBBBBBBBBBBBBBBBBBBBBBB |

\endsprite


\begin{abstract}
This is the fourth part in a series of papers in which we introduce
and develop a natural, general tensor category theory for suitable
module categories for a vertex (operator) algebra.  In this paper
(Part IV), we give constructions of the $P(z)$- and $Q(z)$-tensor
product bifunctors using what we call ``compatibility conditions'' and
certain other conditions.
\end{abstract}


\tableofcontents
\vspace{2em}

In this paper, Part IV of a series of eight papers on logarithmic
tensor category theory, we give constructions of the $P(z)$- and
$Q(z)$-tensor product bifunctors using what we call ``compatibility
conditions'' and certain other conditions.  The sections, equations,
theorems and so on are numbered globally in the series of papers
rather than within each paper, so that for example equation (a.b) is
the b-th labeled equation in Section a, which is contained in the
paper indicated as follows: In Part I \cite{HLZ1}, which contains
Sections 1 and 2, we give a detailed overview of our theory, state our
main results and introduce the basic objects that we shall study in
this work.  We include a brief discussion of some of the recent
applications of this theory, and also a discussion of some recent
literature.  In Part II \cite{HLZ2}, which contains Section 3, we
develop logarithmic formal calculus and study logarithmic intertwining
operators.  In Part III \cite{HLZ3}, which contains Section 4, we
introduce and study intertwining maps and tensor product bifunctors.
The present paper, Part IV, contains Sections 5 and 6.  In Part V
\cite{HLZ5}, which contains Sections 7 and 8, we study products and
iterates of intertwining maps and of logarithmic intertwining
operators and we begin the development of our analytic approach.  In
Part VI \cite{HLZ6}, which contains Sections 9 and 10, we construct
the appropriate natural associativity isomorphisms between triple
tensor product functors.  In Part VII \cite{HLZ7}, which contains
Section 11, we give sufficient conditions for the existence of the
associativity isomorphisms.  In Part VIII \cite{HLZ8}, which contains
Section 12, we construct braided tensor category structure.

\paragraph{Acknowledgments}
The authors gratefully
acknowledge partial support {}from NSF grants DMS-0070800 and
DMS-0401302.  Y.-Z.~H. is also grateful for partial support {}from NSF
grant PHY-0901237 and for the hospitality of Institut des Hautes 
\'{E}tudes Scientifiques in the fall of 2007.

\renewcommand{\theequation}{\thesection.\arabic{equation}}
\renewcommand{\therema}{\thesection.\arabic{rema}}
\setcounter{section}{4}
\setcounter{equation}{0}
\setcounter{rema}{0}

\section{Constructions of the $P(z)$- and $Q(z)$-tensor product
bifunctors; the compatibility conditions}

In Section 4 we defined and studied the $P(z)$- and $Q(z)$-tensor
product bifunctors.  But in order to prove substantial results about
these bifunctors, notably, that under suitable conditions, they give
rise to braided tensor category structure, we will need to give a
useful, general construction of (models for) these bifunctors, when
they in fact exist.  This section is devoted to such constructions.
Our constructions hinge on subtle conditions, including
``compatibility conditions,'' on elements of certain dual spaces.

The results in this section are generalizations to the setting of the
present work of the constructions of the $P(z)$- and $Q(z)$-tensor
product bifunctors in \cite{tensor1}--\cite{tensor3}.  In the earlier
work \cite{tensor1}-\cite{tensor3} of the first two authors, the
$Q(z)$-tensor product of two modules was studied and developed first,
in \cite{tensor1} and \cite{tensor2}. The $P(z)$-tensor product was
then studied systematically in \cite{tensor3}, and many proofs for the
$P(z)$ case were given by using the results established for the $Q(z)$
case in \cite{tensor1} and \cite{tensor2}, rather than by carrying out
the subtle arguments in the $P(z)$ case itself, arguments that are
similar to (but different {}from) those in the $Q(z)$ case.  In the
present section and the next section, instead of following this
approach of \cite{tensor1}--\cite{tensor3}, we shall construct the
$P(z)$-tensor product and $Q(z)$-tensor product of two (generalized)
modules independently. In particular, even for the finitely reductive
case carried out in \cite{tensor1}--\cite{tensor3}, some of the
present results and proofs of the main theorems are completely new.
One new result is Proposition \ref{pz-comm} below, which was not
stated or proved (or needed) in the finitely reductive case in
\cite{tensor3}. This is proved below, by a direct argument in the
$P(z)$ setting, rather than by the use of the $Q(z)$
structure. Theorems \ref{comp=>jcb}, \ref{stable}, \ref{6.1} and
\ref{6.2} formulated below will be proved in the next section.  The
proofs of Theorems \ref{comp=>jcb} and \ref{stable} are new, even in
the finitely reductive case.

Recall the setup and conventions in Assumption \ref{assum}, including
the categories $\mathcal{M}_{sg}$ and $\mathcal{GM}_{sg}$.

We continue to let $z$ be a fixed nonzero complex number.

Later in this section (Assumption \ref{assum-c}), we shall take
$\mathcal{C}$ to be a full subcategory of $\mathcal{M}_{sg}$ or
$\mathcal{GM}_{sg}$, with certain additional properties.

\subsection{Affinizations of vertex algebras and the
opposite-operator map}

Just as in \cite{tensor1}--\cite{tensor3}, 
we shall use the Jacobi identity as a motivation 
to construct tensor products of
(generalized) $V$-modules in a suitable category. To do this, 
we need to study  various
``affinizations'' of a vertex algebra with respect to certain
algebras and vector spaces of formal Laurent series and formal 
rational functions.  The treatment of these matters below is 
very similar to that 
in \cite{tensor1}, but here we must take into account 
the gradings by $A$ and $\tilde{A}$. Here, as in Section 2 above,  
we are replacing the symbol $*$ for the 
``opposite-operator map''  in \cite{tensor1} by $o$. 

In Sections 5.2 and 5.3 below, we will be using the material in this
section to construct certain actions $\tau_{P(z)}$ and $\tau_{Q(z)}$
(see Definitions \ref{deftau} and \ref{deftauQ} as well as
(\ref{LP'(j)}) and (\ref{LQ'(j)})), in order to construct $P(z)$- and
$Q(z)$-tensor products.  These actions will be used to reformulate the
notions of $P(z)$- and $Q(z)$-intertwining maps (see Propositions
\ref{pz} and \ref{qz}), and these reformulations will enable us to
correspondingly reformulate the notions of $P(z)$- and $Q(z)$-tensor
products (see Corollaries \ref{tensorproductusingI'} and
\ref{q-tensorproductusingI'}).  This in turn will lead to the desired
constructions of $P(z)$- and $Q(z)$-tensor products when they in fact
exist (see Propositions \ref{tensor1-13.7} and \ref{tensor1-5.7} and
Theorems \ref{characterizationofbackslash} and
\ref{q-characterizationofbackslash}).

Let $(W, Y_{W})$ be a generalized $V$-module. We adjoin the
formal variable $t$ to our list of commuting formal variables. This
variable will play a special role.
Consider the vector spaces $$V[t,t^{-1}] = V \otimes {\C}[t,t^{-1}]
\subset V \otimes {\C} ((t)) \subset V\otimes {\C} [[t,t^{-1}]]
\subset V[[t,t^{-1}]] $$ (note carefully the distinction between the
last two, since $V$ is typically infinite-dimensional) and $W \otimes
{\C} \{t\}
\subset W \{t\}$
(recall (\ref{formalserieswithcomplexpowers})). The linear map
\begin{eqnarray}\label{tauW}
\tau_W: V[t,t^{-1}] &\to& \mbox{End} \;W\nno\\
v \otimes t^n &\mapsto&
v_n
\end{eqnarray}
($v\in V$, $n\in {\Z}$) extends canonically to 
\begin{eqnarray}\label{tauw}
\tau_W:& V \otimes {\C}((t)) &\to\; \mbox{End} \;W \nno\\
&v \otimes {\dps \sum_{n > N}} a_nt^n
&\mapsto\; \sum_{n > N} a_nv_n
\end{eqnarray}
 (but not to $V((t))$), in view of (\ref{set:wtvn}) and Assumption 
\ref{assum}. 
It further extends canonically to 
\begin{equation}
\tau_W: (V
\otimes {\C}((t)))[[x,x^{-1}]] \to (\mbox{End} \;W)[[x,x^{-1}]],
\end{equation}
where of course $(V
\otimes {\C}((t)))[[x,x^{-1}]]$ can be viewed as the subspace of
$V[[t,t^{-1},x,x^{-1}]]$ such that the coefficient of each power of
$x$ lies in $V \otimes {\C}((t))$.  

Let $v \in V$ and define the ``generic vertex operator''
\begin{equation}\label{3.4}
Y_t(v,x)  = \sum_{n \in {\Z}}(v \otimes
t^n)x^{-n-1} \in  (V \otimes {\C}[t,t^{-1}])[[x,x^{-1}]].
\end{equation}  
Then 
\begin{eqnarray}\label{3.5}
Y_t(v,x) & = & v\otimes x^{-1} \delta \left(\frac{t}{x}\right)\nno\\
 & = & v \otimes t^{-1} \delta \left(\frac{x}{t}\right)\nno\\
 & \in & V \otimes {\C}
[[t,t^{-1},x,x^{-1}]]\nno\\
 & (\subset & V[[t,t^{-1},x,x^{-1}]]),
\end{eqnarray} 
and the linear map 
\begin{eqnarray}\label{3.6}
V& \to &V\otimes {\C}[[t,t^{-1},x,x^{-1}]]\nno\\
v &\mapsto & Y_t(v,x)
\end{eqnarray}
 is simply the  map  given by tensoring by the
``universal element'' $x^{-1} \delta
\left(\frac{t}{x}\right)$.  We have 
\begin{equation}\label{3.7}
\tau_W(Y_t(v,x)) =
Y_W(v,x).
\end{equation}
 For all $f(x) \in {\C} [[x,x^{-1}]]$,
$f(x)Y_t(v,x)$ is defined and 
\begin{equation}
f(x)Y_t(v,x) = f(t) Y_t(v,x);
\end{equation}
it is crucial to keep in mind the delta-function substitution
principles (\ref{Xx1x2=Xx2x2}) and
(\ref{deltafunctionsubstitutionformula}), which we will be using
regularly.

In case $f(x) \in {\C}((x))$, then $\tau_{W}(f(x)Y_{t}(v, x))$ is
also defined, and
\begin{equation}\label{3.9}
f(x)Y_W(v,x)  =  f(x)\tau_W(Y_t(v,x)) = 
\tau_W(f(x)Y_t(v,x)) =  \tau_W(f(t)Y_t(v,x)).
\end{equation}
The expansion coefficients, in powers of $x$, of $Y_t(v,x)$
span $v \otimes {\C}[t,t^{-1}],$ the
$x$-expansion coefficients
of $Y_W(v,x)$ span $\tau_W(v \otimes {\C}[t,t^{-1}])$ and 
for $f(x)\in \C[[x, x^{-1}]]$, the
$x$-expansion coefficients of $f(x)Y_t(v,x)$ span $v \otimes f(t)
{\C} [t,t^{-1}]$. In case $f(x) \in {\C}((x))$, the 
$x$-expansion coefficients of $f(x)Y_W(v,x)$ span $\tau_W(v \otimes
f(t) {\C}[t,t^{-1}])$.
 
Using this viewpoint, we shall examine each of the three terms in
the Jacobi identity (\ref{log:jacobi}) in the definition of logarithmic
intertwining operator. First we consider the formal Laurent series in
$x_{0}$, $x_{1}$, $x_{2}$ and $t$ given by
\begin{eqnarray}\label{3.10}
{\dps x^{-1}_2 \delta\left(\frac{x_1-x_0}{x_2}\right)
Y_t(v,x_0) }&=& {\dps x^{-1}_1\delta \left(\frac{x_2+x_0}{x_1}\right)
Y_t(v,x_0)}\nno\\
&=&{\dps  v \otimes x^{-1}_1\delta\left(\frac{x_2 +
t}{x_1}\right)x^{-1}_0 \delta\left(\frac{t}{x_0}\right)}
\end{eqnarray}
(cf. the right-hand side of (\ref{log:jacobi})). The
expansion coefficients in powers of $x_0$, $x_1$ and $x_2$ of (\ref{3.10}) 
span just the space $v \otimes {\C}[t,t^{-1}]$. However,
the expansion coefficients in $x_0$ and $x_1$ only (but not in $x_{2}$) of 
\begin{eqnarray}\label{3.11}
x^{-1}_1\delta\left(\frac{x_{2}+x_0}{x_1}\right) Y_t(v,x_0)
&  =& v \otimes x^{-
1}_1\delta\left(\frac{x_{2}+t}{x_1}\right)x^{-
1}_0\delta\left(\frac{t}{x_0}\right)\nno\\ 
&  =& v \otimes
\left(\sum_{m \in {\Z}} (x_{2} + t)^m x^{-m-
1}_1\right)\left(\sum_{n \in {\Z}}t^n x^{-n-
1}_0\right)
\end{eqnarray}
 span $$v \otimes \iota_{x_{2},t}{\C}[t,t^{- 1}, x_{2}+t,
(x_{2}+t)^{-1}]\subset v \otimes {\C}[x_{2}, x_{2}^{-1}]((t)),$$
where $\iota_{x_{2}, t}$ is the operation of expanding a formal rational 
function in the indicated algebra as a formal Laurent series
involving only finitely many negative powers of $t$ (cf.
the notation $\iota_{12}$, etc., considered at the end of Section 2). 
We shall use similar
$\iota$-notations below. Specifically, the coefficient of 
$x^{-n-1}_0 x^{-m-1}_1$ $(m,n
\in {\Z})$ in (\ref{3.11}) is $v \otimes (x_{2} + t)^m t^n$.  

We may
specialize $x_2 \mapsto z\in {\C}^{\times}$, and (\ref{3.11}) becomes
\begin{eqnarray}\label{3.12}
z^{-1}\delta\left(\frac{x_{1}-x_0}{z}\right) Y_t(v,x_0)
&=&x^{-1}_1\delta\left(\frac{z+x_0}{x_1}\right) Y_t(v,x_0)\nno\\
&  =& v \otimes x^{-
1}_1\delta\left(\frac{z+t}{x_1}\right)x^{-
1}_0\delta\left(\frac{t}{x_0}\right)\nno\\ 
&  =& v \otimes
\left(\sum_{m \in {\Z}} (z + t)^m x^{-m-
1}_1\right)\left(\sum_{n \in {\Z}}t^n x^{-n-
1}_0\right).
\end{eqnarray}
The coefficient of $x^{-n-1}_0 x^{-m-1}_1$ $(m,n \in {\Z})$ in (\ref{3.12})
is $v \otimes (z + t)^m t^n\in V\otimes {\C}((t))$, and these
coefficients span
\begin{equation}\label{3.13}
v\otimes {\C}[t, t^{-1}, (z+t)^{-1}]\subset v\otimes {\C}((t)).
\end{equation}
Our $Q(z)$-tensor product construction in Section 5.3 below
will be based on a certain action of the 
space $V\otimes {\C}[t, t^{-1}, (z+t)^{-1}]$, and the description of 
this space as the span of the coefficients of the expression (\ref{3.12}) (as
$v\in V$ varies) will be very useful.

Now consider 
\begin{eqnarray}\label{3.14}
x^{-1}_0 \delta\left(\frac{-x_2 +
x_1}{x_0}\right) Y_t(v,x_1)& =&  v
\otimes x^{-1}_0 \delta\left(\frac{-x_2 + t}{x_0}\right) x^{-1}_1
\delta \left(\frac{t}{x_1}\right)\nno\\
&=&{\dps    v\otimes \biggr(\sum_{n \in {\Z}} (-x_2 + t)^n x^{-n-
1}_0\biggr)\biggr(\sum_{m \in {\Z}} t^m x^{-m-1}_1\biggr)}
\end{eqnarray}
(cf. the second term on the left-hand side of (\ref{log:jacobi})). 
The expansion coefficients in powers of $x_0$ and $x_1$ (but not $x_2$)
span 
$$v \otimes \iota_{x_{2}, t}{\C}[t,t^{-1},-x_2 + t, (-x_2 +
t)^{-1}],$$ and in fact the coefficient of $x^{-n-1}_0 x^{-m-1}_1$
$(m,n \in {\Z})$ in (\ref{3.14}) is $v \otimes (-x_2 + t)^n t^m$. Again
specializing $x_2 \mapsto z\in {\C}^{\times}$, we obtain
\begin{eqnarray}\label{3.15}
x^{-1}_0 \delta \left(\frac{-z+x_1}{x_0}\right)
Y_t(v,x_1) 
&=& v\otimes x^{-1}_0 \delta \left(\frac{-
z+t}{x_0}\right) x^{-1}_1 \delta \left(\frac{t}{x_1}\right)\nn
& =& v \otimes \biggr(\sum_{n \in {\Z}} (-z+t)^n x^{-n- 1}_0\biggr)
\biggr(\sum_{m \in
{\Z}}t^m x^{-m-1}_1\biggr).
\end{eqnarray}
 The coefficient of $x^{-n-1}_0 x^{- m-1}_1$ $(m,n \in {\Z})$ in
(\ref{3.15}) is $v \otimes (-z+t)^n t^m$, and these coefficients span
\begin{equation}\label{3.16}
v\otimes {\C}[t, t^{-1}, (-z+t)^{-1}]\subset v\otimes {\C}((t)).
\end{equation}

Finally, consider 
\begin{eqnarray}
{\dps x^{-1}_0 \delta \left(\frac{x_1 -
x_2}{x_0}\right) Y_t (v,x_1)}
&=&{\dps v \otimes x^{-1}_0
\delta \left(\frac{t-x_2}{x_0}\right) x^{-1}_1 \delta
\left(\frac{t}{x_1}\right)}
\end{eqnarray}
(cf. the first term on the left-hand side of (\ref{log:jacobi})).
The coefficient of $x_{0}^{-n-1}x_{1}^{-m-1}$ ($m, n\in {\Z}$) is
$v\otimes (t-x_{2})^{n}t^{m}$, and these expansion coefficients  span
$$v\otimes \iota_{t, x_{2}}{\C}[t, t^{-1}, t-x_{2}, (t-x_{2})^{-1}].$$
If we again specialize $x_2
\mapsto z$, we get
\begin{equation}\label{3.18}
x^{-1}_0\delta\left(\frac{x_1-z}{x_0}\right)Y_t(v,x_1) = v\otimes
x_{0}^{-1}\delta\left(\frac{t-z}{x_{0}}\right)x_{1}^{-1}\delta\left(
\frac{t}{x_{1}}\right),
\end{equation}
whose coefficient of $x^{-n-1}_0 x^{-m-1}_1$  is $v
\otimes (t - z)^n t^m$.  These coefficients span 
\begin{equation}\label{3.19}
v\otimes {\C}[t, t^{-1}, (t-z)^{-1}]\subset v\otimes {\C}((t^{-1}))
\end{equation}
(cf. (\ref{3.13}), (\ref{3.16})).

In the construction of $P(z)$-tensor products in Section 5.2, we 
shall also need the following expression, 
which is slightly different {}from what we have just analyzed:
\begin{eqnarray}\label{y-t-delta}
{\dps x^{-1}_0 \delta \left(\frac{x_1^{-1} -
x_2}{x_0}\right) Y_t (v, x_1)}
&=&{\dps v \otimes x^{-1}_0
\delta \left(\frac{t^{-1}-x_2}{x_0}\right) x^{-1}_1 \delta
\left(\frac{t}{x_1}\right).}
\end{eqnarray}
The coefficient of $x_{0}^{-n-1}x_{1}^{-m-1}$ ($m, n\in {\Z}$) is
$v\otimes (t^{-1}-x_{2})^{n}t^{m}$, and these expansion coefficients  span
$$v\otimes \iota_{t^{-1}, x_{2}}{\C}[t, t^{-1}, t^{-1}-x_{2}, (t^{-1}-x_{2})^{-1}].$$
If we again specialize $x_2
\mapsto z$, we get
\begin{equation}\label{3.18-1}
x^{-1}_0\delta\left(\frac{x_1^{-1}-z}{x_0}\right)Y_t(v,x_1) = v\otimes
x_{0}^{-1}\delta\left(\frac{t^{-1}-z}{x_{0}}\right)x_{1}^{-1}\delta\left(
\frac{t}{x_{1}}\right),
\end{equation}
whose coefficient of $x^{-n-1}_0 x^{-m-1}_1$  is $v
\otimes (t^{-1} - z)^n t^m$.  These coefficients span 
\begin{equation}\label{3.19-1}
v\otimes {\C}[t, t^{-1}, (t^{-1}-z)^{-1}]
=v\otimes {\C}[t, t^{-1}, (z^{-1}-t)^{-1}]\subset v\otimes {\C}((t)).
\end{equation}
Our $P(z)$-tensor product construction in Section 5.2 below
will be based on a certain action of the 
space $V\otimes {\C}[t, t^{-1}, (z^{-1}-t)^{-1}]$ (see Definition
\ref{deftau} below).

In fact, we shall be evaluating the Jacobi identity
(\ref{log:jacobi}), or more specifically, (\ref{im:def}), on the
elements of the contragredient module $W'_{3}$, and this will in
particular allow us to convert the expansion (\ref{3.19}) into an
expansion in positive powers of $t$, as in the previous paragraph; see
(\ref{im:def'}) and Definition \ref{deftau}.  For describing the
action given in Definition \ref{deftau}, it will be useful to examine
the notions of opposite and contragredient vertex operators $Y^o$ and
$Y'$ more closely (recall Section 2, in particular, (\ref{yo}),
(\ref{v'vo}) and Theorem \ref{set:W'}).

We shall interpret the opposite vertex operator map $Y^{o}_{W}$ by means
of an operation on $V\otimes {\C}[[t, t^{-1}]]$ that will convert
vertex operators into their opposites.  We shall write this
``opposite-operator'' map, in various contexts, as ``$o$.''  The
operation $o$ will be an involution.  We proceed as follows: First we
generalize $Y^{o}$ in the following way: Recall that 
by Assumption \ref{assum}, $L(1)$ acts nilpotently on 
any element $v\in V$. In particular, $e^{xL(1)}v$ is a polynomial
in the formal variable $x$. Given any vector space $U$
and any linear map
\begin{eqnarray}
Z(\cdot,x):\ V & \rightarrow &  U[[x,x^{-
1}]]\;\;\biggr(=\prod_{n \in {\Z}} U \otimes x^n\biggr)\nno\\ 
v & \mapsto & Z(v,x)
\end{eqnarray}
{from} $V$ into $U[[x,x^{-1}]]$
(i.e., given any family of linear maps {from} $V$ into the spaces
$U \otimes x^n$), we define $Z^o(\cdot,x):\ V \to U[[x,x^{-1}]]$
by 
\begin{equation}\label{3.21}
Z^o(v,x) = Z(e^{xL(1)}(-x^{-2})^{L(0)}v,x^{-1}),
\end{equation}
 where
we use the obvious linear map $Z(\cdot,x^{-1}):\ V \to U[[x,x^{-
1}]],$ and where we extend $Z(\cdot,x^{-1})$ canonically to a
linear map $Z(\cdot,x^{-1}):\ V[x,x^{-1}] \to U[[x,x^{-1}]].$
Then by formula (5.3.1) in \cite{FHL}  (the proof of Proposition
5.3.1), we have
\begin{eqnarray}\label{Zoo}
Z^{oo}(v,x) & = & Z^o(e^{xL(1)}(-x^{-
2})^{L(0)}v,x^{-1})\nno\\ 
& = & Z(e^{x^{-1}L(1)}(-
x^2)^{L(0)}e^{xL(1)}(-x^{-2})^{L(0)}v,x)\nno\\
&=& Z(v,x).
\end{eqnarray}
That is,
\begin{equation}
Z^{oo}(\cdot,x) = Z(\cdot,x).
\end{equation}
 Moreover, if
$Z(v,x) \in U((x))$, then $Z^o(v,x) \in U((x^{-1}))$ and vice
versa.  
 
Now we expand $Z(v,x)$ and $Z^o(v,x)$ in components. Write
\begin{equation}
Z(v,x) = \sum_{n \in {\Z}}v_{(n)} x^{-n-1},
\end{equation}
 where for all $n\in {\Z}$,
\begin{eqnarray}
V & \to & U\nno\\
v & \mapsto & v_{(n)}
\end{eqnarray}
 is a linear map depending on $Z(\cdot, x)$ (and in fact, as
$Z(\cdot, x)$ varies, these linear maps are arbitrary).
Also write 
\begin{equation}
Z^o(v,x)
= \sum_{n \in {\Z}} v^o_{(n)}x^{-n-1}
\end{equation}
 where
\begin{eqnarray}
V & \to & U\nno\\
v & \mapsto & v^o_{(n)}
\end{eqnarray}
 is a linear map depending on $Z(\cdot, x)$.  We shall compute $v_{(n)}^{o}$.
First note that
\begin{equation}\label{3.32}
\sum_{n
\in {\Z}}v^o_{(n)}x^{-n-1} = \sum_{n \in {\Z}}(e^{xL(1)}(-
x^{-2})^{L(0)}v)_{(n)}x^{n+1}.
\end{equation}
  For convenience, suppose that $v
\in V_{(h)}$, for $h \in {\Z}$.  Then the right-hand side of (\ref{3.32}) 
is equal to
\begin{eqnarray}
\lefteqn{(-1)^h\sum_{n \in {\Z}}(e^{xL(1)}v)_{(-
n)}x^{-n+1-2h}}\nno\\
&&=  (-1)^h \sum_{n \in {\Z}} \sum_{m \in 
{\N}} \frac{1}{m!}(L(1)^m v)_{(-n)}x^{m-n+1-2h}\nno\\
&&= 
(-1)^h\sum_{m
\in {\N}} \frac{1}{m!} \sum_{n \in {\Z}}(L(1)^mv)_{(-n-m-
2+2h)}x^{-n-1}, 
\end{eqnarray} 
that is, 
\begin{equation}\label{vo}
v^o_{(n)} = (-
1)^h\sum_{m \in {\N}}\frac{1}{m!} (L(1)^mv)_{(-n-m-2+2h)}.
\end{equation}
(Recall that by Assumption \ref{assum}, $L(1)^{m}v=0$ when  $m$ 
is sufficiently large, so that these expressions 
are well defined.)
For $v\in V$ not necessarily homogeneous, $v^o_{(n)}$ is given by the
appropriate sum of such expressions.

Now consider the special case where $U = V \otimes {\C}[t,t^{-1}]$ and
where $Z(\cdot, x)$ is the ``generic'' linear map
\begin{eqnarray}
Y_t(\cdot,x): V & \rightarrow & (V \otimes
{\C}[t,t^{-1}])[[x,x^{-1}]]\nno\\
v & \mapsto & Y_t(v,x) = \sum_{n
\in {\Z}}(v \otimes t^n)x^{-n-1}
\end{eqnarray} 
(recall (\ref{3.4})), i.e.,
\begin{equation}
v_{(n)} = v \otimes t^n.
\end{equation}
Then for $v\in V_{(h)}$, 
\begin{equation}
v^o_{(n)} = (-1)^h\sum_{m \in {\N}}
\frac{1}{m!}((L(1))^mv) \otimes t^{-n-m-2+2h}
\end{equation}
 in this case.
 
This motivates defining an $o$-operation on $V \otimes {\C}[t,t^{-1}]$
as follows: For any $n, h\in {\Z}$ and $v\in V_{(h)}$, define
\begin{equation}\label{3.38}
(v\otimes t^n)^o = (-1)^h\sum_{m \in {\N}} \frac{1}{m!}
(L(1)^mv) \otimes t^{-n-m-2+2h}\in V \otimes {\C}[t,t^{-
1}],
\end{equation} 
and extend by linearity to $V \otimes {\C}[t,t^{-1}]$. 
That is, $(v \otimes t^n)^o = v^o_{(n)}$ for the special case $Z(\cdot, x)
= Y_t(\cdot, x)$ discussed above. (Note that for general $Z$, 
we cannot expect to
be able to define an analogous $o$-operation on $U$.)  Also consider the map
\begin{eqnarray}
Y^o_t(\cdot,x) = (Y_t(\cdot,x))^o: V & \rightarrow & (V \otimes
{\C}[t,t^{- 1}])[[x,x^{-1}]]\nno\\ v & \mapsto & Y^o_t(v,x) = \sum_{n
\in {\Z}}(v \otimes t^n)^ox^{-n-1}.
\end{eqnarray} 
Then for general
$Z(\cdot, x)$ as above, we can define a linear map
\begin{eqnarray}\label{3.40}
\varepsilon_{Z}: \ V \otimes {\C}[t,t^{-1}]
& \rightarrow & U\nno\\
v \otimes t^n & \mapsto & v_{(n)}
\end{eqnarray}
(``evaluation with respect to $Z$''), i.e.,
\begin{equation}
\varepsilon_{Z}:\ Y_t(v,x) \mapsto Z(v,x),
\end{equation}
 and a linear map
\begin{eqnarray}
\varepsilon^o_{Z}:\ V \otimes {\C}[t,t^{-1}]
& \rightarrow & U\nno\\
v \otimes t^n & \mapsto & v^o_{(n)},
\end{eqnarray}
 i.e., 
\begin{equation}
\varepsilon^o_{Z}:\ Y_t(v,x) \mapsto
Z^o(v,x).
\end{equation}
 Then 
\begin{equation}
\varepsilon^o_{Z} = \varepsilon_{Z} \circ o,
\end{equation}
that is,
\begin{equation}
\varepsilon_{Z}(Y^o_t(v,x)) = Z^o(v,x),
\end{equation}
or equivalently, the diagram
\begin{eqnarray}
Y_t(v,x) & \stackrel{\varepsilon_{Z}}{\longmapsto}
& Z(v,x)\nno\\
{o} \bar{\downarrow}\hspace{1.5em} & & \hspace{1.5em}
\bar{\downarrow} \;(Z(\cdot, x) \mapsto Z^o(\cdot, x))\nno\\
Y^o_t(v,x) &
\stackrel{\varepsilon_{Z}}{\longmapsto} & Z^o(v,x) 
\end{eqnarray} 
commutes. Note that the components $v^o_{(n)}$ of $Z^{o}(v, x)$ depend on
all the components $v_{(n)}$ of $Z(v,x)$ (for
arbitrary $v$), whereas the component $(v \otimes t^n)^o$ of
$Y^o_t(v,x)$ can be defined generically and abstractly; $(v
\otimes t^n)^o$ depends linearly on $v \in V$ alone.
 
     Since in general $Z^{oo}(v, x) = Z(v, x)$, we know that 
\begin{equation}
Y^{oo}_t(v, x) =Y_t(v, x)
\end{equation}
as a special case, and in particular (and equivalently),
\begin{equation}
(v \otimes t^n)^{oo} = v \otimes t^n
\end{equation}
for all $v\in V$ and $n\in {\Z}$. Thus $o$ is an involution of $V
\otimes {\C}[t,t^{-1}]$.
 
     Furthermore, the involution $o$ of $V \otimes {\C}[t,t^{-1}]$
extends canonically to a linear map
$$V \otimes {\C}[[t,t^{-1}]] \stackrel{o}{\rightarrow} V \otimes
{\C}[[t,t^{-1}]].$$
In fact, consider the restriction of $o$ to $V=V \otimes t^0$:
\begin{eqnarray}
V &\stackrel{o}{\rightarrow} &V \otimes
{\C}[t,t^{-1}]\nno\\
v &\mapsto& v^o = (-1)^h \sum_{m \in {\N}}
\frac{1}{m!}(L(1)^m v) \otimes t^{-m-2+2h},
\end{eqnarray}
 extended by
linearity {from} $V_{(h)}$ to $V$.  Then for $v\in V$, we may write 
\begin{equation}\label{vo1}
v^{o}=e^{t^{-1}L(1)}(-t^{2})^{L(0)}vt^{-2}.
\end{equation}
Also, for $v\in V$ and $n\in {\Z}$,
\begin{equation}
(v \otimes t^n)^o = v^ot^{-n},
\end{equation}
 and it is clear that $o$ extends to $V\otimes {\C}[[t, t^{-1}]]$: 
For $f(t) \in {\C}[[t,t^{-1}]],$ 
\begin{equation}
(v \otimes
f(t))^o = v^of(t^{-1}).
\end{equation}

To see that $o$ is an involution of this larger space, first note that 
\begin{equation}
v^{oo} = v
\end{equation}
(although $v^{o}\not\in V$ in general). (This could of course
alternatively be proved by direct calculation using formula (\ref{3.38}).)
Also, for $g(t) \in {\C}[t,t^{-1}]$ and $f(t) \in {\C}[[t,t^{-1}]],$
\begin{equation}
(v \otimes
g(t)f(t))^o = v^og(t^{-1})f(t^{-1})= (v \otimes
g(t))^of(t^{-1}).
\end{equation}
 Thus for all $u \in V\otimes {\C}[t,t^{-
1}]$ and $f(t) \in {\C}[[t,t^{-1}]],$ 
\begin{equation}
(uf(t))^o = u^of(t^{-
1}).
\end{equation}
  It follows that 
\begin{eqnarray}
(v \otimes f(t))^{oo} &
= & (v^of(t^{-1}))^o\nno\\ 
& = & v^{oo}f(t)\nno\\
& = & vf(t)\nno\\ 
& = & v
\otimes f(t),
\end{eqnarray}
and we have shown that $o$ is an involution of $V
\otimes {\C}[[t,t^{-1}]]$.
We have  
\begin{equation}
o: V \otimes {\C}((t)) \leftrightarrow V
\otimes {\C}((t^{-1})).
\end{equation}

 Note that
\begin{eqnarray}\label{op-y-t}
Y^o_t(v,x) & = & \sum_{n \in {\Z}}(v \otimes
t^n)^ox^{-n-1}\nno\\
& = & v^o \sum_{n \in {\Z}}t^{-n}x^{-n-1}\nno\\
&= & v^ox^{-1} \delta(tx)\nno\\
& = & v^ot \delta(tx)\nno\\
& \in & V
\otimes {\C}[[t,t^{-1},x,x^{-1}]].
\end{eqnarray}
Thus the map $v \mapsto Y^o_t(v,x)$ is the linear map given by
multiplying $v^o$ by the ``universal element'' $t \delta(tx)$ 
(cf. the comment following (\ref{3.6})). By (\ref{vo1}), we also have 
\begin{eqnarray}\label{op-y-t-2}
Y^o_t(v,x)&=&e^{t^{-1}L(1)}(-t^{2})^{L(0)}vt^{-1}\delta(tx)\nn
&=&e^{xL(1)}(-x^{-2})^{L(0)}v\otimes x\delta(tx).
\end{eqnarray}
For all $f(x) \in {\C}[[x,x^{-1}]],$ $f(x)Y^o_t(v,x)$
is defined and 
\begin{eqnarray}
f(x)Y^o_t(v,x) & = & f(t^{-
1})Y^o_t(v,x)\nno\\
& = & v^of(t^{-1})t\delta(tx).
\end{eqnarray}

Now we return to the starting point---the original special case: $U =
\mbox{End}\; W$ and $Z(\cdot,x) = Y_W(\cdot,x):\ V \to (\mbox{End}\;
W)[[x,x^{-1}]].$ The corresponding map
\begin{eqnarray}
\varepsilon_{Z} =
\varepsilon_{Y_W}:\ V[t,t^{-1}] & \to & \mbox{End}\; W\nno\\
v \otimes
t^n & \mapsto & v_{(n)}
\end{eqnarray} 
(recall (\ref{3.40})) is just the map
$\tau_W: v \otimes t^n
\mapsto v_n$ (recall (\ref{tauW})), i.e.,
 $v_{(n)} = v_n$ in this case.  Recall that this
map extends canonically to $V \otimes {\C}((t))$.  The map
$\varepsilon^o_{Z}$ (see above) is 
\[
\tau^o_{W}=\tau_{W}\circ o:\ V \otimes {\C}[t,t^{-1}] \to \mbox{End}\; W,
\]
and this map extends canonically to $V \otimes {\C}((t^{-
1}))$. In addition to (\ref{3.7}),  we have 
\begin{equation}\label{tauw-yto}
\tau_W(Y^o_t(v,x))  =  Y^o_W(v,x)
\end{equation} 
($v^o_{(n)} =
v^o_n$ in this case; recall (\ref{yo1})).  In case $f(x) \in {\C}((x^{-1})),$
$$f(x)Y^o_W(v,x) = \tau_W(f(x)Y^o_t(v,x))$$
 is defined and is equal to
$\tau_W(f(t^{-1})Y^o_t(v,z))$ (which is also defined).
 
     The $x$-expansion coefficients of $f(x)Y^o_t(v,x)$, for
$f(x) \in {\C}[[x,x^{-1}]]$, span 
\begin{eqnarray}
v^of(t^{-
1}){\C}[t,t^{-1}] & = & (v{\C}[t,t^{-1}])^of(t^{-1})\nno\\
& = &
(vf(t){\C}[t,t^{-1}])^o.
\end{eqnarray}
The $x$-expansion coefficients of $Y^o_W(v,x)$ span
\begin{eqnarray}
\tau_W(v^o{\C}[t,t^{-1}])&=& \tau_W((v \otimes
{\C}[t,t^{-1}])^o)\nno\\ &=& \tau^o_W(v \otimes {\C}[t,t^{- 1}]).
\end{eqnarray}
In case $f(x) \in {\C}((x^{-1}))$, the $x$-expansion coefficients of
$f(x)Y^o_W(v,x)$ span 
\[
\tau_W(v^of(t^{-1}){\C}[t,t^{-1}])=
\tau^o_W(vf(t){\C}[t,t^{-1}]).
\]
  (Cf. the comments after (\ref{3.9}).)

We shall need spaces of the forms $V\otimes \iota_{+} {\C}[t, t^{-1},
(z+t)^{-1}]$ and 
$V\otimes \iota_{-} {\C}[t, t^{-1},
(z+t)^{-1}]$, where we use the notations
\begin{eqnarray}\label{iota+-}
\iota_+: {\C}(t) & \hookrightarrow & {\C}((t)) \subset
{\C}[[t,t^{-1}]]\nno\\ \iota_-:{\C}(t) & \hookrightarrow &
{\C}((t^{-1})) \subset {\C} [[t,t^{-1}]]
\end{eqnarray}
to denote the operations of expanding a rational function of the
formal variable $t$ in
the indicated directions (as in Section 8.1 of
\cite{FLM2}).  
We shall also 
need certain translation operations, as well as the $o$-operation. 
For $a \in \mathbb{C}$, we define the translation
isomorphism
\begin{eqnarray}
T_a:\ \mathbb{C}(t) & \stackrel{\sim}{\rightarrow}
& \mathbb{C}(t)\nno\\f(t) & \mapsto & f(t+a)
\end{eqnarray} 
and (for our use below) we also set
\begin{equation}
T^\pm_a = \iota_\pm \circ T_a \circ
\iota^{-1}_+:\
\iota_+\mathbb{C}(t) \hookrightarrow \mathbb{C}((t^{\pm1})).
\end{equation}
(Note that the domains of these maps consist of certain series expansions
of formal rational functions rather than of 
formal rational functions themselves.)
The following lemma will be needed for our action $\tau_{P(z)}$ in 
Section 5.2 below (we shall sometimes write $o(Y_{t}(v, x_{1}))$ for 
$Y_{t}^{o}(v, x_{1})$, etc.):

\begin{lemma}\label{tauP}
Let $z\in \C^{\times}$. Then
\begin{eqnarray}
\lefteqn{o\left(x_0^{-1}\delta\left(\frac{x^{-1}_1-z}{x_0}\right)
Y_{t}(v,x_1)\right)=x_0^{-1}\delta\left(\frac{x^{-1}_1-z}{x_0}\right)
Y^o_{t}(v,x_1),}\label{ztr1}\\
\lefteqn{(\iota_+\circ\iota_-^{-1}\circ o)\left(x_0^{-1}\delta\left(
\frac{x^{-1}_1 -z}{x_0}\right)Y_{t}(v,x_1)\right)=x_0^{-1}\delta\left(
\frac{z-x^{-1}_1}{-x_0}\right) Y^o_{t}(v,x_1),}\label{ztr2}\\
\lefteqn{(\iota_+\circ T_{z}\circ\iota_-^{-1}\circ o)\left(x_0^{-1}
\delta\left(\frac{x^{-1}_1-z}{x_0}\right) Y_{t}(v,x_1)\right)}\nno\\
&&\hspace{4.5em}=z^{-1}\delta \left(\frac{x^{-1}_1-x_0}{z}\right)
Y_{t}(e^{x_{1}L(1)}(-x_{1}^{-2})^{L(0)}v,x_0).\quad\quad\quad\quad\quad
\label{ztr3}
\end{eqnarray}
\end{lemma}
\pf 
Formula (\ref{ztr1}) is immediate {}from the definition 
of the map $o$ (recall (\ref{3.38})). 
By (\ref{ztr1}), (\ref{op-y-t}) and (\ref{Xx1x2=Xx2x2}), we have 
\begin{eqnarray}
\lefteqn{(\iota_+\circ\iota_-^{-1}\circ o)\left(x_0^{-1}\delta\left(
\frac{x^{-1}_1 -z}{x_0}\right)Y_{t}(v,x_1)\right)}\nn
&&=(\iota_+\circ\iota_-^{-1})\left(x_0^{-1}\delta\left(
\frac{x^{-1}_1 -z}{x_0}\right)Y^{o}_{t}(v,x_1)\right)\nn
&&=(\iota_+\circ\iota_-^{-1})\left(x_0^{-1}\delta\left(
\frac{x^{-1}_1 -z}{x_0}\right)v^{o}t\delta(tx_{1})\right)\nn
&&=(\iota_+\circ\iota_-^{-1})\left(x_0^{-1}\delta\left(
\frac{t -z}{x_0}\right)v^{o}t\delta(tx_{1})\right)\nn
&&=x_0^{-1}\delta\left(
\frac{z-t}{-x_0}\right)v^{o}t\delta(tx_{1})\nn
&&=x_0^{-1}\delta\left(
\frac{z-x^{-1}_1}{-x_0}\right)v^{o}t\delta(tx_{1})\nn
&&=x_0^{-1}\delta\left(
\frac{z-x^{-1}_1}{-x_0}\right) Y^o_{t}(v,x_1),
\end{eqnarray}
proving (\ref{ztr2}).
For (\ref{ztr3}), note that by (\ref{op-y-t-2}), 
the coefficient of $x_0^{-n-1}$ in the right-hand side of (\ref{ztr1}) is
\begin{eqnarray*}
\lefteqn{(x_1^{-1}-z)^n\left(e^{x_1L(1)}(-x_1^{-2})^{L(0)}v\otimes
x_1\delta\left(\frac t{x_1^{-1}}\right)\right)}\\
&&=(t-z)^n\left(e^{x_1L(1)}(-x_1^{-2})^{L(0)}v\otimes
x_1\delta\left(\frac t{x_1^{-1}}\right)\right).
\end{eqnarray*}
Acted on by $\iota_+\circ T_z\circ\iota_-^{-1}$, this becomes
\begin{eqnarray*}
\lefteqn{t^n\left(e^{x_1L(1)}(-x_1^{-2})^{L(0)}v\otimes
x_1\delta\left(\frac {z+t}{x_1^{-1}}\right)\right)}\\
&&=z^{-1}\delta\left(\frac {x_1^{-1}-t}{z}\right)\left(e^{x_1L(1)}
(-x_1^{-2})^{L(0)}v\otimes t^n\right),
\end{eqnarray*}
which by (\ref{3.5}) is the  coefficient of $x_0^{-n-1}$ in the right-hand side of
(\ref{ztr3}).  
\epfv

We shall be interested in 
\begin{equation}
T^\pm_{-z}: \iota_+\mathbb{C}[t,t^{-1},(z +
t)^{-1}]
\hookrightarrow \mathbb{C}((t^{\pm1})),
\end{equation}
where $z$ is an arbitrary nonzero complex number, as above.
The images of these two maps are $\iota_{\pm}\mathbb{C}[t,t^{-1},(z-t)^{-1}]$.
 
Extend the maps $T^\pm_{-z}$ to linear isomorphisms 
\begin{equation}\label{Tpm-z}
T^\pm_{-z}:\ V \otimes
\iota_+\mathbb{C}[t,t^{-1},(z+t)^{-1}] \stackrel{\sim}{\rightarrow}
V
\otimes \iota_\pm\mathbb{C}[t,t^{-1},(z-t)^{-1}]
\end{equation}
 given by $1\otimes T^\pm_{-z}$ with $T^\pm_{-z}$ as defined above.
Note that the domain of these two maps is described by 
(\ref{3.12})--(\ref{3.13}),
that the image of the map $T^{+}_{-z}$ is described by 
(\ref{3.15})--(\ref{3.16})
and that the image of the map $T^{-}_{-z}$ is described by
(\ref{3.18})--(\ref{3.19}).

We have the two mutually inverse maps
\begin{eqnarray}
V \otimes \iota_-\mathbb{C}[t,t^{-
1},(z-t)^{-1}] & \stackrel{o}{\rightarrow}& 
V \otimes \iota_+\mathbb{C}[t,t^{-1},(z^{-1}-t)^{-1}]\nno\\
v \otimes f(t) &\mapsto &v^of(t^{-1})  
\end{eqnarray}
 and 
\begin{eqnarray}
V \otimes \iota_+{\C}[t,t^{-1},(z^{-1}-t)^{-1}] &
\stackrel{o}{\rightarrow}& V \otimes
\iota_-\mathbb{C}[t,t^{-1},(z-t)^{-1}]\nno\\ v \otimes f(t) & \mapsto
&v^of(t^{-1}),
\end{eqnarray}
which are both isomorphisms. We form the composition
\begin{equation}\label{To-z}
T^o_{-z} = o
\circ T^-_{-z}
\end{equation}
to obtain another isomorphism
\[
T^o_{-z}: V \otimes \iota_+{\C}[t,t^{-1},(z+t)^{-1}]
\stackrel{\sim}{\rightarrow} V \otimes
\iota_+\mathbb{C}[t,t^{-1},(z^{-1}-t)^{-1}].
\]
The maps $T^{+}_{-z}$ and $T^o_{-z}$ will be the main ingredients of
our action $\tau_{Q(z)}$ (see Section 5.3 below). The following
result asserts that $T^{+}_{-z}$, $T^{-}_{-z}$ and $T^o_{-z}$
transform the expression (\ref{3.12}) into (\ref{3.15}), (\ref{3.18})
and the $o$-transform of (\ref{3.18}), respectively:

\begin{lemma}\label{lemma5.2}
We have
\begin{eqnarray}
T_{-z}^{+}\left(z^{-1}\delta\left(\frac{x_{1}-x_{0}}{z}\right)Y_{t}(v,
x_{0})\right)
=x_{0}^{-1}\delta\left(\frac{z-x_{1}}{-x_{0}}\right) Y_{t}(v, x_{1}),
\label{3.71}\\
T_{-z}^{-}\left(z^{-1}\delta\left(\frac{x_{1}-x_{0}}{z}\right)Y_{t}(v,
x_{0})\right)
=x_{0}^{-1}\delta\left(\frac{x_{1}-z}{x_{0}}\right) Y_{t}(v, x_{1}),
\label{3.72}\\
T_{-z}^{o}\left(z^{-1}\delta\left(\frac{x_{1}-x_{0}}{z}\right)Y_{t}(v,
x_{0})\right)
=x_{0}^{-1}\delta\left(\frac{x_{1}-z}{x_{0}}\right) Y^{o}_{t}(v, x_{1}).
\label{3.73}
\end{eqnarray}
\end{lemma}
\pf
We  prove (\ref{3.71}): {From} (\ref{3.12}), the coefficient of 
$x_{0}^{-n-1}x_{1}^{-m-1}$ in the left-hand side of
(\ref{3.71}) is $T_{-z}^{+}(v\otimes (z+t)^{m}t^{n})$. By the definitions,
\begin{eqnarray}
T^{+}_{-z}(v\otimes (z+t)^{m}t^{n})
=v\otimes t^{m}(-(z-t))^{n}.
\end{eqnarray}
On the other hand, the right-hand side of (\ref{3.71}) can be written as
\begin{eqnarray}\label{3.75}
v\otimes
x_{0}^{-1}\delta\left(\frac{z-x_{1}}{-x_{0}}\right)  x_{1}^{-1}
\delta\left(\frac{t}{x_{1}}\right)
=v\otimes x_{0}^{-1}\delta\left(\frac{z-t}{-x_{0}}\right) x_{1}^{-1}
\delta\left(\frac{t}{x_{1}}\right),
\end{eqnarray}
where we have used (\ref{3.5}) and the fundamental property 
(\ref{Xx1x2=Xx2x2}) of the
formal $\delta$-function. The coefficient of
$x_{0}^{-n-1}x_{1}^{-m-1}$ in the right-hand side of (\ref{3.75}) is also
$v\otimes t^{m}(-(z-t))^{n}$,  proving (\ref{3.71}).  Formula
(\ref{3.72}) is proved similarly, and (\ref{3.73}) is obtained 
{from} (\ref{3.72}) by the
application of the map $o$.
\epf

\subsection{Constructions of $P(z)$-tensor products}

We  proceed to the construction of $P(z)$-tensor products.  While
one can certainly consider categories ${\cal C}$ in Remark
\ref{bifunctor} that are not closed under the contragredient functor,
it is most natural to consider such categories ${\cal C}$ that are
indeed closed under this functor (recall Notation \ref{MGM}).  Our
constructions of $P(z)$-tensor products will in fact use the
contragredient functor; the $P(z)$-tensor product of (generalized)
modules $W_1$ and $W_2$ will arise as the contragredient module of a
certain subspace of the vector space dual $(W_1 \otimes W_2)^*$.  We
now present this ``double-dual'' approach to the construction of
$P(z)$-tensor products, generalizing the double-dual approach carried
out in \cite{tensor1}--\cite{tensor3}.  At first, we need not fix any
subcategory ${\cal C}$ of ${\cal M}_{sg}$ or ${\cal GM}_{sg}$.
As usual, we take $z\in \C^{\times}$.

We shall be constructing an action of the space $V\otimes \mathbb{C}[t,
t^{-1}, (z^{-1}-t)^{-1}]$ on the space $(W_1 \otimes W_2)^*$, given
generalized $V$-modules $W_1$ and $W_2$. This action will be based
on the translation operations and on the $o$-operation discussed in the 
preceding section.  More
precisely, it is the space $V\otimes \iota_{+} \mathbb{C}[t, t^{-1},
(z^{-1}-t)^{-1}]$ whose action we shall define.

Let $I$ be a $P(z)$-intertwining map of type ${W_3\choose W_1\, W_2}$,
as in Definition \ref{im:imdef}. Consider the contragredient generalized 
$V$-module $(W_{3}', Y_{3}')$, recall the opposite vertex operator 
(\ref{yo}) and formula (\ref{y'}), and recall why the ingredients of 
formula (\ref{im:def})  are well defined. For 
$v\in V$, $w_{(1)}\in W_{1}$, $w_{(2)}\in W_{2}$ and $w'_{(3)}\in W'_3$, 
applying $w'_{(3)}$ to (\ref{im:def}), replacing $x_1$ by
$x_1^{-1}$ in the resulting formula 
and then replacing $v$ by $e^{x_1L(1)}(-x_1^{-2})^{L(0)}v$,
we get:
\begin{eqnarray}\label{im:def'}
\lefteqn{\left\langle x_0^{-1}\delta\bigg(\frac{x^{-1}_1-z}{x_0}\bigg)
Y_{3}'(v, x_1)w'_{(3)}, I(w_{(1)}\otimes w_{(2)})\right\rangle}\nno\\
&&=\left\langle w'_{(3)},z^{-1}\delta\bigg(\frac{x^{-1}_1-x_0}{z}\bigg)
I(Y_1(e^{x_1L(1)}(-x_1^{-2})^{L(0)}v, x_0)w_{(1)}\otimes
w_{(2)})\right\rangle\nno\\
&&\quad +\left\langle w'_{(3)}, x^{-1}_0\delta\bigg(\frac{z-x^{-1}_1}{-x_0}\bigg)
I(w_{(1)}\otimes Y_2^{o}(v, x_1)w_{(2)})\right\rangle.
\end{eqnarray}
We shall use this to motivate our action.

As we discussed in the preceding section (see (\ref{3.18-1}) and
(\ref{3.19-1})), in the left-hand side of (\ref{im:def'}), the
coefficients of
\begin{equation}\label{deltaY3'}
x_0^{-1}\delta\bigg(\frac{x^{-1}_1-z}{x_0}\bigg) Y_{3}'(v, x_1)
\end{equation}
in powers of $x_0$ and $x_1$, for all $v\in V$, span
\begin{equation}\label{tausubW3'}
\tau_{W'_3}(V\otimes \iota_{+}{\mathbb C}[t,t^{-1},(z^{-1}-t)^{-1}])
\end{equation}
(recall (\ref{tauw}) and (\ref{3.7})).
Let us now define an action of $V\otimes \iota_{+}{\mathbb
C}[t,t^{-1},(z^{-1}-t)^{-1}]$ on $(W_1\otimes W_2)^*$. 

\begin{defi}\label{deftau}{\rm
Define the linear action $\tau_{P(z)}$ of
\[
V \otimes \iota_{+}{\mathbb C}[t,t^{- 1}, (z^{-1}-t)^{-1}]
\]
on $(W_1 \otimes W_2)^*$ by 
\begin{eqnarray}\label{taudef0}
(\tau_{P(z)}(\xi)\lambda)(w_{(1)}\otimes w_{(2)})&=&\lambda
(\tau_{W_{1}}((\iota_+\circ T_z\circ\iota_-^{-1}\circ o)\xi)w_{(1)}
\otimes w_{(2)})\nno\\
&&+\lambda (w_{(1)}\otimes\tau_{W_{2}}((\iota_+\circ
\iota_-^{-1}\circ o)\xi)w_{(2)})
\end{eqnarray}
for $\xi\in V \otimes \iota_{+}{\mathbb C}[t,t^{- 1},
(z^{-1}-t)^{-1}]$, $\lambda\in (W_1\otimes W_2)^*$, $w_{(1)}\in W_1$
and $w_{(2)}\in W_2$. (The fact that the right-hand side is well
defined follows immediately {}from the generating-function reformulation
of (\ref{taudef0}) given in (\ref{taudef}) below.) Denote by
$Y'_{P(z)}$ the action of $V\otimes{\mathbb C}[t,t^{-1}]$ on
$(W_1\otimes W_2)^*$ thus defined, that is,
\begin{equation}\label{y'-p-z}
Y'_{P(z)}(v,x)=\tau_{P(z)}(Y_t(v,x))
\end{equation}
for $v \in V$.
}
\end{defi}

By Lemma \ref{tauP}, (\ref{3.7}) and (\ref{tauw-yto}), we see that
formula (\ref{taudef0}) can be written in terms of
generating functions as
\begin{eqnarray}\label{taudef}
\lefteqn{\bigg(\tau_{P(z)}
\bigg(x_0^{-1}\delta\bigg(\frac{x^{-1}_1-z}{x_0}\bigg)
Y_{t}(v, x_1)\bigg)\lambda\bigg)(w_{(1)}\otimes w_{(2)})}\nno\\
&&=z^{-1}\delta\bigg(\frac{x^{-1}_1-x_0}{z}\bigg)
\lambda(Y_1(e^{x_1L(1)}(-x_1^{-2})^{L(0)}v, x_0)w_{(1)}\otimes w_{(2)})
\nno\\
&&\quad +x^{-1}_0\delta\bigg(\frac{z-x^{-1}_1}{-x_0}\bigg)
\lambda(w_{(1)}\otimes Y_2^{o}(v, x_1)w_{(2)})
\end{eqnarray}
for $v\in V$, $\lambda\in (W_1\otimes W_2)^{*}$, $w_{(1)}\in W_1$,
$w_{(2)}\in W_2$; note that by (\ref{3.18-1})--(\ref{3.19-1}), the
expansion coefficients in $x_{0}$ and $x_{1}$ of the left-hand side
span the space of elements in the left-hand side of (\ref{taudef0}). Compare 
formula (\ref{taudef}) with the motivating formula (\ref{im:def'}).
  The generating function form of the action
$Y'_{P(z)}$ can be obtained by taking $\res_{x_0}$ of both sides of
(\ref{taudef}), that is,
\begin{eqnarray}\label{Y'def}
\lefteqn{(Y'_{P(z)}(v,x_1)\lambda)(w_{(1)}\otimes w_{(2)})=
\lambda(w_{(1)}\otimes Y_2^{o}(v, x_1)w_{(2)})}\nno\\
&&\quad+\res_{x_0}z^{-1}\delta\bigg(\frac{x^{-1}_1-x_0}{z}\bigg)
\lambda(Y_1(e^{x_1L(1)}(-x_1^{-2})^{L(0)}v, x_0)w_{(1)}\otimes
w_{(2)}).
\end{eqnarray}

\begin{rema}\label{I-intw}{\rm
Using the actions $\tau_{W'_3}$ and $\tau_{P(z)}$, we can write
(\ref{im:def'}) as
\[
\left(x_0^{-1}\delta\left(\frac{x^{-1}_1-z}{x_0}\right)
Y_{3}'(v, x_1)w'_{(3)}\right)\circ I=
\tau_{P(z)}\left(x_0^{-1}\delta\left(\frac{x^{-1}_1-z}{x_0}\right)
Y_t(v, x_1)\right)(w'_{(3)}\circ I)
\]
or equivalently, as
\[
\left(\tau_{W'_3}\left(x_0^{-1}\delta\left(\frac{x^{-1}_1-z}{x_0}\right)
Y_t(v, x_1)\right)w'_{(3)}\right)\circ I=
\tau_{P(z)}\left(x_0^{-1}\delta\left(\frac{x^{-1}_1-z}{x_0}\right)
Y_t(v, x_1)\right)(w'_{(3)}\circ I).
\]
}
\end{rema}

In the spirit of the discussion related to Lemma \ref{4.36}, we find
it natural to introduce subspaces of $(W_{1}\otimes W_{2})^{*}$
homogeneous with respect to $\tilde{A}$.  Since $W_{1}$ and $W_{2}$
are $\tilde{A}$-graded, $W_{1}\otimes W_{2}$ also has a natural
$\tilde{A}$-grading---the tensor product grading, and we shall write
$(W_{1}\otimes W_{2})^{(\beta)}$ for the homogeneous subspace of
degree $\beta\in \tilde{A}$ of $W_{1}\otimes W_{2}$. For $\beta\in
\tilde{A}$, define
\begin{equation}\label{W1W2beta}
((W_{1}\otimes W_{2})^{*})^{(\beta)} 
=\{\lambda \in (W_{1}\otimes W_{2})^{*}\;|\; \lambda(\tilde{w})=0 
\;{\rm for}\; \tilde{w}\in (W_{1}\otimes W_{2})^{(\gamma)}
\;{\rm with}\; \gamma\ne -\beta\}
\end{equation}
(cf. (\ref{W'beta}) and note the minus sign).  Of course, the full
space $(W_{1}\otimes W_{2})^{*}$ is not $\tilde{A}$-graded since it is
not a direct sum of subspaces homogeneous with respect to $\tilde{A}$.

The space 
$V \otimes \iota_{+}{\mathbb C}[t,t^{- 1}, (z^{-1}-t)^{-1}]$
also has an $A$-grading, induced {}from the $A$-grading on $V$:
For $\alpha\in A$, 
\begin{equation}
(V\otimes 
\iota_{+}{\mathbb C}[t,t^{- 1}, (z^{-1}-t)^{-1}])^{(\alpha)}
=V^{(\alpha)}\otimes 
\iota_{+}{\mathbb C}[t,t^{- 1}, (z^{-1}-t)^{-1}].
\end{equation}

Using these gradings, we formulate:

\begin{defi}\label{linearactioncompatible}
{\rm We call a linear action $\tau$ of 
$V \otimes \iota_{+}{\mathbb C}[t,t^{- 1}, (z^{-1}-t)^{-1}]$
on $(W_1 \otimes W_2)^*$
{\it $\tilde{A}$-compatible} if 
for $\alpha\in A$, $\beta\in \tilde{A}$,
$\xi\in 
(V \otimes \iota_{+}{\mathbb C}[t,t^{- 1}, (z^{-1}-t)^{-1}])^{(\alpha)}$
and $\lambda\in ((W_{1}\otimes W_{2})^{*})^{(\beta)}$,
\[
\tau(\xi)\lambda\in ((W_{1}\otimes W_{2})^{*})^{(\alpha+\beta)}.
\]}
\end{defi}

{}From (\ref{taudef0}) or (\ref{taudef}), we have:

\begin{propo}\label{tau-a-comp}
The action $\tau_{P(z)}$ is $\tilde{A}$-compatible. \epf
\end{propo}

\begin{rema}
{\rm Notice that Proposition \ref{tau-a-comp} is analogous to the
condition (\ref{m-v_l-A}) in the definition of the notion of
(generalized) module.  We now proceed to establish several more of the
module-action properties for our action $\tau_{P(z)}$ on $(W_1 \otimes
W_2)^*$, in both the conformal and M\"obius cases.  However, while we
will prove the commutator formula for our action (see Proposition
\ref{pz-comm} below), we will {\em not} be able to prove the Jacobi
identity on an element $\lambda \in (W_1 \otimes W_2)^*$ until we
assume the ``$P(z)$-compatibility condition'' for the element
$\lambda$ (see Theorem \ref{comp=>jcb} below).  We shall be
constructing a certain subspace $W_1\hboxtr_{P(z)} W_2$ of $(W_1
\otimes W_2)^*$ which under suitable conditions will be a generalized
$V$-module and whose contragredient module will be
$W_1\boxtimes_{P(z)} W_2$ (see Remark \ref{motivationofbackslash} and
Proposition \ref{tensor1-13.7}), and we shall use the
$P(z)$-compatibility condition to describe this subspace (see Theorem
\ref{characterizationofbackslash}).  }
\end{rema}

We have the following result generalizing Proposition 13.3 in
\cite{tensor3}:

\begin{propo}\label{id-dev}
The action $Y'_{P(z)}$ has the property
\[
Y'_{P(z)}({\bf 1},x)=1,
\]
where $1$ on the right-hand side is the identity map of $(W_1\otimes
W_2)^*$. It also has the $L(-1)$-derivative property
\[
\frac{d}{dx}Y'_{P(z)}(v,x)=Y'_{P(z)}(L(-1)v,x)
\]
for $v\in V$.
\end{propo}
\pf The first statement follows directly {}from the definition. We
prove the $L(-1)$-derivative property. {}From (\ref{Y'def}), we
obtain, using (\ref{yo-l-1}),
\begin{eqnarray}\label{der-1}
\lefteqn{\left(\frac{d}{dx}Y'_{P(z)}(v, x)\lambda\right)
(w_{(1)}\otimes w_{(2)})}\nno\\
&&=\frac{d}{dx}\lambda(w_{(1)}\otimes Y_{2}^{o}(v, x)w_{(2)})\nn
&&\quad +\frac{d}{dx}\res_{x_{0}}z^{-1}\delta\left(\frac{x^{-1}-x_{0}}{z}\right)
\lambda(Y_{1}(e^{xL(1)}(-x^{-2})^{L(0)}v, x_{0})w_{(1)}\otimes w_{(2)})\nno\\
&&=\lambda\left(w_{(1)}\otimes \frac{d}{dx}Y_{2}^{o}(v, x)w_{(2)}\right)\nn
&&\quad+ \res_{x_{0}}\left(\frac{d}{dx}\left(z^{-1}
\delta\left(\frac{x^{-1}-x_{0}}{z}\right)\right)\right)
\lambda(Y_{1}(e^{xL(1)}(-x^{-2})^{L(0)}v, x_{0})w_{(1)}\otimes w_{(2)})\nno\\
&&\quad +\res_{x_{0}}z^{-1}\delta\left(\frac{x^{-1}-x_{0}}{z}\right)
\frac{d}{dx}\lambda(Y_{1}(e^{xL(1)}(-x^{-2})^{L(0)}v, x_{0})w_{(1)}
\otimes w_{(2)})
\nonumber\\
&&=\lambda(w_{(1)}\otimes Y_{2}^{o}(L(-1)v, x)w_{(2)})\nn
&&\quad+ \res_{x_{0}}\left(\frac{d}{dx}\left(z^{-1}
\delta\left(\frac{x^{-1}-x_{0}}{z}\right)\right)\right)
\lambda(Y_{1}(e^{xL(1)}(-x^{-2})^{L(0)}v, x_{0})w_{(1)}\otimes w_{(2)})\nno\\
&&\quad +\res_{x_{0}}z^{-1}\delta\left(\frac{x^{-1}-x_{0}}{z}\right)
\lambda(Y_{1}(e^{xL(1)}L(1)(-x^{-2})^{L(0)}v, x_{0})w_{(1)}\otimes w_{(2)})
\nonumber\\
&&\quad -2\res_{x_{0}}z^{-1}\delta\left(\frac{x^{-1}-x_{0}}{z}\right)
\lambda(Y_{1}(e^{xL(1)}L(0)x^{-1}(-x^{-2})^{L(0)}v, x_{0})w_{(1)}
\otimes w_{(2)}).
\end{eqnarray}
The second term
on the right-hand side of (\ref{der-1}) is equal to
\begin{eqnarray}\label{der-2}
\lefteqn{-\res_{x_{0}}x^{-2}\left(\frac{d}{dx^{-1}}\left(z^{-1}
\delta\left(\frac{x^{-1}-x_{0}}{z}\right)\right)\right)\cdot}\nno\\
&&\quad \quad\quad\quad \cdot
\lambda(Y_{1}(e^{xL(1)}(-x^{-2})^{L(0)}v, x_{0})w_{(1)}\otimes w_{(2)})\nno\\
&&=\res_{x_{0}}x^{-2}\left(\frac{d}{dx_{0}}\left(z^{-1}
\delta\left(\frac{x^{-1}-x_{0}}{z}\right)\right)\right)\cdot\nno\\
&&\quad \quad\quad\quad \cdot
\lambda(Y_{1}(e^{xL(1)}(-x^{-2})^{L(0)}v, x_{0})w_{(1)}\otimes w_{(2)})\nno\\
&&=-\res_{x_{0}}x^{-2}z^{-1}
\delta\left(\frac{x^{-1}-x_{0}}{z}\right)\cdot\nno\\
&&\quad \quad\quad\quad \cdot
\frac{d}{dx_{0}}\lambda(Y_{1}(e^{xL(1)}(-x^{-2})^{L(0)}v, x_{0})w_{(1)}
\otimes w_{(2)})\nno\\
&&=-\res_{x_{0}}x^{-2}z^{-1}
\delta\left(\frac{x^{-1}-x_{0}}{z}\right)\cdot\nno\\
&&\quad \quad\quad\quad \cdot
\lambda(Y_{1}(L(-1)e^{xL(1)}(-x^{-2})^{L(0)}v, x_{0})w_{(1)}
\otimes w_{(2)}).
\end{eqnarray}
By (\ref{der-2}), (\ref{log:SL2-3}) and (\ref{log:xLx^}),
the right-hand side of (\ref{der-1}) is equal to
\begin{eqnarray*}
\lefteqn{\displaystyle \lambda(w_{(1)}\otimes 
Y_{2}^{o}(L(-1)v, x)w_{(2)})}\nno\\
&&\quad + \res_{x_{0}}z^{-1}\delta\left(\frac{x^{-1}-x_{0}}{z}\right)
\lambda(Y_{1}(e^{xL(1)}(-x^{-2})^{L(0)}L(-1)v, x_{0})w_{(1)}\otimes w_{(2)})
\nonumber\\
&& =(Y'_{P(z)}(L(-1)v, x)\lambda)(w_{(1)}\otimes w_{(2)}),
\end{eqnarray*}
proving the $L(-1)$-derivative property.
\epfv

\begin{propo}\label{pz-comm}
The action $Y'_{P(z)}$ satisfies the commutator formula for vertex
operators, that is, on $(W_1\otimes W_2)^*$,
\begin{eqnarray*}
\lefteqn{[Y'_{P(z)}(v_1,x_1),Y'_{P(z)}(v_2,x_2)]}\\
&&=\res_{x_0}x_2^{-1}\delta\bigg(\frac{x_1-x_0}{x_2}\bigg)
Y'_{P(z)}(Y(v_1,x_0)v_2,x_2)
\end{eqnarray*}
for $v_1,v_2\in V$.
\end{propo}
\pf
In the following proof, the reader should note the
well-definedness of each expression and the justifiability of each use of
a $\delta$-function property.

Let $\lambda
\in (W_{1}\otimes W_{2})^{*}$, $v_{1}, v_{2}\in V$, $w_{(1)}\in W_{1}$
and $w_{(2)}\in W_{2}$. By (\ref{Y'def}),
\bea\label{y-12}
\lefteqn{(Y'_{P(z)}(v_{1}, x_{1})Y'_{P(z)}(v_{2}, x_{2})
\lambda)(w_{(1)}\otimes w_{(2)})}\nno\\
&&=(Y'_{P(z)}(v_{2}, x_{2})\lambda)(w_{(1)}\otimes Y_2^{o}(v_{1}, x_1)w_{(2)})\nno\\
&&\quad+\res_{y_{1}}z^{-1}\delta\bigg(\frac{x^{-1}_1-y_{1}}{z}\bigg)
(Y'_{P(z)}(v_{2}, x_{2})\lambda)(Y_1(e^{x_1L(1)}(-x_1^{-2})^{L(0)}v_{1}, y_{1})w_{(1)}
\otimes w_{(2)})\nn
&&=\lambda(w_{(1)}\otimes Y_2^{o}(v_{2}, x_2)Y_2^{o}(v_{1}, x_1)w_{(2)})\nno\\
&&\quad+\res_{y_{2}}z^{-1}\delta\bigg(\frac{x^{-1}_{2}-y_{2}}{z}\bigg)
\lambda(Y_1(e^{x_{2}L(1)}(-x_{2}^{-2})^{L(0)}v_{2}, y_{2})w_{(1)}\otimes
Y_2^{o}(v_{1}, x_1)w_{(2)})\nn
&&\quad+\res_{y_{1}}z^{-1}\delta\bigg(\frac{x^{-1}_1-y_{1}}{z}\bigg)
\lambda(Y_1(e^{x_1L(1)}(-x_1^{-2})^{L(0)}v_{1}, y_{1})w_{(1)}\otimes
Y_2^{o}(v_{2}, x_2)w_{(2)})\nn
&&\quad+\res_{y_{1}}\res_{y_{2}}z^{-1}\delta\bigg(\frac{x^{-1}_{1}-y_{1}}{z}\bigg)
z^{-1}\delta\bigg(\frac{x^{-1}_2-y_{2}}{z}\bigg)\cdot\nn
&&\quad\quad\quad\quad\cdot\lambda(Y_1(e^{x_2L(1)}(-x_2^{-2})^{L(0)}v_{2}, y_{2})
Y_1(e^{x_1 L(1)}(-x_1^{-2})^{L(0)}v_{1}, y_{1})w_{(1)}\otimes
w_{(2)}).
\eea
Transposing the subscripts $1$ and $2$ of the symbols $v$, $x$ and $y$,
we also have
\bea\label{y-21}
\lefteqn{(Y'_{P(z)}(v_{2}, x_{2})Y'_{P(z)}(v_{1}, x_{1})
\lambda)(w_{(1)}\otimes w_{(2)})}\nno\\
&&=\lambda(w_{(1)}\otimes Y_2^{o}(v_{1}, x_1)Y_2^{o}(v_{2}, x_2)w_{(2)})\nno\\
&&\quad+\res_{y_{1}}z^{-1}\delta\bigg(\frac{x^{-1}_{1}-y_{1}}{z}\bigg)
\lambda(Y_1(e^{x_{1}L(1)}(-x_{1}^{-2})^{L(0)}v_{1}, y_{1})w_{(1)}\otimes
Y_2^{o}(v_{2}, x_2)w_{(2)})\nn
&&\quad+\res_{y_{2}}z^{-1}\delta\bigg(\frac{x^{-1}_2-y_{2}}{z}\bigg)
\lambda(Y_1(e^{x_2L(1)}(-x_2^{-2})^{L(0)}v_{2}, y_{2})w_{(1)}\otimes
Y_2^{o}(v_{1}, x_1)w_{(2)})\nn
&&\quad+\res_{y_{2}}\res_{y_{1}}z^{-1}\delta\bigg(\frac{x^{-1}_{2}-y_{2}}{z}\bigg)
z^{-1}\delta\bigg(\frac{x^{-1}_1-y_{1}}{z}\bigg)\cdot\nn
&&\quad\quad\quad\quad\cdot\lambda(Y_1(e^{x_1L(1)}(-x_1^{-2})^{L(0)}v_{1}, y_{1})
Y_1(e^{x_2 L(1)}(-x_2^{-2})^{L(0)}v_{2}, y_{2})w_{(1)}\otimes
w_{(2)}).
\eea
The equalities (\ref{y-12}) and (\ref{y-21}) give
\bea\label{y-bracket}
\lefteqn{([Y'_{P(z)}(v_{1}, x_{1}), Y'_{P(z)}(v_{2}, x_{2})]
\lambda)(w_{(1)}\otimes w_{(2)})}\nno\\
&&=\lambda(w_{(1)}\otimes [Y_2^{o}(v_{2}, x_2), Y_2^{o}(v_{1}, x_1)]w_{(2)})\nn
&&\quad-\res_{y_{1}}\res_{y_{2}}z^{-1}\delta\bigg(\frac{x^{-1}_{1}-y_{1}}{z}\bigg)
z^{-1}\delta\bigg(\frac{x^{-1}_2-y_{2}}{z}\bigg)\cdot\nn
&&\quad\quad\quad\quad\cdot\lambda([Y_1(e^{x_1 L(1)}(-x_1^{-2})^{L(0)}v_{1}, y_{1}),
Y_1(e^{x_2L(1)}(-x_2^{-2})^{L(0)}v_{2}, y_{2})]w_{(1)}\otimes
w_{(2)})\nn
&&=\res_{x_{0}}x_{2}^{-1}
\delta\left(\frac{x_{1}-x_{0}}{x_{2}}\right)\lambda(w_{(1)}\otimes
Y_2^{o}(Y(v_{1}, x_0)v_{2}, x_2)w_{(2)})\nn
&&\quad-\res_{y_{1}}\res_{y_{2}}\res_{x_{0}}
z^{-1}\delta\bigg(\frac{x^{-1}_{1}-y_{1}}{z}\bigg)
z^{-1}\delta\bigg(\frac{x^{-1}_2-y_{2}}{z}\bigg)
y_{2}^{-1}\delta\left(\frac{y_{1}-x_{0}}{y_{2}}\right)\cdot\nn
&&\quad\quad\quad\quad\cdot\lambda(
Y_{1}(Y(e^{x_1 L(1)}(-x_1^{-2})^{L(0)}v_{1}, x_{0})
e^{x_2L(1)}(-x_2^{-2})^{L(0)}v_{2}, y_{2})w_{(1)}\otimes
w_{(2)})
\eea
(recall (\ref{op-jac-id})).

But we have
\bea\label{delta-idty}
\lefteqn{z^{-1}\delta\bigg(\frac{x^{-1}_{1}-y_{1}}{z}\bigg)
z^{-1}\delta\bigg(\frac{x^{-1}_2-y_{2}}{z}\bigg)
y_{2}^{-1}\delta\left(\frac{y_{1}-x_{0}}{y_{2}}\right)}\nno\\
&&=\left({\displaystyle \sum_{m,n\in {\mathbb Z}}}
\frac{(x_{1}^{-1}-y_{1})^{m}}{z^{m+1}}
\frac{(x_{2}^{-1}-y_{2})^{n}}{z^{n+1}}\right)
y_{2}^{-1}\delta\left(\frac{y_{1}-x_{0}}{y_{2}}\right)\nno\\
&&=\left({\displaystyle \sum_{m,n\in {\mathbb Z}}}(x_{2}^{-1}-y_{2})^{-1}
\left(\frac{x_{1}^{-1}-y_{1}}{x_{2}^{-1}-y_{2}}\right)^{m}
\frac{(x_{2}^{-1}-y_{2})^{m+n+1}}
{z^{m+n+2}} \right)
y_{2}^{-1}\delta\left(\frac{y_{1}-x_{0}}{y_{2}}\right)\nno\\
&&=\left({\displaystyle \sum_{m,k\in {\mathbb Z}}}(x_{2}^{-1}-y_{2})^{-1}
\left(\frac{x_{1}^{-1}-y_{1}}{x_{2}^{-1}-y_{2}}\right)^{m}
z^{-1}\left(\frac{x_{2}^{-1}-y_{2}}{z}\right)^{k}\right)
y_{2}^{-1}\delta\left(\frac{y_{1}-x_{0}}{y_{2}}\right)\nno\\
&&=(x_{2}^{-1}-y_{2})^{-1}\delta\left(\frac{x_{1}^{-1}-y_{1}}{x_{2}^{-1}-y_{2}}\right)
z^{-1}\delta\left(\frac{x_{2}^{-1}-y_{2}}{z}\right)
y_{2}^{-1}\delta\left(\frac{y_{1}-x_{0}}
{y_{2}}\right)\nno\\
&&=x_{2}\delta\left(\frac{x_{1}^{-1}-(y_{1}-y_{2})}{x_{2}^{-1}}\right)
z^{-1}\delta\left(\frac{x_{2}^{-1}-y_{2}}{z}\right)
y_{1}^{-1}\delta\left(\frac{y_{2}+x_{0}}
{y_{1}}\right)\nno\\
&&=x_{2}\delta\left(\frac{x_{1}^{-1}-x_{0}}{x_{2}^{-1}}\right)
z^{-1}\delta\left(\frac{x_{2}^{-1}-y_{2}}{z}\right)y_{1}^{-1}\delta
\left(\frac{y_{2}+x_{0}}{y_{1}}\right).
\eea
By (\ref{log:p2}), (\ref{log:p3}) and (\ref{xe^Lx}), we also have
\bea\label{sl2-idty}
\lefteqn{Y(e^{x_1 L(1)}(-x_1^{-2})^{L(0)}v_{1}, x_{0})
e^{x_2L(1)}(-x_2^{-2})^{L(0)}}\nn
&&=e^{x_2L(1)}Y\left(e^{-x_{2}(1+x_{0}x_{2})L(1)}(1+x_{0}x_{2})^{-2L(0)}
e^{x_1 L(1)}(-x_1^{-2})^{L(0)}v_{1}, \frac{x_{0}}{1+x_{0}x_{2}}\right)
(-x_2^{-2})^{L(0)}\nn
&&=e^{x_2L(1)}(-x_2^{-2})^{L(0)}
Y\Biggl((-x_2^{2})^{L(0)}e^{-x_{2}(1+x_{0}x_{2})L(1)}\cdot\nn
&&\quad\quad\quad\quad\quad\quad\quad\quad\quad\quad\quad\cdot (1+x_{0}x_{2})^{-2L(0)}
e^{x_1 L(1)}(-x_1^{-2})^{L(0)}v_{1}, -\frac{x_{0}x_{2}^{2}}{1+x_{0}x_{2}}\Biggr)
\nn
&&=e^{x_2L(1)}(-x_2^{-2})^{L(0)}
Y\Biggl(e^{-x_{2}(1+x_{0}x_{2})(-x_2^{-2})L(1)}
(-x_2^{2})^{L(0)}\cdot\nn
&&\quad\quad\quad\quad\quad\quad\quad\quad\quad\quad\quad\cdot(1+x_{0}x_{2})^{-2L(0)}
e^{x_1 L(1)}(-x_1^{-2})^{L(0)}v_{1}, -\frac{x_{0}x_{2}^{2}}{1+x_{0}x_{2}}\Biggr)
\nn
&&=e^{x_2L(1)}(-x_2^{-2})^{L(0)}
Y\Biggl(e^{(x_{2}^{-1}+x_{0})L(1)}
(-(x_{2}^{-1}+x_{0})^{2})^{-L(0)}
e^{x_1 L(1)}(-x_1^{-2})^{L(0)}v_{1}, -\frac{x_{0}x_{2}}{x_{2}^{-1}+x_{0}}\Biggr)
\nn
&&=e^{x_2L(1)}(-x_2^{-2})^{L(0)}
Y\Biggl(e^{(x_{2}^{-1}+x_{0})L(1)}e^{-x_1 (x_{2}^{-1}+x_{0})^{2}L(1)}
\cdot\nn
&&\quad\quad\quad\quad\quad\quad\quad\quad\quad\quad\quad\cdot
(-(x_{2}^{-1}+x_{0})^{2})^{-L(0)}(-x_1^{-2})^{L(0)}v_{1},
-\frac{x_{0}x_{2}}{x_{2}^{-1}+x_{0}}\Biggr)
\nn
&&=e^{x_2L(1)}(-x_2^{-2})^{L(0)}
Y\left(e^{(x_{2}^{-1}+x_{0})L(1)}
e^{-x_1 (x_{2}^{-1}+x_{0})^{2}L(1)}
((x_{2}^{-1}+x_{0})x_{1})^{-2L(0)}v_{1},
-\frac{x_{0}x_{2}}{x_{2}^{-1}+x_{0}}\right).\nn
&&
\eea

Using (\ref{delta-idty}), (\ref{sl2-idty}) and the basic
properties of the formal delta function, we see that
(\ref{y-bracket}) becomes
\bea
\lefteqn{([Y'_{P(z)}(v_{1}, x_{1}), Y'_{P(z)}(v_{2}, x_{2})]\lambda)(w_{(1)}
\otimes w_{(2)})}\nno\\
&&=\res_{x_{0}}x_{2}^{-1}
\delta\left(\frac{x_{1}-x_{0}}{x_{2}}\right)\lambda(w_{(1)}\otimes
Y_2^{o}(Y(v_{1}, x_0)v_{2}, x_2)w_{(2)})\nn
&&\quad-\res_{y_{1}}\res_{y_{2}}\res_{x_{0}}
x_{2}\delta\left(\frac{x_{1}^{-1}-x_{0}}{x_{2}^{-1}}\right)
z^{-1}\delta\left(\frac{x_{2}^{-1}-y_{2}}{z}\right)y_{1}^{-1}\delta
\left(\frac{y_{2}+x_{0}}{y_{1}}\right)\cdot\nn
&&\quad\quad\quad\quad\cdot\lambda\Biggl(
Y_{1}\Biggl(e^{x_2L(1)}(-x_2^{-2})^{L(0)}
Y\Biggl(e^{(x_{2}^{-1}+x_{0})L(1)}
e^{-x_1 (x_{2}^{-1}+x_{0})^{2}L(1)}\cdot\nn
&&\quad\quad\quad\quad\quad\quad\quad\quad\quad\quad\quad\cdot
((x_{2}^{-1}+x_{0})x_{1})^{-2L(0)}v_{1},
-\frac{x_{0}x_{2}}{x_{2}^{-1}+x_{0}}\Biggr)v_{2}, y_{2}\Biggr)w_{(1)}\otimes
w_{(2)}\Biggr)\nn
&&=\res_{x_{0}}x_{2}^{-1}
\delta\left(\frac{x_{1}-x_{0}}{x_{2}}\right)\lambda(w_{(1)}\otimes
Y_2^{o}(Y(v_{1}, x_0)v_{2}, x_2)w_{(2)})\nn
&&\quad-\res_{x_{0}}\res_{y_{2}}
x_{2}\delta\left(\frac{x_{1}^{-1}-x_{0}}{x_{2}^{-1}}\right)
z^{-1}\delta\left(\frac{x_{2}^{-1}-y_{2}}{z}\right)\cdot\nn
&&\quad\quad\quad\quad\cdot\lambda(
Y_{1}(e^{x_2L(1)}(-x_2^{-2})^{L(0)}
Y(e^{x_{1}^{-1}L(1)}
e^{-x_1^{-1}L(1)}\cdot\nn
&&\quad\quad\quad\quad\quad\quad\quad\quad\quad\quad\quad\cdot
(x_{1}^{-1}x_{1})^{-2L(0)}v_{1},
-x_{0}x_{1}x_{2})v_{2}, y_{2})w_{(1)}\otimes
w_{(2)})\nn
&&=\res_{x_{0}}x_{2}^{-1}
\delta\left(\frac{x_{1}-x_{0}}{x_{2}}\right)\lambda(w_{(1)}\otimes
Y_2^{o}(Y(v_{1}, x_0)v_{2}, x_2)w_{(2)})\nn
&&\quad-\res_{x_{0}}x_{2}\delta\left(\frac{x_{2}+(-x_{0}x_{1}x_{2})}{x_{1}}\right)
\res_{y_{2}}
z^{-1}\delta\left(\frac{x_{2}^{-1}-y_{2}}{z}\right)\cdot\nn
&&\quad\quad\quad\quad\cdot\lambda(
Y_{1}(e^{x_2L(1)}(-x_2^{-2})^{L(0)}
Y(v_{1},
-x_{0}x_{1}x_{2})v_{2}, y_{2})w_{(1)}\otimes
w_{(2)})\nn
&&=\res_{x_{0}}x_{2}^{-1}
\delta\left(\frac{x_{1}-x_{0}}{x_{2}}\right)\lambda(w_{(1)}\otimes
Y_2^{o}(Y(v_{1}, x_0)v_{2}, x_2)w_{(2)})\nn
&&\quad+\res_{y_{0}}x_{1}^{-1}\delta\left(\frac{x_{2}+y_{0}}{x_{1}}\right)
\res_{y_{2}}
z^{-1}\delta\left(\frac{x_{2}^{-1}-y_{2}}{z}\right)\cdot\nn
&&\quad\quad\quad\quad\cdot\lambda(
Y_{1}(e^{x_2L(1)}(-x_2^{-2})^{L(0)}
Y(v_{1},
y_{0})v_{2}, y_{2})w_{(1)}\otimes
w_{(2)})\nn
&&=\res_{x_{0}}x_{2}^{-1}
\delta\left(\frac{x_{1}-x_{0}}{x_{2}}\right)\lambda(w_{(1)}\otimes
Y_2^{o}(Y(v_{1}, x_0)v_{2}, x_2)w_{(2)})\nn
&&\quad+\res_{x_{0}}x_{2}^{-1}\delta\left(\frac{x_{1}-x_{0}}{x_{2}}\right)
\res_{y_{2}}
z^{-1}\delta\left(\frac{x_{2}^{-1}-y_{2}}{z}\right)\cdot\nn
&&\quad\quad\quad\quad\cdot\lambda(
Y_{1}(e^{x_2L(1)}(-x_2^{-2})^{L(0)}
Y(v_{1},
x_{0})v_{2}, y_{2})w_{(1)}\otimes
w_{(2)})\nn
&&=\res_{x_{0}}x_{2}^{-1}
\delta\left(\frac{x_{1}-x_{0}}{x_{2}}\right)
(Y'_{P(z)}(Y(v_{1}, x_{0})v_{2}, x_{2})\lambda)(w_{(1)}\otimes w_{(2)}).
\eea
Since $\lambda$, $w_{(1)}$ and $w_{(2)}$ are arbitrary,
this  equality gives the commutator formula for $Y'_{P(z)}$.
\epfv

The following observations are analogous to those in Remark 8.1 of
\cite{tensor2} (concerning the case of $Q(z)$ rather than $P(z)$):

\begin{rema}
{\rm The proof of Proposition \ref{pz-comm} suggests the following:  
Using the definitions (\ref{taudef0}) and (\ref{taudef})
as motivation, we define 
a (linear) action
$\sigma_{P(z)}$ of $V \otimes \iota_{+}{\mathbb C}[t,t^{- 1}, 
(z^{-1}-t)^{-1}]$
on the 
vector space $W_{1}\otimes W_{2}$ (as opposed to $(W_{1}\otimes W_{2})^{*}$)
as follows:
\begin{equation}
\sigma_{P(z)}(\xi)(w_{(1)}\otimes w_{(2)})
=
\tau_{W_{1}}((\iota_+\circ T_z\circ\iota_-^{-1}\circ o)\xi)w_{(1)}
\otimes w_{(2)}
+w_{(1)}\otimes\tau_{W_{2}}((\iota_+\circ
\iota_-^{-1}\circ o)\xi)w_{(2)}
\end{equation}
for $\xi\in V \otimes \iota_{+}{\mathbb C}[t,t^{- 1}, 
(z^{-1}-t)^{-1}]$, $w_{(1)}\in
W_{1}$, $w_{(2)}\in W_{2}$, or equivalently,
\begin{eqnarray}\label{sigma-p-z}
\lefteqn{\bigg(\sigma_{P(z)}\bigg(x_0^{-1}
\delta\bigg(\frac{x^{-1}_1-z}{x_0}\bigg)
Y_{t}(v, x_1)\bigg)\bigg)
(w_{(1)}\otimes w_{(2)})}\nno\\
&&=z^{-1}\delta\bigg(\frac{x^{-1}_1-x_0}{z}\bigg)
Y_1(e^{x_1L(1)}(-x_1^{-2})^{L(0)}v, x_0)w_{(1)}\otimes w_{(2)}
\nno\\
&&\quad +x^{-1}_0\delta\bigg(\frac{z-x^{-1}_1}{-x_0}\bigg)
w_{(1)}\otimes Y_2^{o}(v, x_1)w_{(2)}.
\end{eqnarray}
That is, the operators $\sigma_{P(z)}(\xi)$ and $\tau_{P(z)}(\xi)$ are 
mutually adjoint:
\begin{equation}
(\tau_{P(z)}(\xi)\lambda)(w_{(1)}\otimes w_{(2)})
=\lambda(\sigma_{P(z)}(\xi)(w_{(1)}\otimes w_{(2)})).
\end{equation}
While this action on $W_{1}\otimes W_{2}$ is not very useful, it has the 
following three properties:
\begin{equation}\label{sigma-id}
\sigma_{P(z)}(Y_{t}({\bf 1}, x))=1,
\end{equation}
\begin{equation}\label{sigma-dev}
\frac{d}{dx}\sigma_{P(z)}(Y_{t}(v, x))=\sigma_{P(z)}(Y_{t}(L(-1)v, x))
\end{equation}
for $v\in V$,
\begin{eqnarray}\label{sigma-comm}
\lefteqn{[\sigma_{P(z)}(Y_{t}(v_{2}, x_{2})), \sigma_{P(z)}(Y_{t}(v_{1}, 
x_{1}))]}\nno\\
&&=\res_{x_{0}}x_{2}^{-1}\delta\left(\frac{x_{1}-x_{0}}{x_{2}}\right)
\sigma_{P(z)}(Y_{t}(Y(v_{1}, x_{0})v_{2}, x_{2}))
\end{eqnarray}
for $v_{1}, v_{2}\in V$ (the opposite commutator formula). These
follow either {}from the assertions of Propositions \ref{id-dev} and 
\ref{pz-comm} or,
better, {}from the fact that it was actually 
(\ref{sigma-id})--(\ref{sigma-comm}) that the
proofs of these propositions were proving.}
\end{rema}

\begin{rema}
{\rm Taking $\res_{x_{0}}$ of (\ref{sigma-p-z}),
we obtain
\begin{eqnarray}\label{sigma-p-z-1}
\lefteqn{(\sigma_{P(z)}(Y_{t}(v, x_1)))
(w_{(1)}\otimes w_{(2)})}\nno\\
&&=\res_{x_{0}}z^{-1}\delta\bigg(\frac{x^{-1}_1-x_0}{z}\bigg)
Y_1(e^{x_1L(1)}(-x_1^{-2})^{L(0)}v, x_0)w_{(1)}\otimes w_{(2)}
\nno\\
&&\quad +
w_{(1)}\otimes Y_2^{o}(v, x_1)w_{(2)}.
\end{eqnarray}
Substituting first $(-x_1^{-2})^{-L(0)}e^{-x_1L(1)}v$ for $v$ in
(\ref{sigma-p-z-1}) and then $x^{-1}_{1}$ for $x_{1}$ in the same
formula and using (\ref{xe^Lx}), (\ref{3.5}) and (\ref{op-y-t-2}), we
obtain
\begin{eqnarray}\label{sigma-p-z-1.5}
(\sigma_{P(z)}(Y^{o}_{t}(v, x_1)))
(w_{(1)}\otimes w_{(2)})
&=&\res_{x_{0}}z^{-1}\delta\bigg(\frac{x_1-x_0}{z}\bigg)
Y_1(v, x_0)w_{(1)}\otimes w_{(2)}
\nno\\
&& +
w_{(1)}\otimes Y_2(v, x_1)w_{(2)}.
\end{eqnarray}
Using this, we see that $\sigma_{P(z)}$ can actually be viewed as a
map {}from $V[t, t^{-1}]$ to $V((t))\otimes V[t, t^{-1}]$ defined
by
\begin{eqnarray}\label{sigma-p-z-2}
\sigma_{P(z)}(Y^{o}_{t}(v, x_1))
&=&\res_{x_{0}}z^{-1}\delta\bigg(\frac{x_1-x_0}{z}\bigg)
Y_t(v, x_0)\otimes {\mathbf{1}}
\nno\\
&& +
{\mathbf{1}}\otimes Y_t(v, x_1).
\end{eqnarray}
Let $\Delta_{P(z)}=\sigma_{P(z)}\circ o$.  Then (\ref{sigma-p-z-2})
becomes
\begin{eqnarray}\label{sigma-p-z-3}
\Delta_{P(z)}(Y_{t}(v, x_1))
&=&\res_{x_{0}}z^{-1}\delta\bigg(\frac{x_1-x_0}{z}\bigg)
Y_t(v, x_0)\otimes {\mathbf{1}}
\nno\\
&& +
{\mathbf{1}}\otimes Y_t(v, x_1),
\end{eqnarray}
which again can be viewed as a map
\begin{equation}
\Delta_{P(z)}:V[t, t^{-1}] \rightarrow V((t))\otimes V[t, t^{-1}].
\end{equation}
In formula (2.4) of \cite{MS}, Moore and Seiberg introduced a map
$\Delta_{z, 0}$ which in fact corresponds exactly to the map
$\Delta_{P(z)}$ defined by (\ref{sigma-p-z-3}).  They proposed to
define a $V$-module structure (called ``a representation of
$\mathcal{A}$'' in \cite{MS}, where $\mathcal{A}$ corresponds to our
vertex algebra $V$) on $W_{1}\otimes W_{2}$ by using this map, which
can be viewed as a sort of analogue of a coproduct, but they
acknowledged that E. Witten pointed out ``subtleties in this
definition which are related to the fact that [a] representation of
$\mathcal{A}$ obtained this way is not always a highest weight
representation.''  In fact, it is these subtleties that make it
impossible to work with $W_{1}\otimes W_{2}$; in virtually all
interesting cases, $W_{1}\otimes W_{2}$ does not have a natural
(generalized) $V$-module structure.  This is exactly the reason why we
had to use a completely different approach to construct our
$P(z)$-tensor product of $W_{1}$ and $W_{2}$.}
\end{rema}

When $V$ is in fact a conformal (rather than M\"{o}bius) vertex
algebra, we will write
\begin{equation}\label{13.11}
Y'_{P(z)}(\omega,x)=\sum_{n\in {\mathbb Z}} L'_{P(z)}(n)x^{-n-2}.
\end{equation}
Then {}from the last two propositions we see that the coefficient operators
of $Y'_{P(z)}(\omega, x)$ satisfy the Virasoro algebra commutator
relations, that is,
\[
[L'_{P(z)}(m), L'_{P(z)}(n)]
=(m-n)L'_{P(z)}(m+n)+{\displaystyle\frac1{12}}
(m^3-m)\delta_{m+n,0}c,
\]
with $c\in \C$ the central charge of $V$ (recall Definition 
\ref{cva}).
Moreover, in this case, by setting $v=\omega$ in (\ref{Y'def}) and
taking the coefficient of $x_1^{-j-2}$ for $j=-1, 0, 1$, we find that
\begin{eqnarray}\label{LP'(j)}
\lefteqn{(L'_{P(z)}(j)\lambda)(w_{(1)}\otimes w_{(2)})}\nn
&&=\lambda\Biggl(w_{(1)}\otimes L(-j)w_{(2)}+\Biggl(\sum_{i=0}
^{1-j}{1-j\choose i}z^iL(-j-i)\Biggr)w_{(1)}\otimes w_{(2)}
\Biggr),
\end{eqnarray}
by (\ref{Yoppositeomega}). If $V$ is just a M\"obius vertex algebra, we
define the actions $L'_{P(z)}(j)$ on $(W_1\otimes W_2)^*$ by
(\ref{LP'(j)}) for $j=-1, 0$ and $1$. 

\begin{rema}\label{I-intw2}{\rm
In view of the action $L'_{P(z)}(j)$, the ${\mathfrak s}{\mathfrak
l}(2)$-bracket relations (\ref{im:Lj}) for a $P(z)$-intertwining map
$I$, with notation as in Definition \ref{im:imdef}, can be written as
\begin{equation}\label{I-intw2f}
(L'(j)w'_{(3)})\circ I=L'_{P(z)}(j)(w'_{(3)}\circ I)
\end{equation}
for $w'_{(3)}\in W'_{3}$ and  $j=-1$, $0$ and $1$ (cf. 
(\ref{im:def'}) and Remark
\ref{I-intw}).
}
\end{rema}

\begin{rema}\label{L'jpreservesbetaspace}
{\rm We have 
\[
L'_{P(z)}(j)((W_{1}\otimes W_{2})^{*})^{(\beta)}\subset ((W_{1}\otimes
W_{2})^{*})^{(\beta)}
\]
for $j=-1, 0, 1$ and $\beta\in \tilde{A}$ (cf. Proposition \ref{3.6}).}
\end{rema}

When $V$ is a conformal vertex algebra, {}from the commutator formula
for $Y'_{P(z)}(\omega,x)$, we see that $L'_{P(z)}(-1)$, $L'_{P(z)}(0)$
and $L'_{P(z)}(1)$ realize the actions of $L_{-1}$, $L_0$ and $L_1$ in
${\mathfrak s}{\mathfrak l}(2)$ (recall (\ref{L_*})) on $(W_1\otimes
W_2)^*$.  When $V$ is just a M\"obius vertex algebra, the same
conclusion still holds but a proof is needed.  We state this as a
proposition:

\begin{propo}\label{sl-2}
Let $V$ be a M\"obius vertex algebra and let $W_{1}$ and $W_{2}$ be
generalized $V$-modules.  Then the operators $L'_{P(z)}(-1)$,
$L'_{P(z)}(0)$ and $L'_{P(z)}(1)$ realize the actions of $L_{-1}$,
$L_0$ and $L_1$ in ${\mathfrak s}{\mathfrak l}(2)$ on $(W_1\otimes
W_2)^*$, according to (\ref{L_*}).
\end{propo}
\pf
For $\lambda\in (W_{1}\otimes W_{2})^{*}$, $w_{(1)}\in W_{1}$,
$w_{(2)}\in W_{2}$ and $j, k=-1, 0, 1$, we have 
\begin{eqnarray}\label{kj}
\lefteqn{(L'_{P(z)}(j)L'_{P(z)}(k)\lambda)(w_{(1)}\otimes w_{(2)})}\nn
&&=(L'_{P(z)}(k)\lambda)
\Biggl(w_{(1)}\otimes L(-j)w_{(2)}+\Biggl(\sum_{i=0}
^{1-j}{1-j\choose i}z^iL(-j-i)\Biggr)w_{(1)}\otimes w_{(2)}
\Biggr)\nn
&&=(L'_{P(z)}(k)\lambda)
(w_{(1)}\otimes L(-j)w_{(2)})\nn
&&\quad +(L'_{P(z)}(k)\lambda)
\Biggl(\Biggl(\sum_{i=0}
^{1-j}{1-j\choose i}z^iL(-j-i)\Biggr)w_{(1)}\otimes w_{(2)}\Biggr)\nn
&&=\lambda
\Biggl(w_{(1)}\otimes L(-k)L(-j)w_{(2)}+\Biggl(\sum_{l=0}
^{1-k}{1-k\choose l}z^lL(-k-l)\Biggr)w_{(1)}\otimes L(-j)w_{(2)}\Biggr)
\nn
&&\quad +\lambda
\Biggl(\Biggl(\sum_{i=0}
^{1-j}{1-j\choose i}z^iL(-j-i)\Biggr)w_{(1)}\otimes L(-k)w_{(2)}\nn
&&\quad\quad\quad\quad +\Biggl(\sum_{l=0}
^{1-k}{1-k\choose l}z^lL(-k-l)\Biggr)\Biggl(\sum_{i=0}
^{1-j}{1-j\choose i}z^iL(-j-i)\Biggr)w_{(1)}\otimes w_{(2)}\Biggr).\nn
&&
\end{eqnarray}
{}From formula (\ref{kj}) we obtain
\begin{eqnarray}\label{kj-comm}
\lefteqn{([L'_{P(z)}(j), L'_{P(z)}(k)]\lambda)(w_{(1)}\otimes w_{(2)})}\nn
&&=\lambda
\Biggl(w_{(1)}\otimes [L(-k), L(-j)]w_{(2)}\nn
&&\quad\quad\quad +\Biggl(\sum_{l=0}^{1-k}\sum_{i=0}^{1-j}
{1-k\choose l}{1-j\choose i}z^{l+i}
[L(-k-l), L(-j-i)]\Biggr)w_{(1)}\otimes w_{(2)}\Biggr)\nn
&&=\lambda
\Biggl(w_{(1)}\otimes (j-k)L(-k-j)w_{(2)}\nn
&&\quad\quad\quad +\Biggl(\sum_{l=0}^{1-k}\sum_{i=0}^{1-j}
{1-k\choose l}{1-j\choose i}z^{l+i}(j+i-k-l)
L(-k-l-j-i)\Biggr)w_{(1)}\otimes w_{(2)}\Biggr)\nn
\end{eqnarray}

Taking $j=1$ and $k=-1, 0$ in (\ref{kj-comm}), 
we obtain 
\begin{eqnarray*}
\lefteqn{([L'_{P(z)}(1), L'_{P(z)}(k)]\lambda)(w_{(1)}\otimes w_{(2)})}\nn
&&=\lambda
\Biggl(w_{(1)}\otimes (1-k)L(-k-1)w_{(2)}\nn
&&\quad\quad\quad +\Biggl(\sum_{l=0}^{1-k}
{1-k\choose l}z^{l}(1-k-l)
L(-k-l-1)\Biggr)w_{(1)}\otimes w_{(2)}\Biggr)\nn
&&=\lambda
\Biggl(w_{(1)}\otimes (1-k)L(-k-1)w_{(2)}\nn
&&\quad\quad\quad +\Biggl(\sum_{l=0}^{1-(k+1)}
{1-k\choose l}z^{l}(1-k-l)
L(-k-l-1)\Biggr)w_{(1)}\otimes w_{(2)}\Biggr)\nn
&&=(1-k)\lambda
\Biggl(w_{(1)}\otimes L(-1-k)w_{(2)}\nn
&&\quad\quad\quad\quad\quad\quad\quad +\Biggl(\sum_{l=0}^{1-(k+1)}
{1-(k+1)\choose l}z^{l}
L(-(k+1)-l)\Biggr)w_{(1)}\otimes w_{(2)}\Biggr)\nn
&&=((1-k)L'_{P(z)}(1+k)\lambda)(w_{(1)}\otimes w_{(2)}),
\end{eqnarray*}
proving the commutator formula
\[
[L'_{P(z)}(1), L'_{P(z)}(k)]=(1-k)L'_{P(z)}(1+k)
\]
for $k=-1, 0$. 

Taking $j=0$ and $k=-1$ in (\ref{kj-comm}), 
we obtain 
\begin{eqnarray*}
\lefteqn{([L'_{P(z)}(0), L'_{P(z)}(-1)]\lambda)(w_{(1)}\otimes w_{(2)})}\nn
&&=\lambda
\Biggl(w_{(1)}\otimes L(1)w_{(2)} +\Biggl(\sum_{l=0}^{2}\sum_{i=0}^{1}
{2\choose l}z^{l+i}(i+1-l)
L(1-l-i)\Biggr)w_{(1)}\otimes w_{(2)}\Biggr)\nn
&&=\lambda(w_{(1)}\otimes L(1)w_{(2)} +
(L(1)+2zL(0)+z^{2}L(-1))w_{(1)}\otimes w_{(2)})\nn
&&=\lambda\Biggl(w_{(1)}\otimes L(1)w_{(2)} 
+\Biggl(\sum_{m=0}^{1-(-1)}{2\choose m}z^{m}L(-(-1)-m)\Biggr)
w_{(1)}\otimes w_{(2)}\Biggr)\nn
&&=(L'_{P(z)}(-1)\lambda)(w_{(1)}\otimes w_{(2)}),
\end{eqnarray*}
proving that
\[
[L'_{P(z)}(0), L'_{P(z)}(-1)]=L'_{P(z)}(-1).
\]
The three commutator formulas we have proved show that 
$L'_{P(z)}(-1)$,
$L'_{P(z)}(0)$ and $L'_{P(z)}(1)$ indeed realize the actions of $L_{-1}$,
$L_0$ and $L_1$ in ${\mathfrak s}{\mathfrak l}(2)$.
\epfv

The commutator formulas corresponding to (\ref{sl2-1})--(\ref{sl2-3})
(recall Definition \ref{moduleMobius})) also need to be proved in the
M\"obius case:

\begin{propo}\label{pz-l-y-comm}
Let $V$ be a M\"obius vertex algebra and let $W_{1}$ and $W_{2}$ be
generalized $V$-modules.  Then for $v\in V$, we have the following
commutator formulas:
\begin{eqnarray}
{[L(-1), Y'_{P(z)}(v, x)]}&=&Y'_{P(z)}(L(-1)v, x),\label{pz-sl-2-pz-y--2}\\
{[L(0), Y'_{P(z)}(v, x)]}&=&Y'_{P(z)}(L(0)v, x)+xY'_{P(z)}(L(-1)v, x),
\label{pz-sl-2-pz-y--1}\\
{[L(1), Y'_{P(z)}(v, x)]}&=&Y'_{P(z)}(L(1)v, x)
+2xY'_{P(z)}(L(0)v, x)+x^{2}Y'_{P(z)}(L(-1)v, x),\label{pz-sl-2-pz-y}\nn
&&
\end{eqnarray}
where for brevity we write $L'_{P(z)}(j)$ acting on $(W_1\otimes
W_2)^*$ as $L(j)$.
\end{propo}
\pf Let $\lambda\in (W_{1}\otimes W_{2})^{*}$, $w_{(1)}\in W_{1}$ and
$w_{(2)}\in W_{2}$. Using (\ref{LP'(j)}), (\ref{Y'def}), the
commutator formulas for $L(j)$ and $Y_{1}(v, x_{0})$ for $j=-1, 0, 1$
and $v\in V$ (recall Definition \ref{moduleMobius}), and the
commutator formulas for $L(j)$ and $Y_{2}^{o}(v, x)$ for $j=-1, 0, 1$
and $v\in V$ (recall Lemma \ref{sl2opposite}), we obtain, for $j=-1,
0, 1$,
\begin{eqnarray}\label{pz-sl-2-pz-y-1}
\lefteqn{([L(j), Y'_{P(z)}(v, x)]\lambda)(w_{(1)}\otimes w_{(2)})}\nn
&&=(L(j)Y'_{P(z)}(v, x)\lambda)(w_{(1)}\otimes w_{(2)})
-(Y'_{P(z)}(v, x)L(j)\lambda)(w_{(1)}\otimes w_{(2)})\nn
&&=(Y'_{P(z)}(v, x)\lambda)(w_{(1)}\otimes L(-j)w_{(2)})\nn
&&\quad +\sum_{i=0}^{1-j}{1-j\choose i}
(Y'_{P(z)}(v, x)\lambda)(z^iL(-j-i)w_{(1)}\otimes w_{(2)})\nn
&&\quad -(L(j)\lambda)(w_{(1)}\otimes Y_{2}^{o}(v, x)w_{(2)})\nn
&&\quad -\res_{x_{0}}z^{-1}\delta\left(\frac{x^{-1}-x_{0}}{z}\right)
(L(j)\lambda)(Y_{1}(e^{xL(1)}(-x^{-2})^{L(0)}v, x_{0})w_{(1)}
\otimes w_{(2)})\nn
&&=\lambda(w_{(1)}\otimes Y^{o}_{2}(v, x)L(-j)w_{(2)})\nn
&&\quad +\res_{x_{0}}z^{-1}\delta\left(\frac{x^{-1}-x_{0}}{z}\right)
\lambda(Y_{1}(e^{xL(1)}(-x^{-2})^{L(0)}v, x_{0})w_{(1)}
\otimes L(-j)w_{(2)})\nn
&&\quad +\sum_{i=0}^{1-j}{1-j\choose i}\lambda(z^iL(-j-i)w_{(1)}
\otimes Y_{2}^{o}(v, x)w_{(2)})\nn
&&\quad +\sum_{i=0}^{1-j}{1-j\choose i}
\res_{x_{0}}z^{-1}\delta\left(\frac{x^{-1}-x_{0}}{z}\right)\cdot\nn
&&\quad\quad\quad\quad\quad\quad\quad\quad\quad\quad\quad\quad \cdot
\lambda(Y_{1}(e^{xL(1)}(-x^{-2})^{L(0)}v, x_{0})
z^iL(-j-i)w_{(1)}\otimes w_{(2)})\nn
&&\quad -\lambda(w_{(1)}\otimes L(-j)Y_{2}^{o}(v, x)w_{(2)})\nn
&&\quad -\sum_{i=0}^{1-j}{1-j\choose i}
\lambda(z^iL(-j-i)w_{(1)}\otimes Y_{2}^{o}(v, x)w_{(2)})\nn
&&\quad -\res_{x_{0}}z^{-1}\delta\left(\frac{x^{-1}-x_{0}}{z}\right)
\lambda(Y_{1}(e^{xL(1)}(-x^{-2})^{L(0)}v, x_{0})w_{(1)}
\otimes L(-j)w_{(2)})\nn
&&\quad -\sum_{i=0}^{1-j}{1-j\choose i}
\res_{x_{0}}z^{-1}\delta\left(\frac{x^{-1}-x_{0}}{z}\right)\cdot\nn
&&\quad\quad\quad\quad\quad\quad\quad\quad\quad\quad\quad\quad \cdot
\lambda(z^iL(-j-i)Y_{1}(e^{xL(1)}(-x^{-2})^{L(0)}v, x_{0})w_{(1)}
\otimes w_{(2)})\nn
&&=\lambda(w_{(1)}\otimes [Y^{o}_{2}(v, x), L(-j)]w_{(2)})\nn
&&\quad -\sum_{i=0}^{1-j}{1-j\choose i}
\res_{x_{0}}z^{-1}\delta\left(\frac{x^{-1}-x_{0}}{z}\right)\cdot\nn
&&\quad\quad\quad\quad\quad\quad\quad\quad\quad\quad\quad\quad \cdot
\lambda([z^iL(-j-i), Y_{1}(e^{xL(1)}(-x^{-2})^{L(0)}v, x_{0})]
w_{(1)}\otimes w_{(2)})\nn
&&=\sum_{k=0}^{j+1}{j+1\choose k}x^{j+1-k}
\lambda(w_{(1)}\otimes Y^{o}_{2}(L(k-1)v, x)w_{(2)})\nn
&&\quad -\res_{x_{0}}z^{-1}\delta\left(\frac{x^{-1}-x_{0}}{z}\right)
\sum_{i=0}^{1-j}\sum_{k=0}^{-j-i+1}{1-j\choose i}
{-j-i+1\choose k}z^ix_{0}^{-j-i+1-k}
\cdot\nn
&&\quad\quad\quad\quad\quad\quad\quad\quad\quad\quad\quad\quad \cdot
\lambda(Y_{1}(L(k-1)e^{xL(1)}(-x^{-2})^{L(0)}v, x_{0})
w_{(1)}\otimes w_{(2)}).\nn
\end{eqnarray}
Using (\ref{2termdeltarelation}) and
(\ref{deltafunctionsubstitutionformula}) when necessary, we see that
the second term on the right-hand side of (\ref{pz-sl-2-pz-y-1}) is
equal to the following expressions for $j=1, 0$ and $-1$, respectively:
\begin{eqnarray}\label{pz-sl-2-pz-y-2}
-\res_{x_{0}}z^{-1}\delta\left(\frac{x^{-1}-x_{0}}{z}\right)
\lambda(Y_{1}(L(-1)e^{xL(1)}(-x^{-2})^{L(0)}v, x_{0})
w_{(1)}\otimes w_{(2)}),
\end{eqnarray}
\begin{eqnarray}\label{pz-sl-2-pz-y-3}
\lefteqn{-\res_{x_{0}}z^{-1}\delta\left(\frac{x^{-1}-x_{0}}{z}\right)
\sum_{i=0}^{1}\sum_{k=0}^{-i+1}{1\choose i}
{1-i\choose k}z^ix_{0}^{1-i-k}
\cdot}\nn
&&\quad\quad\quad\quad\quad\quad\quad\cdot
\lambda(Y_{1}(L(k-1)e^{xL(1)}(-x^{-2})^{L(0)}v, x_{0})
w_{(1)}\otimes w_{(2)})\nn
&&=-\res_{x_{0}}z^{-1}\delta\left(\frac{x^{-1}-x_{0}}{z}\right)
x_{0}
\lambda(Y_{1}(L(-1)e^{xL(1)}(-x^{-2})^{L(0)}v, x_{0})
w_{(1)}\otimes w_{(2)})\nn
&&\quad -\res_{x_{0}}z^{-1}\delta\left(\frac{x^{-1}-x_{0}}{z}\right)
\lambda(Y_{1}(L(0)e^{xL(1)}(-x^{-2})^{L(0)}v, x_{0})
w_{(1)}\otimes w_{(2)})\nn
&&\quad -\res_{x_{0}}z^{-1}\delta\left(\frac{x^{-1}-x_{0}}{z}\right)
z\lambda(Y_{1}(L(-1)e^{xL(1)}(-x^{-2})^{L(0)}v, x_{0})
w_{(1)}\otimes w_{(2)})\nn
&&=-\res_{x_{0}}x\delta\left(\frac{z+x_{0}}{x^{-1}}\right)
(x_{0}+z)
\lambda(Y_{1}(L(-1)e^{xL(1)}(-x^{-2})^{L(0)}v, x_{0})
w_{(1)}\otimes w_{(2)})\nn
&&\quad -\res_{x_{0}}z^{-1}\delta\left(\frac{x^{-1}-x_{0}}{z}\right)
\lambda(Y_{1}(L(0)e^{xL(1)}(-x^{-2})^{L(0)}v, x_{0})
w_{(1)}\otimes w_{(2)})\nn
&&=-\res_{x_{0}}z^{-1}\delta\left(\frac{x^{-1}-x_{0}}{z}\right)x^{-1}
\lambda(Y_{1}(L(-1)e^{xL(1)}(-x^{-2})^{L(0)}v, x_{0})
w_{(1)}\otimes w_{(2)})\nn
&&\quad -\res_{x_{0}}z^{-1}\delta\left(\frac{x^{-1}-x_{0}}{z}\right)
\lambda(Y_{1}(L(0)e^{xL(1)}(-x^{-2})^{L(0)}v, x_{0})
w_{(1)}\otimes w_{(2)})\nn
\end{eqnarray}
and
\begin{eqnarray}\label{pz-sl-2-pz-y-4}
\lefteqn{-\res_{x_{0}}z^{-1}\delta\left(\frac{x^{-1}-x_{0}}{z}\right)
\sum_{i=0}^{2}\sum_{k=0}^{2-i}{2\choose i}
{2-i\choose k}z^ix_{0}^{2-i-k}
\cdot}\nn
&&\quad\quad\quad\quad\quad\quad\quad\quad\quad\quad\quad\quad \cdot
\lambda(Y_{1}(L(k-1)e^{xL(1)}(-x^{-2})^{L(0)}v, x_{0})
w_{(1)}\otimes w_{(2)})\nn
&&=-\res_{x_{0}}z^{-1}\delta\left(\frac{x^{-1}-x_{0}}{z}\right)
x_{0}^{2}
\lambda(Y_{1}(L(-1)e^{xL(1)}(-x^{-2})^{L(0)}v, x_{0})
w_{(1)}\otimes w_{(2)})\nn
&&\quad -\res_{x_{0}}z^{-1}\delta\left(\frac{x^{-1}-x_{0}}{z}\right)
2x_{0}
\lambda(Y_{1}(L(0)e^{xL(1)}(-x^{-2})^{L(0)}v, x_{0})
w_{(1)}\otimes w_{(2)})\nn
&&\quad -\res_{x_{0}}z^{-1}\delta\left(\frac{x^{-1}-x_{0}}{z}\right)
\lambda(Y_{1}(L(1)e^{xL(1)}(-x^{-2})^{L(0)}v, x_{0})
w_{(1)}\otimes w_{(2)})\nn
&&\quad -\res_{x_{0}}z^{-1}\delta\left(\frac{x^{-1}-x_{0}}{z}\right)
2zx_{0}
\lambda(Y_{1}(L(-1)e^{xL(1)}(-x^{-2})^{L(0)}v, x_{0})
w_{(1)}\otimes w_{(2)})\nn
&&\quad -\res_{x_{0}}z^{-1}\delta\left(\frac{x^{-1}-x_{0}}{z}\right)
2z\lambda(Y_{1}(L(0)e^{xL(1)}(-x^{-2})^{L(0)}v, x_{0})
w_{(1)}\otimes w_{(2)})\nn
&&\quad -\res_{x_{0}}z^{-1}\delta\left(\frac{x^{-1}-x_{0}}{z}\right)
z^2\lambda(Y_{1}(L(-1)e^{xL(1)}(-x^{-2})^{L(0)}v, x_{0})
w_{(1)}\otimes w_{(2)})\nn
&&=-\res_{x_{0}}z^{-1}\delta\left(\frac{x^{-1}-x_{0}}{z}\right)
x^{-2}
\lambda(Y_{1}(L(-1)e^{xL(1)}(-x^{-2})^{L(0)}v, x_{0})
w_{(1)}\otimes w_{(2)})\nn
&&\quad -\res_{x_{0}}z^{-1}\delta\left(\frac{x^{-1}-x_{0}}{z}\right)
2x^{-1}
\lambda(Y_{1}(L(0)e^{xL(1)}(-x^{-2})^{L(0)}v, x_{0})
w_{(1)}\otimes w_{(2)})\nn
&&\quad -\res_{x_{0}}z^{-1}\delta\left(\frac{x^{-1}-x_{0}}{z}\right)
\lambda(Y_{1}(L(1)e^{xL(1)}(-x^{-2})^{L(0)}v, x_{0})
w_{(1)}\otimes w_{(2)}).\nn
\end{eqnarray}
Using (\ref{log:SL2-3}) and (\ref{log:xLx^}), we see that
(\ref{pz-sl-2-pz-y-2}), (\ref{pz-sl-2-pz-y-3}) and
(\ref{pz-sl-2-pz-y-4}) are respectively equal to
\begin{eqnarray}\label{pz-sl-2-pz-y-5}
\lefteqn{-\res_{x_{0}}z^{-1}\delta\left(\frac{x^{-1}-x_{0}}{z}\right)
\cdot}\nn
&&\quad\quad
\cdot\lambda(Y_{1}(e^{xL(1)}(x^{2}L(1)-2xL(0)+L(-1))(-x^{-2})^{L(0)}v, x_{0})
w_{(1)}\otimes w_{(2)})\nn
&&=\res_{x_{0}}z^{-1}\delta\left(\frac{x^{-1}-x_{0}}{z}\right)\cdot\nn
&&\quad\quad
\cdot
\lambda(Y_{1}(e^{xL(1)}(-x^{-2})^{L(0)}(L(1)+2xL(0)+x^{2}L(-1))v, x_{0})
w_{(1)}\otimes w_{(2)}),\nn
\end{eqnarray}
\begin{eqnarray}\label{pz-sl-2-pz-y-6}
\lefteqn{-\res_{x_{0}}z^{-1}\delta\left(\frac{x^{-1}-x_{0}}{z}\right)x^{-1}
\cdot}\nn
&&\quad\quad\cdot
\lambda(Y_{1}(e^{xL(1)}(x^{2}L(1)-2xL(0)+L(-1))(-x^{-2})^{L(0)}v, x_{0})
w_{(1)}\otimes w_{(2)})\nn
&&\quad -\res_{x_{0}}z^{-1}\delta\left(\frac{x^{-1}-x_{0}}{z}\right)
\lambda(Y_{1}(e^{xL(1)}(-xL(1)+L(0))(-x^{-2})^{L(0)}v, x_{0})
w_{(1)}\otimes w_{(2)})\nn
&&=\res_{x_{0}}z^{-1}\delta\left(\frac{x^{-1}-x_{0}}{z}\right)x^{-1}\cdot\nn
&&\quad\quad\cdot
\lambda(Y_{1}(e^{xL(1)}(-x^{-2})^{L(0)}(L(1)+2xL(0)+x^{2}L(-1))v, x_{0})
w_{(1)}\otimes w_{(2)})\nn
&&\quad +\res_{x_{0}}z^{-1}\delta\left(\frac{x^{-1}-x_{0}}{z}\right)
\lambda(Y_{1}(e^{xL(1)}(-x^{-2})^{L(0)}(-x^{-1}L(1)-L(0))v, x_{0})
w_{(1)}\otimes w_{(2)})\nn
&&=\res_{x_{0}}z^{-1}\delta\left(\frac{x^{-1}-x_{0}}{z}\right)
\lambda(Y_{1}(e^{xL(1)}(-x^{-2})^{L(0)}(L(0)+xL(-1))v, x_{0})
w_{(1)}\otimes w_{(2)})\nn
\end{eqnarray}
and
\begin{eqnarray}\label{pz-sl-2-pz-y-7}
\lefteqn {-\res_{x_{0}}z^{-1}\delta\left(\frac{x^{-1}-x_{0}}{z}\right)
x^{-2}\cdot}\nn
&&\quad\quad\cdot
\lambda(Y_{1}(e^{xL(1)}(x^{2}L(1)-2xL(0)+L(-1))(-x^{-2})^{L(0)}v, x_{0})
w_{(1)}\otimes w_{(2)})\nn
&&\quad -\res_{x_{0}}z^{-1}\delta\left(\frac{x^{-1}-x_{0}}{z}\right)
2x^{-1}
\lambda(Y_{1}(e^{xL(1)}(-xL(1)+L(0))(-x^{-2})^{L(0)}v, x_{0})
w_{(1)}\otimes w_{(2)})\nn
&&\quad -\res_{x_{0}}z^{-1}\delta\left(\frac{x^{-1}-x_{0}}{z}\right)
\lambda(Y_{1}(e^{xL(1)}L(1)(-x^{-2})^{L(0)}v, x_{0})
w_{(1)}\otimes w_{(2)})\nn
&&=\res_{x_{0}}z^{-1}\delta\left(\frac{x^{-1}-x_{0}}{z}\right)
x^{-2}\cdot\nn
&&\quad\quad
\cdot
\lambda(Y_{1}(e^{xL(1)}(-x^{-2})^{L(0)}
(L(1)+2xL(0)+x^{2}L(-1))v, x_{0})
w_{(1)}\otimes w_{(2)})\nn
&&\quad +\res_{x_{0}}z^{-1}\delta\left(\frac{x^{-1}-x_{0}}{z}\right)
2x^{-1}
\lambda(Y_{1}(e^{xL(1)}(-x^{-2})^{L(0)}(-x^{-1}L(1)-L(0))v, x_{0})
w_{(1)}\otimes w_{(2)})\nn
&&\quad +\res_{x_{0}}z^{-1}\delta\left(\frac{x^{-1}-x_{0}}{z}\right)
\lambda(Y_{1}(e^{xL(1)}(-x^{-2})^{L(0)}x^{-2}L(1)v, x_{0})
w_{(1)}\otimes w_{(2)})\nn
&&=\res_{x_{0}}z^{-1}\delta\left(\frac{x^{-1}-x_{0}}{z}\right)
\lambda(Y_{1}(e^{xL(1)}(-x^{-2})^{L(0)}
L(-1)v, x_{0})
w_{(1)}\otimes w_{(2)}).
\end{eqnarray}
The right-hand sides of (\ref{pz-sl-2-pz-y-5}),
(\ref{pz-sl-2-pz-y-6}) and   (\ref{pz-sl-2-pz-y-7}) can be written as 
\[
\sum_{k=0}^{j+1}{j+1\choose k}x^{j+1-k}
\res_{x_{0}}z^{-1}\delta\left(\frac{x^{-1}-x_{0}}{z}\right)
\lambda(Y_{1}(e^{xL(1)}(-x^{-2})^{L(0)}
L(k-1)v, x_{0})
w_{(1)}\otimes w_{(2)}),
\]
for $j=1, 0, -1$, respectively. Thus the right-hand side of 
(\ref{pz-sl-2-pz-y-1}) is equal to 
\begin{eqnarray}\label{pz-sl-2-pz-y-8}
\lefteqn{\sum_{k=0}^{j+1}{j+1\choose k}x^{j+1-k}
\lambda(w_{(1)}\otimes Y^{o}_{2}(L(k-1)v, x)w_{(2)})}\nn
&&\quad +\sum_{k=0}^{j+1}{j+1\choose k}x^{j+1-k}
\res_{x_{0}}z^{-1}\delta\left(\frac{x^{-1}-x_{0}}{z}\right)\cdot\nn
&&\quad\quad\quad\quad\quad\quad\quad\quad\quad\quad\quad\quad\quad\quad
\cdot
\lambda(Y_{1}(e^{xL(1)}(-x^{-2})^{L(0)}
L(k-1)v, x_{0})
w_{(1)}\otimes w_{(2)})\nn
&&=\sum_{k=0}^{j+1}{j+1\choose k}x^{j+1-k}
(Y'_{P(z)}(L(k-1)v, x)\lambda)(w_{(1)}\otimes w_{(2)}),
\end{eqnarray}
proving the proposition.
\epfv

We have seen in (\ref{tauw}), (\ref{deltaY3'}) and (\ref{tausubW3'})
that for a generalized $V$-module $(W,Y_W)$, the space $V \otimes
{\C}((t))$, and in particular, the space $V \otimes \iota_{+}{\mathbb
C}[t,t^{- 1}, (z^{-1}-t)^{-1}]$, acts naturally on $W$ via the action
$\tau_W$, in view of (\ref{set:wtvn}) and Assumption \ref{assum};
recall that $v \otimes t^n$ ($v \in V$, $n \in {\mathbb Z}$) acts as
the component $v_n$ of $Y_W(v,x)$, and that more generally,
\begin{equation}\label{tau-w-comp}
\tau_W \left(v \otimes \sum_{n > N} a_nt^n\right) = \sum_{n > N} a_nv_n
\end{equation}
for $a_n \in {\mathbb C}$.  For generalized $V$-modules $W_1$, $W_2$
and $W_3$, we shall next relate the $P(z)$-intertwining maps of type
${W_3\choose W_1\, W_2}$ to certain linear maps {}from $W'_{3}$ to
$(W_1\otimes W_2)^{*}$ intertwining the actions of $V \otimes
\iota_{+}{\mathbb C}[t,t^{- 1}, (z^{-1}-t)^{-1}]$ and of ${\mathfrak
s} {\mathfrak l}(2)$ on $W'_{3}$ and on $(W_1\otimes W_2)^{*}$ (see
Proposition \ref{pz} and Notation \ref{scriptN} below).  For this, as
is suggested by Lemma \ref{4.36} and Proposition \ref{tau-a-comp}, we
need to consider $\tilde{A}$-compatibility for linear maps {}from $W_3$
to $(W_1\otimes W_2)^{*}$:

\begin{defi}\label{defJAtildecompat}
{\rm We call a map $J \in {\rm Hom}(W_3,(W_1\otimes W_2)^{*})$ 
{\it $\tilde{A}$-compatible} if
\begin{equation}\label{JAtildecompat}
J((W_{3})^{(\beta)}) \subset ((W_{1}\otimes W_{2})^{*})^{(\beta)}
\end{equation}
for $\beta \in {\tilde A}$.
}
\end{defi}

As in the discussion preceding Lemma \ref{4.36}, we see that an
element $\lambda$ of $(W_{1}\otimes W_{2}\otimes W_{3})^{*}$ amounts
exactly to a linear map
\[
J_{\lambda}: W_3 \rightarrow (W_{1}\otimes W_{2})^{*}.
\]
If $\lambda$ is $\tilde{A}$-compatible (see that discussion), then for
$w_{(1)}\in W_{1}^{(\beta)}$, $w_{(2)}\in W_{2}^{(\gamma)}$ and
$w_{(3)}\in W_{3}^{(\delta)}$ such that $\beta +\gamma +\delta \ne 0$,
\[
J_{\lambda}(w_{(3)})(w_{(1)}\otimes w_{(2)})=\lambda(w_{(1)}\otimes
w_{(2)}\otimes w_{(3)}) = 0,
\]
so that
\[
J_{\lambda}(w_{(3)}) \in ((W_{1}\otimes W_{2})^{*})^{(\delta)},
\]
and so $J_{\lambda}$ is $\tilde{A}$-compatible.  Similarly, if
$J_{\lambda}$ is $\tilde{A}$-compatible, then so is $\lambda$.  Thus
using Lemma \ref{4.36} we have:

\begin{lemma}\label{IlambdatoJlambda}
The linear functional $\lambda\in (W_{1}\otimes W_{2}\otimes
W_{3})^{*}$ is $\tilde{A}$-compatible if and only if $J_{\lambda}$ is
$\tilde{A}$-compatible.  The map given by $\lambda\mapsto J_{\lambda}$
is the unique linear isomorphism {}from the space of
$\tilde{A}$-compatible elements of $(W_{1}\otimes W_{2}\otimes
W_{3})^{*}$ to the space of $\tilde{A}$-compatible linear maps {}from
$W_3$ to $(W_{1}\otimes W_{2})^{*}$ such that
\[
J_{\lambda}(w_{(3)})(w_{(1)}\otimes w_{(2)})
=\lambda(w_{(1)}\otimes w_{(2)}\otimes w_{(3)})
\]
for $w_{(1)}\in W_{1}$, $w_{(2)}\in W_{2}$ and $w_{(3)}\in W_{3}$.  In
particular, the correspondence $I_{\lambda} \mapsto J_{\lambda}$
defines a (unique) linear isomorphism {}from the space of
$\tilde{A}$-compatible linear maps
\[
I = I_{\lambda}:W_{1}\otimes W_{2} \rightarrow \overline{W_{3}'}
\]
to the space of $\tilde{A}$-compatible linear maps
\[
J = J_{\lambda}: W_3 \rightarrow (W_{1}\otimes W_{2})^{*}
\]
such that
\[
\langle w_{(3)}, I(w_{(1)}\otimes w_{(2)})\rangle
=J(w_{(3)})(w_{(1)}\otimes w_{(2)})
\]
for $w_{(1)}\in W_{1}$, $w_{(2)}\in W_{2}$ and $w_{(3)}\in W_{3}$.
\epf
\end{lemma}

\begin{rema}\label{alternateformoflemma}{\rm
{}From Lemma \ref{IlambdatoJlambda} (with $W_3$ replaced by $W'_{3}$) we
have a canonical isomorphism {}from the space of
$\tilde{A}$-compatible linear maps
\[
I:W_{1}\otimes W_{2} \rightarrow \overline{W_3}
\]
to the space of $\tilde{A}$-compatible linear maps
\[
J:W'_{3} \rightarrow (W_{1}\otimes W_{2})^{*},
\]
determined by:
\begin{equation}\label{IcorrespondstoJ}
\langle w'_{(3)}, I(w_{(1)}\otimes w_{(2)})\rangle
=J(w'_{(3)})(w_{(1)}\otimes w_{(2)})
\end{equation}
for $w_{(1)}\in W_{1}$, $w_{(2)}\in W_{2}$ and $w'_{(3)}\in W'_{3}$,
or equivalently,
\begin{equation}\label{IcorrespondstoJalternateform}
w'_{(3)}\circ I = J(w'_{(3)})
\end{equation}
for $w'_{(3)} \in W'_{3}$.
}
\end{rema}

We introduce another notion, corresponding to the lower truncation
condition (\ref{im:ltc}) for $P(z)$-intertwining maps:

\begin{defi}\label{gradingrestrictedmapJ}{\rm
We call a map $J\in \hom(W_3, (W_1\otimes W_2)^{*})$ {\em grading
restricted} if for $n\in {\mathbb C}$, $w_{(1)}\in W_1$ and
$w_{(2)}\in W_2$,
\begin{equation}\label{Jgradingrestr}
J((W_3)_{[n-m]})(w_{(1)}\otimes w_{(2)})=0\;\;\mbox{ for }\;m\in
{\mathbb N}\;\mbox{ sufficiently large.}
\end{equation}
}
\end{defi}

\begin{rema}\label{Jcompatimpliesgradingrestr}
{\rm
If $J \in {\rm Hom}(W_3,(W_1\otimes W_2)^{*})$ is
$\tilde{A}$-compatible, then $J$ is also grading restricted, as we see
using (\ref{set:dmltc}).
}
\end{rema}

\begin{rema}\label{Jlowerbounded}
{\rm If in addition $W_3$ is lower bounded (recall
(\ref{set:dmltc-1})), then the stronger condition
\begin{equation}\label{Jlowerbdd}
J((W_3)_{[n]})(w_{(1)}\otimes w_{(2)})=0\;\;\mbox{ for }
\;\Re{(n)}\;\mbox{ sufficiently negative}
\end{equation}
holds.}
\end{rema}

\begin{rema}{\rm
Under the natural isomorphism given in Remark
\ref{alternateformoflemma} (see (\ref{IcorrespondstoJ})) in the
$\tilde{A}$-compatible setting, the map $J : W'_{3} \rightarrow
(W_{1}\otimes W_{2})^{*}$ is grading restricted (recall Definition
\ref{gradingrestrictedmapJ}) if and only if the map $I:W_{1}\otimes
W_{2} \rightarrow \overline{W_{3}}$ satisfies the lower truncation
condition (\ref{im:ltc}).  But notice also that in this
$\tilde{A}$-compatible setting, we have seen that both $I$ and $J$
automatically have these properties.
}
\end{rema}

\begin{rema}{\rm Analogous comments hold for the stronger conditions
(\ref{PpinI=0}) and (\ref{Jlowerbounded}) in case $W_3$ is lower
bounded.}
\end{rema}

Using the above together with Remarks \ref{I-intw} and \ref{I-intw2},
we now have the following result, generalizing Proposition 13.1 
in \cite{tensor3}:

\begin{propo}\label{pz}
Let $W_1$, $W_2$ and $W_3$ be generalized $V$-modules.  Under the
natural isomorphism described in Remark \ref{alternateformoflemma}
between the space of $\tilde{A}$-compatible linear maps
\[
I:W_{1}\otimes W_{2} \rightarrow \overline{W_{3}}
\]
and the space of $\tilde{A}$-compatible linear maps
\[
J:W'_{3} \rightarrow (W_{1}\otimes W_{2})^{*}
\]
determined by (\ref{IcorrespondstoJ}), the $P(z)$-intertwining maps
$I$ of type ${W_3\choose W_1\, W_2}$ correspond exactly to the
(grading restricted) $\tilde{A}$-compatible maps $J$ that intertwine
the actions of both
\[
V \otimes \iota_{+}{\mathbb C}[t,t^{- 1}, (z^{-1}-t)^{-1}]
\]
and ${\mathfrak s} {\mathfrak l}(2)$ on $W'_{3}$ and on $(W_1\otimes
W_2)^{*}$.  If $W_3$ is lower bounded, we may replace the grading
restrictions by (\ref{PpinI=0}) and (\ref{Jlowerbdd}).
\end{propo}
\pf 
In view of (\ref{IcorrespondstoJalternateform}), Remark
\ref{I-intw} asserts that (\ref{im:def'}), or equivalently,
(\ref{im:def}), is equivalent to the condition
\begin{equation}\label{j-tau}
J\left(\tau_{W'_3}\left(x_0^{-1}\delta\left(\frac{x^{-1}_1-z}{x_0}\right)
Y_{t}(v, x_1)\right)w'_{(3)}\right)
=\tau_{P(z)}\left(x_0^{-1}\delta\left(\frac{x^{-1}_1-z}{x_0}\right)
Y_{t}(v, x_1)\right)J(w'_{(3)}),
\end{equation}
that is, the condition that $J$ intertwines the actions of $V \otimes
\iota_{+}{\mathbb C}[t,t^{- 1}, (z^{-1}-t)^{-1}]$ on $W'_{3}$ and on
$(W_1\otimes W_2)^{*}$ (recall (\ref{3.18-1})--(\ref{3.19-1})).
Similarly, Remark \ref{I-intw2} asserts that (\ref{im:Lj}) is
equivalent to the condition
\begin{equation}\label{j-lj}
J(L'(j)w'_{(3)}) = L'_{P(z)}(j)J(w'_{(3)})
\end{equation}
for $j=-1$, $0$, $1$, that is, the condition that $J$ intertwines the
actions of ${\mathfrak s} {\mathfrak l}(2)$ on $W'_{3}$ and on
$(W_1\otimes W_2)^{*}$.
\epfv

\begin{nota}\label{scriptN}
{\rm Given generalized $V$-modules $W_1$, $W_2$ and $W_3$, we shall
write ${\cal N}[P(z)]_{W'_3}^{(W_1 \otimes W_2)^{*}}$, or ${\cal
N}_{W'_3}^{(W_1 \otimes W_2)^{*}}$ if there is no ambiguity, for the
space of (grading restricted) $\tilde{A}$-compatible linear maps
\[
J:W'_{3} \rightarrow (W_{1}\otimes W_{2})^{*}
\]
that intertwine the actions of both
\[
V \otimes \iota_{+}{\mathbb C}[t,t^{- 1}, (z^{-1}-t)^{-1}]
\]
and ${\mathfrak s} {\mathfrak l}(2)$ on $W'_{3}$ and on $(W_1\otimes
W_2)^{*}$.
Note that Proposition \ref{pz} gives a natural linear isomorphism
\begin{eqnarray*}
{\cal M}[P(z)]^{W_3}_{W_1 W_2} = {\cal M}^{W_3}_{W_1 W_2} &
\stackrel{\sim}{\longrightarrow} & {\cal N}_{W'_3}^{(W_1 \otimes
W_2)^{*}}\nno\\ I & \mapsto & J
\end{eqnarray*}
(recall {}from Definition \ref{im:imdef} the notations for the space
of $P(z)$-intertwining maps), and if $W_3$ is lower bounded, the
spaces satisfy the stronger grading restrictions (\ref{PpinI=0}) and
(\ref{Jlowerbdd}).  Let us use the symbol ``prime'' to denote this
isomorphism in both directions:
\begin{eqnarray*}
{\cal M}^{W_3}_{W_1 W_2} & \stackrel{\sim}{\longrightarrow} & {\cal
N}_{W'_3}^{(W_1 \otimes W_2)^{*}}\nno\\
I & \mapsto & I'\nno\\
J' & \leftarrow\!\!\!{\scriptstyle |} & J,
\end{eqnarray*}
so that in particular,
\[
I'' = I \;\;\mbox{ and }\;\; J'' = J
\]
for $I \in {\cal M}^{W_3}_{W_1 W_2}$ and $J \in {\cal N}_{W'_3}^{(W_1
\otimes W_2)^{*}}$, and the relation between $I$ and $I'$ is
determined by
\[
\langle w'_{(3)}, I(w_{(1)}\otimes w_{(2)})\rangle
=I'(w'_{(3)})(w_{(1)}\otimes w_{(2)})
\]
for $w_{(1)}\in W_{1}$, $w_{(2)}\in W_{2}$ and $w'_{(3)}\in W'_{3}$,
or equivalently,
\[
w'_{(3)}\circ I = I'(w'_{(3)}).
\]
}
\end{nota}

\begin{rema}\label{NisotoV}
{\rm Combining Proposition \ref{pz} with Proposition
\ref{im:correspond}, we see that for any integer $p$, we also have a
natural linear isomorphism
\[
{\cal N}_{W'_3}^{(W_1 \otimes W_2)^{*}} \stackrel{\sim}{\longrightarrow}
{\cal V}^{W_3}_{W_1 W_2}
\]
{}from ${\cal N}_{W'_3}^{(W_1 \otimes W_2)^{*}}$ to the space of
logarithmic intertwining operators of type ${W_3\choose W_1\,W_2}$.
In particular, given a logarithmic intertwining operator ${\cal Y}$ of
type ${W_3\choose W_1\,W_2}$, the map
\[
I'_{{\cal Y},p}: W'_3\to (W_1\otimes W_2)^*
\]
defined by
\[
I'_{{\cal Y},p}(w'_{(3)})(w_{(1)}\otimes w_{(2)})=\langle
w'_{(3)},{\cal Y}(w_{(1)},e^{l_p(z)})w_{(2)}\rangle_{W_3}
\]
is $\tilde{A}$-compatible and intertwines both actions on both spaces.
If $W_3$ is lower bounded, we have the stronger grading restrictions
(recall Remark \ref{lowerbddcorrespondence}).}
\end{rema}

Recall that we have formulated the notions of $P(z)$-product and
$P(z)$-tensor product using $P(z)$-intertwining maps (Definitions
\ref{pz-product} and \ref{pz-tp}).  Since we now know that
$P(z)$-intertwining maps can be interpreted as in Proposition \ref{pz}
(and Notation \ref{scriptN}), we can easily reformulate the notions of
$P(z)$-product and $P(z)$-tensor product correspondingly:
\begin{propo}\label{productusingI'}
Let ${\cal C}_1$ be either of the categories ${\cal M}_{sg}$ or ${\cal
GM}_{sg}$, as in Definition \ref{pz-product}.  For $W_1, W_2\in
\ob{\cal C}_1$, a $P(z)$-product $(W_3;I_3)$ of $W_1$ and $W_2$
(recall Definition \ref{pz-product}) amounts to an object $(W_3,Y_3)$
of ${\cal C}_1$ equipped with a map $I'_3 \in {\cal N}_{W'_3}^{(W_1
\otimes W_2)^{*}}$, that is, equipped with an $\tilde{A}$-compatible
map
\[
I'_3:W'_{3} \rightarrow (W_{1}\otimes W_{2})^{*}
\]
that intertwines the two actions of $V \otimes \iota_{+}{\mathbb
C}[t,t^{- 1}, (z^{-1}-t)^{-1}]$ and of ${\mathfrak s} {\mathfrak
l}(2)$.  The map $I'_3$ corresponds to the $P(z)$-intertwining map
\[
I_3:W_{1}\otimes W_{2} \rightarrow \overline{W_{3}}
\]
as above:
\[
I'_3(w'_{(3)}) = w'_{(3)}\circ I_3
\]
for $w'_{(3)} \in W'_{3}$ (recall
\ref{IcorrespondstoJalternateform})).  Denoting this structure by
$(W_3,Y_3;I'_3)$ or simply by $(W_3;I'_3)$, let $(W_4;I'_4)$ be
another such structure.  Then a morphism of $P(z)$-products {}from $W_3$
to $W_4$ amounts to a module map $\eta: W_3 \to W_4$ such that the
diagram
\begin{center}
\begin{picture}(100,60)
\put(-2,0){$W'_4$}
\put(13,4){\vector(1,0){104}}
\put(119,0){$W'_3$}
\put(38,50){$(W_1\otimes W_2)^*$}
\put(13,12){\vector(3,2){50}}
\put(118,12){\vector(-3,2){50}}
\put(65,8){$\eta'$}
\put(23,27){$I'_4$}
\put(98,27){$I'_3$}
\end{picture}
\end{center}
commutes, where $\eta'$ is the natural map given by (\ref{fprime}).
\end{propo}
\pf All we need to check is that the diagram in Definition
\ref{pz-product} commutes if and only if the diagram above commutes.
But this follows {}from the definitions and the fact that for
$\overline{w_{(3)}} \in \overline{W_{3}}$ and $w'_{(4)} \in W'_{4}$,
\[
\langle \eta'(w'_{(4)}),\overline{w_{(3)}} \rangle=\langle
w'_{(4)},\overline{\eta}(\overline{w_{(3)}})\rangle,
\]
which in turn follows {}from (\ref{fprime}).  \epfv

\begin{corol}\label{tensorproductusingI'}
Let ${\cal C}$ be a full subcategory of either ${\cal M}_{sg}$ or
${\cal GM}_{sg}$, as in Definition \ref{pz-tp}.  For $W_1, W_2\in
\ob{\cal C}$, a $P(z)$-tensor product $(W_0; I_0)$ of $W_1$ and $W_2$
in ${\cal C}$, if it exists, amounts to an object $W_0 =
W_1\boxtimes_{P(z)} W_2$ of ${\cal C}$ and a structure $(W_0 =
W_1\boxtimes_{P(z)} W_2; I'_0)$ as in Proposition
\ref{productusingI'}, with
\[
I'_0: (W_1\boxtimes_{P(z)} W_2)' \longrightarrow (W_1\otimes W_2)^*
\]
in ${\cal N}_{(W_1\boxtimes_{P(z)} W_2)'}^{(W_1 \otimes W_2)^{*}}$,
such that for any such pair $(W; I')$ $(W\in \ob \mathcal{C})$, with
\[
I': W' \longrightarrow (W_1\otimes W_2)^*
\]
in ${\cal N}_{W'}^{(W_1 \otimes W_2)^{*}}$, there is a unique module
map
\[
\chi: W' \longrightarrow (W_1\boxtimes_{P(z)} W_2)'
\]
such that the diagram
\begin{center}
\begin{picture}(100,60)
\put(-2,0){$W'$}
\put(13,4){\vector(1,0){104}}
\put(119,0){$(W_1\boxtimes_{P(z)} W_2)'$}
\put(38,50){$(W_1\otimes W_2)^*$}
\put(13,12){\vector(3,2){50}}
\put(118,12){\vector(-3,2){50}}
\put(65,8){$\chi$}
\put(23,27){$I'$}
\put(98,27){$I'_0$}
\end{picture}
\end{center}
commutes.  Here $\chi = \eta'$, where $\eta$ is a correspondingly
unique module map
\[
\eta: W_1\boxtimes_{P(z)} W_2 \longrightarrow W.
\]
Also, the map $I_0'$, which is $\tilde{A}$-compatible and which
intertwines the two actions of $V \otimes \iota_{+}{\mathbb C}[t,t^{-
1}, (z^{-1}-t)^{-1}]$ and of ${\mathfrak s} {\mathfrak l}(2)$, is
related to the $P(z)$-intertwining map
\[
I_0 = \boxtimes_{P(z)}: W_1\otimes W_2 \longrightarrow 
\overline{W_1\boxtimes_{P(z)} W_2}
\]
by
\[
I_0'(w') = w' \circ \boxtimes_{P(z)}
\]
for $w' \in (W_1\boxtimes_{P(z)} W_2)'$, that is,
\[
I_0'(w')(w_{(1)}\otimes w_{(2)}) = \langle w',w_{(1)}\boxtimes_{P(z)}
w_{(2)} \rangle
\]
for $w_{(1)}\in W_{1}$ and $w_{(2)}\in W_{2}$, using the notation 
(\ref{boxtensorofelements}).  
\epf
\end{corol}

\begin{rema}\label{motivationofbackslash}{\rm
{}From Corollary \ref{tensorproductusingI'} we see that it is natural to
try to construct $W_1\boxtimes_{P(z)} W_2$, when it exists, as the
contragredient of a suitable natural substructure of $(W_1\otimes
W_2)^*$.  We shall now proceed to do this.  Under suitable assumptions, 
we shall in fact construct
a module-like structure
\[
W_1\hboxtr_{P(z)} W_2 \subset (W_1\otimes W_2)^*
\]
for $W_1,W_2 \in \ob \mathcal{C}$, and we will show that 
$W_1\hboxtr_{P(z)} W_2$ is an object of $\mathcal{C}$ if and only if 
$W_1\boxtimes_{P(z)} W_2$ exists in $\mathcal{C}$, in which case we
will have 
\[
W_1\boxtimes_{P(z)} W_2 = (W_1\hboxtr_{P(z)} W_2)'
\]
(observe the notation
\[
\boxtimes = \hboxtr',
\]
as in the special cases studied in \cite{tensor1}--\cite{tensor3}).
It is important to keep in mind that the space $W_1\hboxtr_{P(z)} W_2$
will depend on the category $\mathcal{C}$.}
\end{rema}

We formalize certain of the properties of the category $\mathcal{C}$ that we 
have been using, and some new ones, as follows:

\begin{assum}\label{assum-c}
Throughout the remainder of this work, we shall assume that
$\mathcal{C}$ is a full subcategory of the category $\mathcal{M}_{sg}$
or $\mathcal{G}\mathcal{M}_{sg}$ closed under the contragredient
functor (recall Notation \ref{MGM}; for now, we are not assuming that
$V\in \ob \mathcal{C}$). We shall also assume that $\mathcal{C}$ is
closed under taking finite direct sums.
\end{assum}

\begin{defi}\label{def-hboxtr}
{\rm For $W_{1}, W_{2}\in \ob \mathcal{C}$, define the subset 
\[
W_{1}\hboxtr_{P(z)}W_{2}\subset (W_{1}\otimes W_{2})^{*}
\]
of $(W_{1}\otimes W_{2})^{*}$ to be the union of the images
\[
I'(W')\subset (W_{1}\otimes W_{2})^{*}
\]
as $(W; I)$ ranges through all the $P(z)$-products of $W_{1}$ and
$W_{2}$ with $W\in \ob \mathcal{C}$. Equivalently,
$W_{1}\hboxtr_{P(z)}W_{2}$ is the union of the images $I'(W')$ as $W$
(or $W'$) ranges through $\ob \mathcal{C}$ and $I'$ ranges through
$\mathcal{N}_{W'}^{(W_{1}\otimes W_{2})^{*}}$---the space of
$\tilde{A}$-compatible linear maps
\[
W'\to (W_{1}\otimes W_{2})^{*}
\]
intertwining the actions of both 
\[
V\otimes \iota_{+}\C[t, t^{-1}, (z^{-1}-t)^{-1}]
\]
and $\mathfrak{s}\mathfrak{l}(2)$ on both spaces.}
\end{defi}

\begin{rema}\label{backslash=sumunion}
{\rm Since $\mathcal{C}$ is closed under finite direct sums (Assumption 
\ref{assum-c}), it is clear that $W_{1}\hboxtr_{P(z)}W_{2}$ is
in fact a linear subspace of $(W_{1}\otimes W_{2})^{*}$, and in particular,
it can be defined alternatively as the sum of all the images $I'(W')$:
\begin{equation}\label{hboxtr-sum}
W_{1}\hboxtr_{P(z)}W_{2}=\sum I'(W') = \bigcup I'(W')\subset 
(W_{1}\otimes W_{2})^{*},
\end{equation}
where the sum and union both range over $W\in \ob
\mathcal{C}$, $I\in \mathcal{M}_{W_{1}W_{2}}^{W}$.}
\end{rema}

For any generalized $V$-modules $W_{1}$ and $W_{2}$,
using the operator $L'_{P(z)}(0)$ (recall (\ref{LP'(j)}))
on $(W_{1}\otimes W_{2})^{*}$ we define
the generalized $L'_{P(z)}(0)$-eigenspaces 
$((W_{1}\otimes W_{2})^{*})_{[n]}$ for $n\in \C$ in the usual way:
\begin{equation}
((W_{1}\otimes W_{2})^{*})_{[n]}=\{w\in (W_{1}\otimes W_{2})^{*}\;|\;
(L'_{P(z)}(0)-n)^{m}w=0 \;{\rm for}\; m\in \N \;
\mbox{\rm sufficiently large}\}.
\end{equation}
Then we have the (proper) subspace 
\begin{equation}
\coprod_{n\in \C}((W_{1}\otimes W_{2})^{*})_{[n]}\subset 
(W_{1}\otimes W_{2})^{*}.
\end{equation}
We also define the ordinary $L'_{P(z)}(0)$-eigenspaces
$((W_{1}\otimes W_{2})^{*})_{(n)}$ in the usual way:
\begin{equation}
((W_{1}\otimes W_{2})^{*})_{(n)}=\{w\in (W_{1}\otimes W_{2})^{*}\;|\;
L'_{P(z)}(0)w=nw \}.
\end{equation}
Then we have the (proper) subspace
\begin{equation}
\coprod_{n\in \C}((W_{1}\otimes W_{2})^{*})_{(n)}\subset 
(W_{1}\otimes W_{2})^{*}.
\end{equation}

\begin{propo}\label{im:abc}
Let $W_{1}, W_{2}\in \ob \mathcal{C}$. 

(a)
The elements of $W_1\hboxtr_{P(z)} W_2$ are exactly the linear functionals
on $W_{1}\otimes W_{2}$ of the form $w'\circ
I(\cdot\otimes \cdot)$ for some $P(z)$-intertwining map $I$ of type
${W\choose W_1\,W_2}$ and some $w'\in W'$, $W\in\ob{\cal C}$.

(b) Let $(W; I)$ be any $P(z)$-product of $W_{1}$ and $W_{2}$, with 
$W$ any generalized $V$-module. Then for $n\in \C$,
\[
I'(W'_{[n]}) \subset  ((W_{1}\otimes W_{2})^{*})_{[n]}
\]
and 
\[
I'(W'_{(n)})\subset ((W_{1}\otimes W_{2})^{*})_{(n)}.
\]

(c) The structure $(W_1\hboxtr_{P(z)} W_2,Y'_{P(z)})$ (recall
(\ref{y'-p-z})) satisfies all the axioms in the definition of
(strongly $\tilde{A}$-graded) generalized $V$-module except perhaps
for the two grading conditions (\ref{set:dmltc}) and
(\ref{set:dmfin}).

(d) Suppose that the objects of the category $\mathcal{C}$ consist
only of (strongly $\tilde{A}$-graded) {\em ordinary}, as opposed to
{\em generalized}, $V$-modules. Then the structure $(W_1\hboxtr_{P(z)}
W_2,Y'_{P(z)})$ satisfies all the axioms in the definition of
(strongly $\tilde{A}$-graded ordinary) $V$-module except perhaps for
(\ref{set:dmltc}) and (\ref{set:dmfin}).
\end{propo}
\pf 
Part (a) is clear {}from the definition of $W_1\hboxtr_{P(z)} W_2$, and 
(b) follows {}from (\ref{j-lj}) with $j=0$.

For (c), let $(W; I)$ be any any $P(z)$-product of $W_{1}$ and $W_{2}$, with 
$W$ any generalized $V$-module. Then $(I'(W'), Y'_{P(z)})$ satisfies all
the conditions in the definition of (strongly $\tilde{A}$-graded) generalized 
$V$-module since $I'$ is $\tilde{A}$-compatible and intertwines the actions of 
$V\otimes {\mathbb C}[t,t^{-1}]$ and of ${\mathfrak s}{\mathfrak
l}(2)$; the $\C$-grading follows {}from Part (b). Note that
\begin{equation}\label{I'W'}
I':W' \rightarrow I'(W')
\end{equation}
is a map of generalized $V$-modules.  Since $W_1\hboxtr_{P(z)} W_2$
is the sum of these structures $I'(W')$ over $W\in \ob \mathcal{C}$ (recall
(\ref{hboxtr-sum})), we see that $(W_1\hboxtr_{P(z)} W_2, Y'_{P(z)})$ satisfies
all the conditions in the definition of generalized module except perhaps for 
(\ref{set:dmltc}) and (\ref{set:dmfin}).

Finally, Part (d) is proved by the same argument as for (c). In fact, for 
$(W; I)$ any $P(z)$-product of possibly generalized $V$-modules $W_{1}$ and $W_{2}$,
with $W$ any ordinary $V$-module, $(I'(W'), Y'_{P(z)})$ satisfies all the conditions 
in the definition of (strongly $\tilde{A}$-graded) ordinary $V$-module; the 
$\C$-grading (this time, by ordinary $L'_{P(z)}(0)$-eigenspaces) 
again follows {}from Part (b).
\epfv

\begin{rema}
{\rm Later, it will be convenient to introduce the term ``doubly
graded'' generalized module (or module) for the module structures in
Parts (c) and (d) (see Definition \ref{doublygraded} below).}
\end{rema}

The following terminology will be useful:

\begin{defi}\label{closedunderimages}
{\rm The category $\mathcal{C}$ is {\it closed under images} if, given
an object $M$ of $\mathcal{C}$, a generalized $V$-module $N$, and a
map $\phi:M \to N$ of generalized $V$-modules, the image $\phi(M)$ is
an object of $\mathcal{C}$.}
\end{defi}

Under this assumption, we have the following alternative description
of $W_1\hboxtr_{P(z)} W_2$:

\begin{propo}\label{backslash=union}
Suppose that $\mathcal{C}$ is closed under images (as well as under
contragredients and finite direct sums (Assumption \ref{assum-c})).
Let $W_1,W_2 \in \ob{\cal C}$.  Then the subspace $W_1\hboxtr_{P(z)}
W_2$ of $(W_{1}\otimes W_{2})^{*}$ is equal to the union and also to
the sum of the objects of $\mathcal{C}$ lying in $(W_{1}\otimes
W_{2})^{*}$:
\[
W_{1}\hboxtr_{P(z)}W_{2}=\bigcup W = \sum W \subset 
(W_{1}\otimes W_{2})^{*},
\]
where, in the union and in the sum, $W$ ranges through the subspaces
of $(W_{1}\otimes W_{2})^{*}$ that are objects of $\mathcal{C}$ when
equipped with the action $Y'_{P(z)}(\cdot,x)$ of $V$ and the
corresponding action of ${\mathfrak s}{\mathfrak l}(2)$ on
$(W_{1}\otimes W_{2})^{*}$.  In particular, every object of
$\mathcal{C}$ lying in $(W_{1}\otimes W_{2})^{*}$ is a subspace of
$W_{1}\hboxtr_{P(z)}W_{2}$ (and for this assertion, the assumption
that $\mathcal{C}$ is closed under images is not needed).
\end{propo}
\pf In view of Remark \ref{backslash=sumunion}, it suffices to show
that the subspaces $I'(W')$ of $(W_{1}\otimes W_{2})^{*}$, as $W$
ranges through $\ob{\cal C}$ and $I'$ ranges through the maps
indicated in Definition \ref{def-hboxtr}, coincide with the objects of
$\mathcal{C}$ lying in $(W_{1}\otimes W_{2})^{*}$.  But each $I'(W')
\in \ob{\cal C}$ because $\mathcal{C}$ is closed under images and
contragredients, and $I':W' \rightarrow I'(W')$ is a map of
generalized $V$-modules (recall (\ref{I'W'})).  Conversely, consider
an arbitrary object of $\mathcal{C}$ lying in $(W_{1}\otimes
W_{2})^{*}$ and write it as $W'$ for a suitable $W \in \ob{\cal C}$.
Then since the embedding
\[
J:W' \hookrightarrow (W_{1}\otimes W_{2})^{*}
\]
has the properties indicated in Notation \ref{scriptN}, it can be
written as $I'$ (with $I=J'$).  \epfv

The next result characterizes $W_{1}\boxtimes_{P(z)}W_{2}$, including
its existence, in terms of $W_1\hboxtr_{P(z)} W_2$; this result
generalizes Proposition 13.7 in \cite{tensor3}:

\begin{propo}\label{tensor1-13.7}
Let $W_{1}, W_{2}\in \ob \mathcal{C}$.  If $(W_1\hboxtr_{P(z)} W_2,
Y'_{P(z)})$ is an object of ${\cal C}$, denote by
$(W_1\boxtimes_{P(z)} W_2, Y_{P(z)})$ its contragredient module:
\[
W_1\boxtimes_{P(z)} W_2 = (W_1\hboxtr_{P(z)} W_2)'.
\]
Then the $P(z)$-tensor product of $W_{1}$ and $W_{2}$ in ${\cal C}$
exists and is
\[
(W_1\boxtimes_{P(z)} W_2, Y_{P(z)}; i'),
\]
where $i$ is
the natural inclusion {}from $W_1\hboxtr_{P(z)} W_2$ to $(W_1\otimes
W_2)^*$ (recall Notation \ref{scriptN}).  Conversely, let us assume
that $\mathcal{C}$ is closed under images.  If the $P(z)$-tensor
product of $W_1$ and $W_2$ in ${\cal C}$ exists, then
$(W_1\hboxtr_{P(z)} W_2, Y'_{P(z)})$ is an object of ${\cal C}$.
\end{propo}
\pf 
Suppose that $(W_1\hboxtr_{P(z)} W_2, Y'_{P(z)})$ is an object of
${\cal C}$ and take $(W_1\boxtimes_{P(z)} W_2, Y_{P(z)})$ and 
the map $i$ as indicated. Then 
\[
i\in \mathcal{N}_{W_1\hboxtr_{P(z)} W_2}^{(W_{1}\otimes W_{2})^{*}},
\]
and
\[
i'\in \mathcal{M}_{W_{1}W_{2}}^{W_1\boxtimes_{P(z)} W_2}.
\]
In the notation of Corollary \ref{tensorproductusingI'}, we take
$I_{0}=i'$, $I_{0}'=i$. For any pair $(W; I')$ as in Corollary
\ref{tensorproductusingI'}, we have
\[
I'(W')\subset W_1\hboxtr_{P(z)} W_2
\]
(which is the union of all such images), so that there is
certainly a unique module map
\[
\chi: W'\to W_1\hboxtr_{P(z)} W_2
\]
such that 
\[
i\circ \chi=I',
\]
namely, $I'$ itself, viewed as a module map. Thus by Corollary
\ref{tensorproductusingI'}, $W_1\boxtimes_{P(z)} W_2$ exists as
indicated.

Conversely, if the $P(z)$-tensor product of $W_1$ and $W_2$ in ${\cal
C}$ exists and is $(W_0; I_0)$, then for any $P(z)$-product $(W; I)$
with $W \in \ob \mathcal{C}$, we have a unique module map $\chi: W'\to
W'_0$ as in Corollary \ref{tensorproductusingI'} such that
$I'=I'_0\circ \chi$, so that $I'(W')\subset I'_0(W'_0)$, proving that
\[
W_1\hboxtr_{P(z)} W_2\subset I'_0(W'_0).
\]
On the other hand, $(W_0; I_0)$ is itself a $P(z)$-product, so that
\[
I'_0(W'_0)\subset W_1\hboxtr_{P(z)} W_2.
\]
Thus
\[
W_1\hboxtr_{P(z)} W_2=I'_0(W'_0),
\]
and so $W_1\hboxtr_{P(z)} W_2$ is a generalized
$V$-module and is the image of the module map
\[
I_{0}': W_{0}'\to W_1\hboxtr_{P(z)} W_2.
\]
Since $\mathcal{C}$ is closed under images by assumption, we have that
$W_1\hboxtr_{P(z)} W_2 \in \ob \mathcal{C}$.  \epfv

\begin{rema}
{\rm Suppose that $W_1\hboxtr_{P(z)} W_2$ is an object of
$\mathcal{C}$.  {}From Corollary \ref{tensorproductusingI'} and
Proposition \ref{tensor1-13.7} we see that
\begin{equation}\label{boxpair}
\langle\lambda, w_{(1)}\boxtimes_{P(z)}w_{(2)}\rangle
_{W_1\boxtimes_{P(z)} W_2}=
\lambda(w_{(1)}\otimes w_{(2)})
\end{equation}
for $\lambda\in W_1\hboxtr_{P(z)} W_2\subset (W_1\otimes W_2)^*$,
$w_{(1)}\in W_1$ and $w_{(2)}\in W_2$.}
\end{rema}

Our next goal is to present a crucial alternative description of the
subspace $W_1\hboxtr_{P(z)} W_2$ of $(W_1\otimes W_2)^*$. The main
ingredient of this description will be the ``$P(z)$-compatibility
condition,'' which was a cornerstone of the development of tensor
product theory in the special cases treated in
\cite{tensor1}--\cite{tensor3} and \cite{tensor4}.

Assume now that $W_{1}$ and $W_{2}$ are arbitrary generalized $V$-modules. 
Let $(W; I)$ ($W$ a generalized $V$-module) 
be a $P(z)$-product of $W_{1}$ and $W_{2}$ and let 
$w'\in W'$. Then {}from (\ref{j-tau}), Proposition \ref{productusingI'},
(\ref{tau-w-comp}), (\ref{3.7}) and (\ref{y'-p-z}), we have, for all
$v\in V$, 
\begin{eqnarray}\label{5.18-p}
\lefteqn{\tau_{P(z)}\left(x_0^{-1}\delta\bigg(\frac{x^{-1}_1-z}{x_0}
\bigg)
Y_{t}(v, x_1)\right)I'(w')}\nno\\
&&=I'\left(\tau_{W'}\left(x_0^{-1}\delta\bigg(\frac{x^{-1}_1-z}{x_0}
\bigg)
Y_{t}(v, x_1)\right)w'\right)\nno\\
&&=I'\left(x_0^{-1}\delta\bigg(\frac{x^{-1}_1-z}{x_0}
\bigg)Y_{W'}(v, x_1)w'\right)\nno\\
&&=x_0^{-1}\delta\bigg(\frac{x^{-1}_1-z}{x_0}
\bigg)I'(Y_{W'}(v,
x_1)w')\nno\\
&&=x_0^{-1}\delta\bigg(\frac{x^{-1}_1-z}{x_0}
\bigg)I'(\tau_{W'}(Y_{t}(v,
x_1))w')\nno\\
&&=x_0^{-1}\delta\bigg(\frac{x^{-1}_1-z}{x_0}
\bigg)\tau_{P(z)}(Y_{t}(v,
x_1))I'(w')\nn
&&=x_0^{-1}\delta\bigg(\frac{x^{-1}_1-z}{x_0}
\bigg)
Y'_{P(z)}(v, x_1)I'(w').
\end{eqnarray}
That is, 
$I'(w')$ satisfies the following nontrivial and subtle condition on
\[
\lambda \in (W_1\otimes W_2)^{*}:
\]

\begin{description}
\item{\bf The $P(z)$-compatibility condition}

(a) The {\em $P(z)$-lower truncation condition}: For all $v\in V$, the formal
Laurent series $Y'_{P(z)}(v, x)\lambda$ involves only finitely many
negative powers of $x$.

(b) The following formula holds:
\begin{eqnarray}\label{cpb}
\lefteqn{\tau_{P(z)}\bigg(x_0^{-1}\delta\bigg(\frac{x^{-1}_1-z}{x_0}
\bigg)
Y_{t}(v, x_1)\bigg)\lambda}\nno\\
&&=x_0^{-1}\delta\bigg(\frac{x^{-1}_1-z}{x_0}\bigg)
Y'_{P(z)}(v, x_1)\lambda  \;\;\mbox{ for all }\;v\in V.
\end{eqnarray}
(Note that the two sides of (\ref{cpb}) are not {\it a priori} equal
for general $\lambda\in (W_1\otimes W_2)^{*}$. Note also that Condition 
(a) insures that the right-hand side in Condition (b) is 
well defined.)
\end{description}

\begin{nota}
{\rm Note that the set of elements of $(W_1\otimes W_2)^*$ satisfying
either  the full $P(z)$-compatibility
condition or Part (a) of this condition forms a subspace. 
We shall denote the space of elements of $(W_1\otimes W_2)^*$ satisfying
the $P(z)$-compatibility
condition by
\[
\comp_{P(z)}((W_1\otimes W_2)^*).
\]}
\end{nota}

Recall {}from (\ref{W1W2beta}) that for each $\beta\in \tilde{A}$ we
have the subspace $((W_1\otimes W_2)^*)^{(\beta)}$ of $(W_1\otimes
W_2)^*$.  The sum of these subspaces is of course direct, and we
denote it by $((W_1\otimes W_2)^*)^{(\tilde{A})}$:
\[
((W_1\otimes W_2)^*)^{(\tilde{A})}
=\sum_{\beta\in \tilde{A}}((W_1\otimes W_2)^*)^{(\beta)}
=\coprod_{\beta\in \tilde{A}}((W_1\otimes W_2)^*)^{(\beta)}.
\]
Each space $((W_1\otimes W_2)^*)^{(\beta)}$ is $L'_{P(z)}(0)$-stable
(recall Proposition \ref{tau-a-comp} and Remark
\ref{L'jpreservesbetaspace}), so that we may consider the subspaces
\[
\coprod_{n\in \C}((W_1\otimes W_2)^*)_{[n]}^{(\beta)} \subset
((W_1\otimes W_2)^*)^{(\beta)}
\]
and 
\[
\coprod_{n\in \C}((W_1\otimes W_2)^*)_{(n)}^{(\beta)} \subset
((W_1\otimes W_2)^*)^{(\beta)}
\]
(recall Remark \ref{generalizedeigenspacedecomp}).  We now define the
two subspaces
\begin{equation}\label{W1W2_[C]^Atilde}
((W_1\otimes W_2)^*)_{[{\mathbb C}]}^{( \tilde A )}=
\coprod_{n\in
\C}\coprod_{\beta\in \tilde{A}}((W_1\otimes W_2)^*)_{[n]}^{(\beta)}
\subset ((W_1\otimes W_2)^*)^{(\tilde{A})}
\subset (W_1\otimes W_2)^*
\end{equation}
and 
\begin{equation}\label{W1W2_(C)^Atilde}
((W_1\otimes W_2)^*)_{({\mathbb C})}^{( \tilde A )}=
\coprod_{n\in
\C}\coprod_{\beta\in \tilde{A}}((W_1\otimes W_2)^*)_{(n)}^{(\beta)}
\subset ((W_1\otimes W_2)^*)^{(\tilde{A})}
\subset (W_1\otimes W_2)^*.
\end{equation}

\begin{rema}\label{singleanddoublegraded}
{\rm Any $L'_{P(z)}(0)$-stable subspace of $((W_1\otimes
W_2)^*)_{[{\mathbb C}]}^{( \tilde A )}$ is graded by generalized
eigenspaces (again recall Remark \ref{generalizedeigenspacedecomp}),
and if such a subspace is also $\tilde A$-graded, then it is doubly
graded; similarly for subspaces of $((W_1\otimes W_2)^*)_{({\mathbb
C})}^{( \tilde A )}$.}
\end{rema}

We have:

\begin{lemma}\label{a-tilde-comp}
Suppose that $\lambda\in ((W_1\otimes W_2)^*)^{(\tilde A)}$ satisfies
the $P(z)$-compatibility condition. Then every $\tilde{A}$-homogeneous
component of $\lambda$ also satisfies this condition.  In particular,
the space
\[
(\comp_{P(z)}((W_1\otimes W_2)^*)) \cap ((W_1\otimes W_2)^*)^{(\tilde A)}
\]
is $\tilde{A}$-graded.
\end{lemma}
\pf
When $v\in V$ is $\tilde{A}$-homogeneous, 
\[
\tau_{P(z)}\bigg(x_0^{-1}\delta\bigg(\frac{x^{-1}_1-z}{x_0} \bigg)
Y_{t}(v, x_1)\bigg)\;\;\mbox{ and }\;\;
x_0^{-1}\delta\bigg(\frac{x^{-1}_1-z}{x_0}\bigg) Y'_{P(z)}(v, x_1)
\]
are both $\tilde{A}$-homogeneous as operators, in the obvious sense.
By comparing the $\tilde{A}$-homogeneous components of both sides of
(\ref{cpb}), we see that the $\tilde{A}$-homogeneous components of
$\lambda$ also satisfy the $P(z)$-compatibility condition.  \epfv

\begin{rema}\label{stableundercomponentops}
{\rm Both the spaces $((W_1\otimes W_2)^*)_{[{\mathbb C}]}^{( \tilde A
)}$ and $((W_1\otimes W_2)^*)_{({\mathbb C})}^{( \tilde A )}$ are
stable under the component operators $\tau_{P(z)}(v\otimes t^m)$ of
the operators $Y'_{P(z)}(v,x)$ for $v\in V$, $m\in {\mathbb Z}$, and
under the operators $L'_{P(z)}(-1)$, $L'_{P(z)}(0)$ and
$L'_{P(z)}(1)$.  For the $\tilde A$-grading, this follows {}from
Proposition \ref{tau-a-comp} and Remark \ref{L'jpreservesbetaspace},
and for the $\C$-gradings, we simply follow the proof of Proposition
\ref{gweight}, using Propositions \ref{id-dev} and \ref{pz-comm}
together with (\ref{pz-sl-2-pz-y--1}).}
\end{rema}

Again let $(W; I)$ ($W$ a generalized $V$-module) 
be a $P(z)$-product of $W_{1}$ and $W_{2}$ and let 
$w'\in W'$. 
Since $I'$ in particular intertwines the actions of
$V\otimes{\mathbb C}[t, t^{-1}]$ and of $\mathfrak{s}\mathfrak{l}(2)$,
and is $\tilde{A}$-compatible, $I'(W')$ is a generalized $V$-module,
as we have seen in the proof of Proposition \ref{im:abc}.
Therefore, for each $w'\in W'$, $I'(w')$ also satisfies
the following condition on 
\[
\lambda \in (W_1\otimes W_2)^*:
\]
\begin{description}
\item{\bf The $P(z)$-local grading restriction condition}

(a) The {\em $P(z)$-grading condition}: $\lambda$ is a (finite) sum of
generalized eigenvectors  for the operator
$L'_{P(z)}(0)$ on $(W_1\otimes W_2)^*$ that 
are also homogeneous with respect to $\tilde A$, that is, 
\[
\lambda\in ((W_1\otimes W_2)^*)_{[{\mathbb C}]}^{( \tilde A )}.
\]
\label{homo}

(b) Let $W_{\lambda}$ be the smallest doubly graded (or equivalently,
$\tilde A$-graded; recall Remark \ref{singleanddoublegraded}) subspace
of $((W_1\otimes W_2)^*)_{[ {\mathbb C} ]}^{( \tilde A )}$ containing
$\lambda$ and stable under the component operators
$\tau_{P(z)}(v\otimes t^m)$ of the operators $Y'_{P(z)}(v,x)$ for
$v\in V$, $m\in {\mathbb Z}$, and under the operators $L'_{P(z)}(-1)$,
$L'_{P(z)}(0)$ and $L'_{P(z)}(1)$.  (In view of Remark
\ref{stableundercomponentops}, $W_{\lambda}$ indeed exists.)  Then
$W_{\lambda}$ has the properties
\begin{eqnarray}
&\dim(W_{\lambda})^{(\beta)}_{[n]}<\infty,&\label{lgrc1}\\
&(W_{\lambda})^{(\beta)}_{[n+k]}=0\;\;\mbox{ for }\;k\in {\mathbb Z}
\;\mbox{ sufficiently negative},&\label{lgrc2}
\end{eqnarray}
for any $n\in {\mathbb C}$ and $\beta\in \tilde A$, where as usual the
subscripts denote the ${\mathbb C}$-grading and the superscripts
denote the $\tilde A$-grading.
\end{description}

In the case that $W$ is an (ordinary) $V$-module and $w'\in W'$,
$I'(w')$ also satisfies the following $L(0)$-semisimple version of
this condition on $\lambda \in (W_1\otimes W_2)^*$:

\begin{description}
\item{\bf The $L(0)$-semisimple $P(z)$-local grading restriction condition}

(a) The {\em $L(0)$-semisimple 
$P(z)$-grading condition}: $\lambda$ is a (finite) sum of
eigenvectors  for the operator
$L'_{P(z)}(0)$ on $(W_1\otimes W_2)^*$ that 
are also homogeneous with respect to $\tilde A$, that is,
\[
\lambda\in ((W_1\otimes W_2)^*)_{({\mathbb C})}^{( \tilde A )}.
\]
\label{semi-homo}

(b) Consider $W_\lambda$ as above, which in this case is in fact the
smallest doubly graded (or equivalently, $\tilde A$-graded) subspace
of $((W_1\otimes W_2)^*)_{({\mathbb C})}^{( \tilde A )}$ containing
$\lambda$ and stable under the component operators
$\tau_{P(z)}(v\otimes t^m)$ of the operators $Y'_{P(z)}(v,x)$ for
$v\in V$, $m\in {\mathbb Z}$, and under the operators $L'_{P(z)}(-1)$,
$L'_{P(z)}(0)$ and $L'_{P(z)}(1)$.  Then $W_\lambda$ has the
properties
\begin{eqnarray}
&\dim(W_\lambda)^{(\beta)}_{(n)}<\infty,&\label{semi-lgrc1}\\
&(W_\lambda)^{(\beta)}_{(n+k)}=0\;\;\mbox{ for }\;k\in {\mathbb Z}
\;\mbox{ sufficiently negative},&\label{semi-lgrc2}
\end{eqnarray}
for any $n\in {\mathbb C}$ and $\beta\in \tilde A$, where  the
subscripts denote the ${\mathbb C}$-grading and the superscripts
denote the $\tilde A$-grading.

\end{description}

\begin{nota}
{\rm Note that the set of elements of $(W_1\otimes W_2)^*$ satisfying
either of these two $P(z)$-local grading restriction conditions, or
either of the Part (a)'s in these conditions, forms a subspace.  We
shall denote the space of elements of $(W_1\otimes W_2)^*$ satisfying
the $P(z)$-local grading restriction condition and the
$L(0)$-semisimple $P(z)$-local grading restriction condition by
\[
\lgr_{[\C]; P(z)}((W_1\otimes W_2)^*)
\]
and 
\[
\lgr_{(\C); P(z)}((W_1\otimes W_2)^*),
\]
respectively.}
\end{nota}

The following theorems are among the most important in this work.
Note that even in the finitely reductive case studied in
\cite{tensor3}, they are stronger and more general than (the last
assertion of) Theorem 13.9 in \cite{tensor3}.  The proofs of these two
theorems will be given in the next section.

\begin{theo}\label{comp=>jcb}
Let $\lambda$ be an element of $(W_{1}\otimes W_{2})^{*}$ satisfying
the $P(z)$-compatibility condition. Then when acting on $\lambda$, the
Jacobi identity for $Y'_{P(z)}$ holds, that is,
\begin{eqnarray}
\lefteqn{x_{0}^{-1}\delta
\left({\displaystyle\frac{x_{1}-x_{2}}{x_{0}}}\right)Y'_{P(z)}(u, x_{1})
Y'_{P(z)}(v, x_{2})\lambda}\nno\\
&&\hspace{2ex}-x_{0}^{-1} \delta
\left({\displaystyle\frac{x_{2}-x_{1}}{-x_{0}}}\right)Y'_{P(z)}(v, x_{2})
Y'_{P(z)}(u, x_{1})\lambda\nonumber \\
&&=x_{2}^{-1} \delta
\left({\displaystyle\frac{x_{1}-x_{0}}{x_{2}}}\right)Y'_{P(z)}(Y(u, x_{0})v,
x_{2})\lambda\label{cjcb}
\end{eqnarray}
for $u, v\in V$.
\end{theo}

\begin{theo}\label{stable}
The subspace $\comp_{P(z)}((W_1\otimes W_2)^*)$ of $(W_{1}\otimes
W_{2})^{*}$ is stable under the operators $\tau_{P(z)}(v\otimes
t^{n})$ for $v\in V$ and $n\in {\mathbb Z}$, and in the M\"obius case,
also under the operators $L'_{P(z)}(-1)$, $L'_{P(z)}(0)$ and
$L'_{P(z)}(1)$; similarly for the subspaces $\lgr_{[\C];
P(z)}((W_1\otimes W_2)^*)$ and $\lgr_{(\C); P(z)}((W_1\otimes W_2)^*)$.
\end{theo}

\begin{rema}{\rm
The converse of Theorem \ref{comp=>jcb} is not true. One can see this
in the tensor product theory of the ``trivial'' case where $V$ is a
vertex operator algebra associated with a finite-dimensional unital
commutative associative algebra $(A,\cdot,1)$ with derivation $D=0$
(cf.\ Remark \ref{va>cva}). In this case, $(V,Y,{\bf 1},\omega)=
(A,\cdot,1,0)$ and the Jacobi identity for a $V$-module $W$ reduces to
\begin{eqnarray*}
&{\dps x^{-1}_0\delta \bigg( {x_1-x_2\over x_0}\bigg) u\cdot (v\cdot
w) - x^{-1}_0\delta \bigg( {x_2-x_1\over -x_0}\bigg) v\cdot (u\cdot
w) }&\nno\\
&{\dps = x^{-1}_2\delta \bigg( {x_1-x_0\over x_2}\bigg)
(u\cdot v)\cdot w}
\end{eqnarray*}
for $u,v\in A$ and $w\in W$, where we also use ``$\cdot$'' to denote
the action of $A$ on its modules.  In particular, a $V$-module is just
a finite-dimensional module for the associative algebra $A$. Given
$V$-modules $W_1$ and $W_2$, the action $Y'_{P(z)}$ given in
(\ref{Y'def}) now becomes
\[
(Y'_{P(z)}(v,x)\lambda)(w_{(1)}\otimes w_{(2)})=
\lambda(w_{(1)}\otimes v\cdot w_{(2)}).
\]
{}From this it is clear that (\ref{cjcb}) holds for any element
$\lambda\in (W_1\otimes W_2)^*$. However, the $P(z)$-compatibility
condition (\ref{cpb}) in this case reduces to
\[
\lambda(v\cdot w_{(1)}\otimes w_{(2)})=
\lambda(w_{(1)}\otimes v\cdot w_{(2)})
\]
for all $v\in A$, $w_{(1)}\in W_1$ and $w_{(2)}\in W_2$, which is not
necessarily true for every $\lambda$.  This example is discussed
further in Remark 2.20 of \cite{HLLZ}, which treats a range of issues
related to the compatibility condition, intertwining operators, and
tensor product theory.}
\end{rema}

We now generalize the notion of ``weak module'' for a vertex operator 
algebra to our M\"obius or
conformal vertex algebra $V$:

\begin{defi}
{\rm A {\it weak module for $V$} (or {\it weak $V$-module}) is a
vector space $W$ equipped with a vertex operator map
\[
Y_{W}: V\otimes W\to W[[x,x^{-1}]]
\]
satisfying (only) the axioms (\ref{ltc-w}),
(\ref{m-1left}), (\ref{m-Jacobi}) and (\ref{L-1}) in Definition
\ref{cvamodule} (note that there is no grading given on $W$) and in
case $V$ is M\"obius, also the existence of a representation of
$\mathfrak{s}\mathfrak{l}(2)$ on $W$, as in Definition
\ref{moduleMobius}, satisfying the conditions
(\ref{sl2-1})--(\ref{sl2-3}).  }
\end{defi}

Then
we have:

\begin{theo}\label{wk-mod}
The space $\comp_{P(z)}((W_1\otimes W_2)^*)$, equipped with the vertex
operator map $Y'_{P(z)}$ and, in case $V$ is M\"{o}bius, also equipped
with the operators $L_{P(z)}'(-1)$, $L_{P(z)}'(0)$ and $L_{P(z)}'(1)$,
is a weak $V$-module; similarly for the spaces
\begin{equation}\label{COMPintLGR[]}
(\comp_{P(z)}((W_1\otimes W_2)^*))\cap 
(\lgr_{[\C]; P(z)}((W_1\otimes W_2)^*))
\end{equation}
and 
\begin{equation}\label{COMPintLGR()}
(\comp_{P(z)}((W_1\otimes W_2)^*))\cap 
(\lgr_{(\C); P(z)}((W_1\otimes W_2)^*)).
\end{equation}
\end{theo}
\pf By Theorem \ref{stable}, $Y'_{P(z)}$ is a map {}from the tensor
product of $V$ with any of these three subspaces to the space of
formal Laurent series with elements of the subspace as coefficients.
By Proposition \ref{id-dev} and Theorem \ref{comp=>jcb} and, in the
case that $V$ is a M\"{o}bius vertex algebra, also by Propositions
\ref{sl-2} and \ref{pz-l-y-comm}, we see that all the axioms for weak
$V$-module are satisfied.  \epfv

We also have:

\begin{theo}\label{generation}
Let 
\[
\lambda\in 
\comp_{P(z)}((W_1\otimes W_2)^*)\cap \lgr_{[\C]; P(z)}((W_1\otimes W_2)^*).
\] 
Then $W_{\lambda}$ (recall Part (b) of the $P(z)$-local grading
restriction condition) equipped with the vertex operator map
$Y_{P(z)}'$ and, in case $V$ is M\"{o}bius, also equipped with the
operators $L'_{P(z)}(-1)$, $L'_{P(z)}(0)$ and $L'_{P(z)}(1)$, is a
(strongly-graded) generalized $V$-module. If in addition
\[
\lambda\in 
\comp_{P(z)}((W_1\otimes W_2)^*)\cap \lgr_{(\C); P(z)}((W_1\otimes W_2)^*),
\] 
that is, $\lambda$ is a sum of eigenvectors of $L_{P(z)}'(0)$, then
$W_{\lambda}$ ($\subset ((W_{1}\otimes
W_{2})^{*})_{(\C)}^{(\tilde{A})}$) is a (strongly-graded) $V$-module.
\end{theo}
\pf Decompose $\lambda$ as
\[
\lambda = \sum_{\beta\in \tilde{A}}\lambda^{(\beta)}
\]
(finite sum), where
\[
\lambda^{(\beta)} \in ((W_1\otimes W_2)^*)^{(\beta)}.
\]
By Lemma \ref{a-tilde-comp}, each
$\lambda^{(\beta)}$ satisfies the $P(z)$-compatibility condition.
Also, each $\lambda^{(\beta)}$ satisfies the $P(z)$-grading condition
(and in the semisimple case, the $L(0)$-semisimple $P(z)$-grading
condition), and each $W_{{\lambda}^{(\beta)}}$ is simply the smallest
subspace containing $\lambda^{(\beta)}$ and stable under the operators
listed above (without the $\tilde A$-gradedness condition).  Moreover,
each $W_{{\lambda}^{(\beta)}} \subset W_{\lambda}$ and in fact
\[
W_{\lambda} = \sum_{\beta\in \tilde{A}}W_{{\lambda}^{(\beta)}}.
\] 
Thus each $\lambda^{(\beta)}$ lies in the space (\ref{COMPintLGR[]})
(or (\ref{COMPintLGR()})).  (Note that we have reduced Theorem
\ref{generation} to the $\tilde A$-homogeneous case.)  By Theorem
\ref{wk-mod}, each $W_{{\lambda}^{(\beta)}}$ is a weak submodule of
the weak module (\ref{COMPintLGR[]}) (or (\ref{COMPintLGR()})), and
hence is a (strongly-graded) generalized module (or module).  Thus
$W_{\lambda}$ has the same properties.
\epfv

Now we can give an alternative description of $W_1\hboxtr_{P(z)} W_2$
by characterizing the elements of $W_1\hboxtr_{P(z)} W_2$ using the
$P(z)$-compatibility condition and the $P(z)$-local grading
restriction conditions, generalizing Theorem 13.10 in
\cite{tensor3}. This description combines the results we have just
presented; these results will be crucial in later sections, especially
in the construction of the associativity isomorphisms, specifically,
in the proof of Theorem \ref{9.7-1} below.

\begin{theo}\label{characterizationofbackslash}
Suppose that for every element 
\[
\lambda
\in \comp_{P(z)}((W_1\otimes W_2)^*)\cap \lgr_{[\C]; P(z)}((W_1\otimes W_2)^*)
\]
the (strongly-graded) generalized module $W_{\lambda}$ given in
Theorem \ref{generation} is a generalized submodule of some object of
${\cal C}$ included in $(W_1\otimes W_2)^*$ (this of course holds in
particular if ${\cal C} = {\cal GM}_{sg}$). Then
\[
W_1\hboxtr_{P(z)}W_2=\comp_{P(z)}((W_1\otimes W_2)^*)
\cap \lgr_{[\C]; P(z)}((W_1\otimes W_2)^*).
\]
Suppose that ${\cal C}$ is a category of strongly-graded $V$-modules
(that is, ${\cal C}\subset {\cal M}_{sg}$) and that for every element
\[
\lambda
\in \comp_{P(z)}((W_1\otimes W_2)^*)\cap \lgr_{(\C); P(z)}((W_1\otimes W_2)^*)
\]
the (strongly-graded) $V$-module $W_{\lambda}$ given in Theorem
\ref{generation} is a submodule of some object of ${\cal C}$ included
in $(W_1\otimes W_2)^*$ (which of course holds in particular if ${\cal
C} = {\cal M}_{sg}$).  Then
\[
W_1\hboxtr_{P(z)}W_2=\comp_{P(z)}((W_1\otimes W_2)^*)
\cap \lgr_{(\C); P(z)}((W_1\otimes W_2)^*).
\]
\end{theo}
\pf 
We have seen that
\[
W_1\hboxtr_{P(z)}W_2\subset 
\comp_{P(z)}((W_1\otimes W_2)^*)\cap \lgr_{[\C]; P(z)}((W_1\otimes W_2)^*)
\] 
and, in case ${\cal C}\subset {\cal M}_{sg}$,
\[
W_1\hboxtr_{P(z)}W_2\subset 
\comp_{P(z)}((W_1\otimes W_2)^*)\cap \lgr_{(\C); P(z)}((W_1\otimes W_2)^*).
\] 
On the other hand, by the assumptions, every element $\lambda$ of
\[
\comp_{P(z)}((W_1\otimes W_2)^*)\cap \lgr_{[\C]; P(z)}((W_1\otimes W_2)^*)
\] 
and, in case ${\cal C}\subset {\cal M}_{sg}$, every element $\lambda$
of
\[
\comp_{P(z)}((W_1\otimes W_2)^*)\cap \lgr_{(\C); P(z)}((W_1\otimes W_2)^*),
\] 
is contained in $W_{\lambda}$ and thus in some object of ${\cal C}$,
and for any such (generalized) module, the inclusion map into
$(W_1\otimes W_2)^*$ satisfies the intertwining conditions in
Proposition \ref{pz}.  Thus $\lambda$ lies in $W_1\hboxtr_{P(z)}W_2$,
proving the desired inclusion.  \epfv

\subsection{Constructions of $Q(z)$-tensor products}

We now give the construction of $Q(z)$-tensor products.  It is
analogous to that of $P(z)$-tensor products, and the formulations,
results and proofs in this section largely parallel those in
Section 5.2.  As usual, $z\in \C^{\times}$.

Given generalized $V$-modules $W_1$ and $W_2$, we shall be
constructing an action of the space $V \otimes \iota_{+}{\mathbb
C}[t,t^{- 1},(z+t)^{-1}]$ on the space $(W_1 \otimes W_2)^*$.

Let $I$ be a $Q(z)$-intertwining map of type ${W_3\choose W_1\, W_2}$,
as in Definition \ref{im:qimdef}. Consider the contragredient
generalized $V$-module $(W_{3}', Y_{3}')$, recall the opposite vertex
operator (\ref{yo}) and formula (\ref{y'}), and recall why the
ingredients of formula (\ref{imq:def}) are well defined. For $v\in V$,
$w_{(1)}\in W_{1}$, $w_{(2)}\in W_{2}$ and $w'_{(3)}\in W'_3$,
applying $w'_{(3)}$ to (\ref{imq:def}) we obtain
\begin{eqnarray}\label{imq:def'}
\lefteqn{\left\langle z^{-1}\delta\bigg(\frac{x_1-x_0}{z}\bigg)
Y_{3}'(v, x_0)w'_{(3)}, I(w_{(1)}\otimes w_{(2)})\right\rangle}\nno\\
&&=\left\langle w'_{(3)},x_{0}^{-1}\delta\bigg(\frac{x_1-z}{x_0}\bigg)
I(Y_1^{o}(v, x_1)w_{(1)}\otimes
w_{(2)})\right\rangle\nno\\
&&\quad -\left\langle w'_{(3)}, x^{-1}_0\delta\bigg(\frac{z-x_1}{-x_0}\bigg)
I(w_{(1)}\otimes Y_2(v, x_1)w_{(2)})\right\rangle.
\end{eqnarray}
We shall use this to motivate our action.

As we discussed in Section 5.1 (see (\ref{3.12}) and
(\ref{3.13})), in the left-hand side of (\ref{imq:def'}), the
coefficients of
\begin{equation}\label{qdeltaY3'}
z^{-1}\delta\bigg(\frac{x_1-x_0}{z}\bigg) Y_{3}'(v, x_0)
\end{equation}
in powers of $x_0$ and $x_1$, for all $v\in V$, span
\begin{equation}\label{qtausubW3'}
\tau_{W'_3}(V
\otimes \iota_{+}{\mathbb C}[t,t^{- 1},(z+t)^{-1}])
\end{equation}
(recall (\ref{tauw}) and (\ref{3.7})).  We now define a linear action
of $V \otimes \iota_{+}{\mathbb C}[t,t^{- 1},(z+t)^{-1}]$ on $(W_1
\otimes W_2)^*$, that is, a linear map $$\tau_{Q(z)}: V\otimes
\iota_{+}{\mathbb C}[t, t^{-1}, (z+t)^{-1}]\to{\rm End}\;(W_1\otimes
W_2)^{*}.$$ Recall the notations $T_{-z}^{+}$ and $T_{-z}^{o}$ {}from
Section 5.1 ((\ref{Tpm-z}) and (\ref{To-z})).

\begin{defi}\label{deftauQ}
{\rm We define the linear action $\tau_{Q(z)}$ of 
\[
V
\otimes \iota_{+}{\mathbb C}[t,t^{- 1},(z+t)^{-1}]
\]
on $(W_1 \otimes
W_2)^*$
by
\begin{equation}\label{(5.1)}
(\tau_{Q(z)}(\xi)\lambda)(w_{(1)}\otimes w_{(2)})
=\lambda(\tau_{W_1}(T_{-z}^{o}\xi)w_{(1)}\otimes w_{(2)})
-\lambda(w_{(1)}\otimes \tau_{W_2}(T_{-z}^{+}\xi)w_{(2)})
\end{equation}
for $\xi\in V\otimes \iota_{+}{\mathbb C}[t, t^{-1}, (z+t)^{-1}]$,
$\lambda\in (W_1\otimes W_2)^{*}$, $w_{(1)}\in W_1$, $w_{(2)}\in W_2$,
and denote by $Y'_{Q(z)}$ the action of $V\otimes{\mathbb C}[t,t^{-1}]$
on $(W_1\otimes W_2)^*$ thus defined, that is,
\begin{equation}\label{y'-q-z}
Y'_{Q(z)}(v, x)=\tau_{Q(z)}(Y_{t}(v, x))
\end{equation}
for $v\in V$.}
\end{defi}

Using Lemma \ref{lemma5.2}, (\ref{3.7}) and (\ref{tauw-yto}),
we see
that (\ref{(5.1)}) can be written using generating
functions as
\begin{eqnarray}\label{5.2}
\lefteqn{\left(\tau_{Q(z)}
\left(z^{-1}\delta\left(\frac{x_1-x_0}{z}\right) Y_{t}(v,
x_0)\right)\lambda\right)(w_{(1)}\otimes w_{(2)})}\nno\\
&&=x^{-1}_0\delta\left(\frac{x_1-z}{x_0}\right) \lambda(Y_1^{o}(v,
x_1)w_{(1)}\otimes w_{(2)})\nno\\
&&\hspace{2em}-x_0^{-1}\delta\left(\frac{z-x_1}{-x_0}\right)
\lambda(w_{(1)}\otimes Y_2(v, x_1)w_{(2)})
\end{eqnarray}
for $v\in V$, $\lambda\in (W_1\otimes W_2)^{*}$, $w_{(1)}\in W_1$,
$w_{(2)}\in W_2$; compare this with (\ref{imq:def'}).  The generating
function form of the action $Y'_{Q(z)}$ can be obtained by taking
$\res_{x_1}$ of both sides of (\ref{5.2}):
\begin{eqnarray}\label{Y'qdef}
\lefteqn{(Y'_{Q(z)}(v,x_0)\lambda)(w_{(1)} \otimes
w_{(2)})}\nno\\
&&= \res_{x_1} x^{-1}_0 \delta\left(\frac{x_1-z}{x_0}\right)
\lambda(Y^o_1(v,x_1)w_{(1)} \otimes w_{(2)})\nno\\
&&\quad-
\res_{x_1}x^{-1}_0 \delta \left(\frac{z-x_1}{-x_0}\right)\lambda(w_{(1)}
\otimes Y_2(v,x_1)w_{(2)})\nno\\
&&= \lambda(Y^o_1(v,x_0 + z)w_{(1)}
\otimes w_{(2)}) \nno\\
&&\quad -
\res_{x_1}x^{-1}_0 \delta \left(\frac{z-x_1}{-x_0}\right)\lambda(w_{(1)}
\otimes Y_2(v,x_1)w_{(2)}).
\end{eqnarray}

\begin{rema}\label{I-intw-q}{\rm
Using the actions $\tau_{W'_3}$ and $\tau_{Q(z)}$, we can write
(\ref{imq:def'}) as
\[
\left(z^{-1}\delta\left(\frac{x_1-x_0}{z}\right)
Y_{3}'(v, x_0)w'_{(3)}\right)\circ I=
\tau_{Q(z)}\left(z^{-1}\delta\left(\frac{x_1-x_0}{z}\right)
Y_t(v, x_0)\right)(w'_{(3)}\circ I)
\]
or equivalently, as
\[
\left(\tau_{W'_3}\left(z^{-1}\delta\left(\frac{x_1-x_0}{z}\right)
Y_t(v, x_0)\right)w'_{(3)}\right)\circ I=
\tau_{Q(z)}\left(z^{-1}\delta\left(\frac{x_1-x_0}{z}\right)
Y_t(v, x_0)\right)(w'_{(3)}\circ I).
\]
}
\end{rema}

Recall the $\tilde{A}$-grading on $(W_1\otimes W_2)^{*}$ and the $A$-grading 
on $V
\otimes \iota_{+}{\mathbb C}[t,t^{- 1},(z^{-1}-t)^{-1}]$. 
Similarly, we also have an $A$-grading 
on $V
\otimes \iota_{+}{\mathbb C}[t,t^{- 1},(z+t)^{-1}]$. 
Definition \ref{linearactioncompatible} also 
applies to a linear action of 
$V \otimes \iota_{+}{\mathbb C}[t,t^{- 1}, (z+t)^{-1}]$
on $(W_1 \otimes W_2)^*$.
{}From 
(\ref{(5.1)}) or (\ref{5.2}), we have:

\begin{propo}\label{tau-q-a-comp}
The action $\tau_{Q(z)}$ is $\tilde{A}$-compatible. \epf
\end{propo}

We also have:

\begin{propo}\label{5.1}
The action $Y'_{Q(z)}$ has the property
\begin{equation}\label{Q-id}
Y'_{Q(z)}({\bf 1}, x)=1
\end{equation}
and the $L(-1)$-derivative property
\begin{equation}\label{QL-1}
\frac{d}{dx}Y'_{Q(z)}(v, x)=Y'_{Q(z)}(L(-1)v, x)
\end{equation}
for $v\in V$. 
\end{propo}
\pf
{From} (\ref{Y'qdef}), (\ref{yo}) and (\ref{3termdeltarelation}), 
\begin{eqnarray}
(Y'_{Q(z)}({\bf 1},x)\lambda)(w_{(1)} \otimes
w_{(2)}) 
&=& \res_{x_1}x^{-1} \delta\left(\frac{x_1-
z}{x}\right)\lambda(w_{(1)} \otimes w_{(2)})\nno\\
&&-\res_{x_1}x^{-1}
\delta\left(\frac{z-x_1}{-x}\right)\lambda(w_{(1)} \otimes w_{(2)})\nno\\
&=&
\res_{x_1}x^{-1}_1 \delta\left(\frac{z+x}{x_1}\right)\lambda(w_{(1)}
\otimes w_{(2)})\nno\\
&=& \lambda(w_{(1)} \otimes
w_{(2)}),
\end{eqnarray}
proving (\ref{Q-id}).  To prove the $L(-1)$-derivative property,
observe that {from} (\ref{Y'qdef}),
\begin{eqnarray}\label{5.8}
\lefteqn{\left(\left(\frac{d}{dx}
Y'_{Q(z)}(v,x)\right)\lambda\right)(w_{(1)} \otimes w_{(2)})}\nno\\
&&= \frac{d}{dx}
\lambda(Y^o_1(v,x + z)w_{(1)} \otimes w_{(2)})\nno\\
&&\quad
+\res_{x_1}\left(\frac{d}{dx} z^{-1}\delta\left(\frac{-x + x_1}{z}\right)
\right)\lambda(w_{(1)}
\otimes Y_2(v,x_1)w_{(2)}).
\end{eqnarray} 
But for any formal Laurent series $f(x)$, we have
\begin{equation}
\frac{d}{dx}f\left(\frac{-x+x_{1}}{z}\right)
=-\frac{d}{dx_{1}}f\left(\frac{-x+x_{1}}{z}\right)
\end{equation}
and if $f(x)$ involves only finitely many negative powers of $x$,
\begin{equation}
\res_{x_{1}}\left(\frac{d}{dx_{1}}z^{-1}\delta\left(
\frac{-x+x_{1}}{z}\right)\right)f(x_{1})=
-\res_{x_{1}}z^{-1}\delta\left(
\frac{-x+x_{1}}{z}\right)\frac{d}{dx_{1}}f(x_{1})
\end{equation}
(since the residue of a derivative is $0$).
We also have the $L(-1)$-derivative property (\ref{yo-l-1}) 
for $Y^o$.
Thus the right-hand side of (\ref{5.8}) equals
\begin{eqnarray}
\lefteqn{\lambda(Y^o_1(L(-
1)v,x+z)w_{(1)} \otimes w_{(2)})}\nno\\
&&\quad +
\res_{x_1}z^{-1}\delta\left(
\frac{-x+x_{1}}{z}\right)\frac{d}{dx_1}\lambda(w_{(1)} \otimes
Y_2(v,x_1)w_{(2)})\nno\\
&&=  \lambda(Y^o_1(L(-1)v,x+z)w_{(1)}
\otimes w_{(2)})\nno\\
&&\quad + \res_{x_1}z^{-1}\delta\left(
\frac{-x+x_{1}}{z}\right)\lambda(w_{(1)}
\otimes Y_2(L(-1)v,x_1)w_{(2)})\nno\\
&&=  (Y'_{Q(z)}(L(-1)v,x)\lambda)(w_{(1)}
\otimes w_{(2)}),
\end{eqnarray}
proving (\ref{QL-1}).  \epfv

\begin{propo}\label{qz-comm}
The action $Y'_{Q(z)}$ satisfies the commutator formula for vertex
operators: On $(W_1\otimes W_2)^{*}$,
\begin{eqnarray}\label{commu-q-z}
\lefteqn{[Y'_{Q(z)}(v_1, x_1), Y'_{Q(z)}(v_2, x_2)]}\nn
&&=\res_{x_0}x_2^{-1}\delta\left(\frac{x_1-x_0}{x_2}\right)
Y'_{Q(z)}(Y(v_1, x_0)v_2, x_2)
\end{eqnarray}
for $v_1, v_2\in V$.
\end{propo}
\pf As usual, the reader should note the well-definedness of each
expression and the justifiability of each use of a $\delta$-function
property in the argument that follows.  This argument is the same as
the proof of Proposition 5.2 of \cite{tensor1}, given in Section 8 of
\cite{tensor2}.  Let $\lambda \in (W_{1}\otimes W_{2})^{*}$, $v_{1},
v_{2}\in V$, $w_{(1)}\in W_{1}$ and $w_{(2)}\in W_{2}$. By
(\ref{Y'qdef}),
\bea\label{8.1}
\lefteqn{(Y'_{Q(z)}(v_{1}, x_{1})Y'_{Q(z)}(v_{2}, x_{2})
\lambda)(w_{(1)}\otimes w_{(2)})}\nno\\
&&=\mbox{\rm Res}_{y_{1}}x_{1}^{-1}\delta
\left(\frac{y_{1}-z}{x_{1}}\right)(Y'_{Q(z)}(v_{2}, x_{2})
\lambda)(Y_{1}^{o}(v_{1}, y_{1})w_{(1)}\otimes w_{(2)})\nno\\
&&\quad -\mbox{\rm Res}_{y_{1}}x_{1}^{-1}\delta
\left(\frac{z-y_{1}}{-x_{1}}\right)(Y'_{Q(z)}(v_{2}, x_{2})
\lambda)(w_{(1)}\otimes Y_{2}(v_{1}, y_{1})w_{(2)})\nno\\
&&=\mbox{\rm Res}_{y_{1}}\mbox{\rm Res}_{y_{2}}x_{1}^{-1}\delta
\left(\frac{y_{1}-z}{x_{1}}\right)x_{2}^{-1}\delta
\left(\frac{y_{2}-z}{x_{2}}\right)\cdot \nno\\
&&\hspace{4em}\cdot \lambda(Y_{1}^{o}(v_{2}, y_{2})Y_{1}^{o}
(v_{1}, y_{1})w_{(1)}\otimes w_{(2)})
\nno\\
&&\quad -\mbox{\rm Res}_{y_{1}}\mbox{\rm Res}_{y_{2}}x_{1}^{-1}\delta
\left(\frac{y_{1}-z}{x_{1}}\right)x_{2}^{-1}\delta
\left(\frac{z-y_{2}}{-x_{2}}\right)\cdot \nno\\
&&\hspace{4em}\cdot \lambda(Y_{1}^{o}(v_{1}, y_{1})w_{(1)}
\otimes Y_{2}(v_{2}, y_{2})w_{(2)})\nno\\
&&\quad -\mbox{\rm Res}_{y_{1}}\mbox{\rm Res}_{y_{2}}x_{1}^{-1}\delta
\left(\frac{z-y_{1}}{-x_{1}}\right)x_{2}^{-1}\delta
\left(\frac{y_{2}-z}{x_{2}}\right)\cdot \nno\\
&&\hspace{4em} \cdot \lambda(Y_{1}^{o}(v_{2}, y_{2})w_{(1)}
\otimes Y_{2}(v_{1}, y_{1})w_{(2)})\nno\\
&&\quad +\mbox{\rm Res}_{y_{1}}\mbox{\rm Res}_{y_{2}}x_{1}^{-1}\delta
\left(\frac{z-y_{1}}{-x_{1}}\right)x_{2}^{-1}\delta
\left(\frac{z-y_{2}}{-x_{2}}\right)\cdot \nno\\
&&\hspace{4em}\cdot \lambda(w_{(1)}
\otimes Y_{2}(v_{2}, y_{2})Y_{2}(v_{1}, y_{1})w_{(2)}).
\eea
Transposing the subscripts $1$ and $2$ of the symbols $v$, $x$ and $y$, 
we have
\bea\label{8.2}
\lefteqn{(Y'_{Q(z)}(v_{2}, x_{2})Y'_{Q(z)}(v_{1}, x_{1})\lambda)(w_{(1)}
\otimes w_{(2)})}\nno\\
&&=\mbox{\rm Res}_{y_{2}}\mbox{\rm Res}_{y_{1}}x_{2}^{-1}
\delta\left(\frac{y_{2}-z}{x_{2}}\right)x_{1}^{-1}
\delta\left(\frac{y_{1}-z}{x_{1}}\right)\cdot \nno\\
&&\hspace{4em}\cdot \lambda(Y_{1}^{o}
(v_{1}, y_{1})Y_{1}^{o}(v_{2}, y_{2})w_{(1)}\otimes w_{(2)})\nno\\
&&\quad -\mbox{\rm Res}_{y_{2}}\mbox{\rm Res}_{y_{1}}x_{2}^{-1}
\delta\left(\frac{y_{2}-z}{x_{2}}\right)x_{1}^{-1}
\delta\left(\frac{z-y_{1}}{-x_{1}}\right)\cdot \nno\\
&&\hspace{4em}\cdot \lambda(Y_{1}^{o}(v_{2}, y_{2})w_{(1)}
\otimes Y_{2}(v_{1}, y_{1})w_{(2)})\nno\\
&&\quad -\mbox{\rm Res}_{y_{2}}\mbox{\rm Res}_{y_{1}}x_{2}^{-1}
\delta\left(\frac{z-y_{2}}{-x_{2}}\right)x_{1}^{-1}\delta
\left(\frac{y_{1}-z}{x_{1}}\right)\cdot \nno\\
&&\hspace{4em}\cdot \lambda(Y_{1}^{o}(v_{1}, y_{1})w_{(1)}
\otimes Y_{2}(v_{2}, y_{2})w_{(2)})\nno\\
&&\quad +\mbox{\rm Res}_{y_{2}}\mbox{\rm Res}_{y_{1}}x_{2}^{-1}\delta
\left(\frac{z-y_{2}}{-x_{2}}\right)x_{1}^{-1}\delta\left
(\frac{z-y_{1}}{-x_{1}}\right)\cdot \nno\\
&&\hspace{4em}\cdot \lambda(w_{(1)}\otimes Y_{2}(v_{1}, y_{1})
Y_{2}(v_{2}, y_{2})w_{(2)}).
\eea
Formulas (\ref{8.1}) and (\ref{8.2}) give
\bea\label{8.3}
\lefteqn{([Y'_{Q(z)}(v_{1}, x_{1}), Y'_{Q(z)}(v_{2}, x_{2})]\lambda)
(w_{(1)}\otimes w_{(2)})}\nno\\
&&=\mbox{\rm Res}_{y_{2}}\mbox{\rm Res}_{y_{1}}x_{1}^{-1}
\delta\left(\frac{y_{1}-z}{x_{1}}\right)x_{2}^{-1}\delta
\left(\frac{y_{2}-z}{x_{2}}\right)\cdot \nno\\
&&\hspace{4em}\cdot \lambda([Y_{1}^{o}(v_{2}, y_{2}), Y_{1}^{o}(v_{1}, y_{1})]w_{(1)}
\otimes w_{(2)})\nno\\
&&\quad -\mbox{\rm Res}_{y_{2}}\mbox{\rm Res}_{y_{1}}x_{1}^{-1}
\delta\left(\frac{z-y_{1}}{-x_{1}}\right)x_{2}^{-1}
\delta\left(\frac{z-y_{2}}{-x_{2}}\right)\cdot \nno\\
&&\hspace{4em}\cdot \lambda(w_{(1)}\otimes [Y_{2}(v_{1}, y_{1}), Y_{2}(v_{2}, y_{2})]
w_{(2)})\nno\\
&&=\mbox{\rm Res}_{y_{2}}\mbox{\rm Res}_{y_{1}}x_{1}^{-1}
\delta\left(\frac{y_{1}-z}{x_{1}}\right)x_{2}^{-1}
\delta\left(\frac{y_{2}-z}{x_{2}}\right)\cdot \nno\\
&&\hspace{4em}\cdot \lambda\left(\mbox{\rm Res}_{x_{0}}y_{2}^{-1}
\delta\left(\frac{y_{1}-x_{0}}{y_{2}}\right)
Y^{o}_{1}(Y(v_{1}, x_{0})v_{2}, y_{2})w_{(1)}\otimes w_{(2)}\right) \nno\\
&&\quad -\mbox{\rm Res}_{y_{2}}\mbox{\rm Res}_{y_{1}}x_{1}^{-1}
\delta\left(\frac{z-y_{1}}{-x_{1}}\right)x_{2}^{-1}
\delta\left(\frac{z-y_{2}}{-x_{2}}\right)\cdot \nno\\
&&\hspace{4em}\cdot \lambda\left(
w_{(1)}\otimes \mbox{\rm Res}_{x_{0}}y_{2}^{-1}
\delta\left(\frac{y_{1}-x_{0}}{y_{2}}\right)Y_{2}(Y(v_{1}, x_{0})v_{2}, 
y_{2})w_{(2)}\right)\nno\\
&&=\mbox{\rm Res}_{x_{0}}\mbox{\rm Res}_{y_{2}}
\mbox{\rm Res}_{y_{1}}x_{1}^{-1}\delta\left(\frac{y_{1}-z}{x_{1}}\right)
x_{2}^{-1}\delta\left(\frac{y_{2}-z}{x_{2}}\right)y_{2}^{-1}
\delta\left(\frac{y_{1}-x_{0}}{y_{2}}\right)\cdot \nno\\
&&\hspace{4em}\cdot \lambda(
Y^{o}_{1}(Y(v_{1}, x_{0})v_{2}, y_{2})w_{(1)}\otimes w_{(2)})\nno\\
&&\quad -\mbox{\rm Res}_{x_{0}}\mbox{\rm Res}_{y_{2}}\mbox{\rm Res}_{y_{1}}
x_{1}^{-1}\delta\left(\frac{z-y_{1}}{-x_{1}}\right)x_{2}^{-1}
\delta\left(\frac{z-y_{2}}{-x_{2}}\right)y_{2}^{-1}\delta
\left(\frac{y_{1}-x_{0}}{y_{2}}\right)\cdot \nno\\
&&\hspace{4em}\cdot \lambda(w_{(1)}\otimes Y_{2}(Y(v_{1}, x_{0})v_{2}, y_{2})w_{(2)}).
\eea
But
\bea
\lefteqn{x_{1}^{-1}\delta\left(\frac{y_{1}-z}{x_{1}}\right)x_{2}^{-1}
\delta\left(\frac{y_{2}-z}{x_{2}}\right)y_{2}^{-1}
\delta\left(\frac{y_{1}-x_{0}}{y_{2}}\right)}\nno\\
&&=y_{1}^{-1}\delta\left(\frac{x_{1}+z}{y_{1}}\right)y_{2}^{-1}
\delta\left(\frac{x_{2}+z}{y_{2}}\right)(x_{2}+z)^{-1}
\delta\left(\frac{(x_{1}+z)-x_{0}}{x_{2}+z}\right)\nno\\
&&=y_{1}^{-1}\delta\left(\frac{x_{1}+z}{y_{1}}\right)y_{2}^{-1}
\delta\left(\frac{x_{2}+z}{y_{2}}\right)x_{2}^{-1}
\delta\left(\frac{x_{1}-x_{0}}{x_{2}}\right)\nno\\
&&=y_{1}^{-1}\delta\left(\frac{x_{1}+z}{y_{1}}\right)x_{2}^{-1}
\delta\left(\frac{y_{2}-z}{x_{2}}\right)x_{2}^{-1}
\delta\left(\frac{x_{1}-x_{0}}{x_{2}}\right)
\eea
and
\bea
\lefteqn{x_{1}^{-1}\delta\left(\frac{z-y_{1}}{-x_{1}}\right)x_{2}^{-1}
\delta\left(\frac{z-y_{2}}{-x_{2}}\right)y_{2}^{-1}
\delta\left(\frac{y_{1}-x_{0}}{y_{2}}\right)}\nno\\
&&=z^{-1}\delta\left(\frac{-x_{1}+y_{1}}{z}\right)z^{-1}
\delta\left(\frac{-x_{2}+y_{2}}{z}\right)y_{2}^{-1}
\delta\left(\frac{y_{1}-x_{0}}{y_{2}}\right)\nno\\
&&=\left({\displaystyle \sum_{m, n\in {\Z}}}
\frac{(-x_{1}+y_{1})^{m}}{z^{m+1}}
\frac{(-x_{2}+y_{2})^{n}}{z^{n+1}}\right) 
y_{2}^{-1}\delta\left(\frac{y_{1}-x_{0}}{y_{2}}\right)\nno\\
&&=\left({\displaystyle \sum_{m, n\in {\Z}}}(-x_{2}+y_{2})^{-1}
\left(\frac{-x_{1}+y_{1}}{-x_{2}+y_{2}}\right)^{m}\frac{(-x_{2}+y_{2})^{m+n+1}}
{z^{m+n+2}} \right)
y_{2}^{-1}\delta\left(\frac{y_{1}-x_{0}}{y_{2}}\right)\nno\\
&&=\left({\displaystyle \sum_{m, k\in {\Z}}}(-x_{2}+y_{2})^{-1}
\left(\frac{-x_{1}+y_{1}}{-x_{2}+y_{2}}\right)^{m}
z^{-1}\left(\frac{-x_{2}+y_{2}}{z}\right)^{k}\right) 
y_{2}^{-1}\delta\left(\frac{y_{1}-x_{0}}{y_{2}}\right)\nno\\
&&=(-x_{2}+y_{2})^{-1}\delta\left(\frac{-x_{1}+y_{1}}{-x_{2}+y_{2}}\right)
z^{-1}\delta\left(\frac{-x_{2}+y_{2}}{z}\right)
y_{2}^{-1}\delta\left(\frac{y_{1}-x_{0}}
{y_{2}}\right)\nno\\
&&=(-x_{2})^{-1}\delta\left(\frac{x_{1}-(y_{1}-y_{2})}{x_{2}}\right)
z^{-1}\delta\left(\frac{-x_{2}+y_{2}}{z}\right)
y_{1}^{-1}\delta\left(\frac{y_{2}+x_{0}}
{y_{1}}\right)\nno\\
&&=x_{2}^{-1}\delta\left(\frac{x_{1}-x_{0}}{x_{2}}\right)x_{2}^{-1}
\delta\left(\frac{z-y_{2}}{-x_{2}}\right)y_{1}^{-1}\delta
\left(\frac{y_{2}+x_{0}}{y_{1}}\right).
\eea
Thus (\ref{8.3}) becomes
\bea
\lefteqn{([Y'_{Q(z)}(v_{1}, x_{1}), Y'_{Q(z)}(v_{2}, x_{2})]\lambda)(w_{(1)}
\otimes w_{(2)})}\nno\\
&&=\mbox{\rm Res}_{x_{0}}\mbox{\rm Res}_{y_{2}}
\mbox{\rm Res}_{y_{1}}y_{1}^{-1}\delta\left(\frac{x_{1}+z}{y_{1}}\right)
x_{2}^{-1}\delta\left(\frac{y_{2}-z}{x_{2}}\right)x_{2}^{-1}
\delta\left(\frac{x_{1}-x_{0}}{x_{2}}\right)\cdot \nno\\
&&\hspace{4em}\cdot \lambda(
Y^{o}_{1}(Y(v_{1}, x_{0})v_{2}, y_{2})w_{(1)}\otimes w_{(2)})\nno\\
&&\quad -\mbox{\rm Res}_{x_{0}}
\mbox{\rm Res}_{y_{2}}\mbox{\rm Res}_{y_{1}}x_{2}^{-1}
\delta\left(\frac{x_{1}-x_{0}}{x_{2}}\right)x_{2}^{-1}
\delta\left(\frac{z-y_{2}}{-x_{2}}\right)y_{1}^{-1}
\delta\left(\frac{y_{2}+x_{0}}{y_{1}}\right)\cdot \nno\\
&&\hspace{4em}\cdot \lambda(w_{(1)}\otimes Y_{2}(Y(v_{1}, x_{0})v_{2},
y_{2})w_{(2)})\nno\\
&&=\mbox{\rm Res}_{x_{0}}x_{2}^{-1}\delta\left(\frac{x_{1}-x_{0}}
{x_{2}}\right)\cdot\nno\\
&&\hspace{2em}\cdot \biggl(\mbox{\rm Res}_{y_{2}}x_{2}^{-1}\delta\left(\frac{y_{2}-z}
{x_{2}}\right)\lambda(
Y^{o}_{1}(Y(v_{1}, x_{0})v_{2}, y_{2})w_{(1)}\otimes w_{(2)})\nno\\
&&\hspace{4em}-\mbox{\rm Res}_{y_{2}}x_{2}^{-1}\delta\left(\frac{z-y_{2}}
{-x_{2}}\right)\lambda(w_{(1)}\otimes Y_{2}(Y(v_{1}, x_{0})v_{2}, y_{2})w_{(2)})\biggr)
\nno\\
&&=\mbox{\rm Res}_{x_{0}}x_{2}^{-1}\delta\left(\frac{x_{1}-x_{0}}{x_{2}}\right)
(Y'_{Q(z)}(Y(v_{1}, x_{0})v_{2}, x_{2})\lambda)(w_{(1)}\otimes w_{(2)}).
\eea
Since $\lambda$, $w_{(1)}$ and $w_{(2)}$ are arbitrary, 
this  equality gives the commutator formula (\ref{commu-q-z})
for $Y'_{Q(z)}$. \epfv

When $V$ is in fact a conformal vertex algebra, we write
\begin{equation}\label{13.11-qz}
Y'_{Q(z)}(\omega, x)=\sum_{n\in {\mathbb Z}}L'_{Q(z)}(n)x^{-n-2}.
\end{equation}
Then {}from the last two propositions we see that the coefficient operators
of $Y'_{Q(z)}(\omega, x)$ satisfy the Virasoro algebra commutator
relations:
\begin{equation}\label{5.14}
[L'_{Q(z)}(m), L'_{Q(z)}(n)]
=(m-n)L'_{Q(z)}(m+n)+{\displaystyle\frac1{12}}
(m^3-m)\delta_{m+n,0}c.
\end{equation}
Moreover, in this case, by setting $v=\omega$ in (\ref{Y'qdef}) and
taking $\res_{x_0}x_0^{j+1}$ for $j=-1,0,1$, we see that
\begin{eqnarray}\label{LQ'(j)}
\lefteqn{(L'_{Q(z)}(j)\lambda)(w_{(1)}\otimes w_{(2)})}\nno\\
&&=\res_{x_1}(x_1-z)^{j+1}\lambda(Y^o_1(\omega,x_1)w_{(1)} \otimes
w_{(2)})\nno\\ 
&&\quad- \res_{x_1}
(-z+x_1)^{j+1}\lambda(w_{(1)} \otimes Y_2(\omega,x_1)w_{(2)})\nn
&&=\sum_{i=0}^{j+1}{j+1\choose i}(-z)^{i}\lambda(L(i-j)w_{(1)}
\otimes w_{(2)})\nn
&&\quad -\sum_{i=0}^{j+1}{j+1\choose i}(-z)^{i}
\lambda(w_{(1)}
\otimes L(j-i)w_{(2)})
\end{eqnarray}
for $j=-1,0,1$. If $V$ is just a M\"obius vertex algebra, we
define the actions $L'_{Q(z)}(j)$ on $(W_1\otimes W_2)^*$ by
the right-hand side of 
(\ref{LQ'(j)}) for $j=-1, 0$ and $1$. 

\begin{rema}\label{I-q-intw2}{\rm
In view of the action $L'_{Q(z)}(j)$, the ${\mathfrak s}{\mathfrak
l}(2)$-bracket relations (\ref{imq:Lj}) for a $Q(z)$-intertwining map
can be written as
\begin{equation}\label{I-q-intw2f}
(L'(j)w'_{(3)})\circ I=L'_{Q(z)}(j)(w'_{(3)}\circ I)
\end{equation}
for $w'_{(3)} \in W'_3$ and $j=-1$, $0$, and $1$.  }
\end{rema}

\begin{rema}\label{L'qjpreservesbetaspace}
{\rm  We have 
\[
L'_{Q(z)}(j)((W_{1}\otimes W_{2})^{*})^{(\beta)}\subset 
((W_{1}\otimes W_{2})^{*})^{(\beta)}
\]
for $j=-1, 0, 1$ and $\beta\in \tilde{A}$ (cf. Proposition \ref{tau-q-a-comp}).
}
\end{rema}

In the case that $V$ is a conformal vertex algebra, $L'_{Q(z)}(-1)$,
$L'_{Q(z)}(0)$ and $L'_{Q(z)}(1)$ realize the actions of $L_{-1}$,
$L_0$ and $L_1$ in ${\mathfrak s}{\mathfrak l}(2)$ (cf.\ (\ref{L_*}))
on $(W_1\otimes W_2)^*$.  In the case that $V$ is just a M\"{o}bius
vertex algebra, we now state this fact as a proposition.  This
proposition is needed in the proof of Theorem \ref{q-wk-mod} and
therefore also for Theorems \ref{q-generation} and
\ref{q-characterizationofbackslash}, but neither this proposition nor
any of these three theorems are needed anywhere else in this work, so
we omit the proof of this proposition.  Of course, however, the proof
is straightforward, as is the case with all the ${\mathfrak
s}{\mathfrak l}(2)$ formulas.

\begin{propo}\label{q-sl-2}
Let $V$ be a M\"{o}bius vertex algebra and let $W_{1}$ and $W_{2}$ be
generalized $V$-modules.  Then the operators $L'_{Q(z)}(-1)$,
$L'_{Q(z)}(0)$ and $L'_{Q(z)}(1)$ realize the actions of $L_{-1}$,
$L_0$ and $L_1$ in ${\mathfrak s}{\mathfrak l}(2)$ on $(W_1\otimes
W_2)^*$.
\end{propo}

We also have:

\begin{propo}\label{qz-l-y-comm}
Let $V$ be a M\"obius vertex algebra and let $W_{1}$ and $W_{2}$ be
generalized $V$-modules.  Then for $v\in V$,
\begin{eqnarray}
{[L(-1), Y'_{Q(z)}(v, x)]}&=&Y'_{Q(z)}(L(-1)v, x),\label{qz-sl-2-qz-y-1}\\
{[L(0), Y'_{Q(z)}(v, x)]}&=&Y'_{Q(z)}(L(0)v, x)+xY'_{Q(z)}(L(-1)v, x),
\label{qz-sl-2-qz-y-2}\\
{[L(1), Y'_{Q(z)}(v, x)]}&=&Y'_{Q(z)}(L(1)v, x)
+2xY'_{Q(z)}(L(0)v, x)+x^{2}Y'_{Q(z)}(L(-1)v, x),\label{qz-sl-2-qz-y-3}\nn
&&
\end{eqnarray}
where for brevity we write $L'_{Q(z)}(j)$ acting on 
$(W_{1}\otimes W_{2})^{*}$ as $L(j)$. 
\end{propo}
\pf We prove only (\ref{qz-sl-2-qz-y-2}) since it is needed for Remark
\ref{q-stableundercomponentops} and in Section 6.  We omit the proofs
of (\ref{qz-sl-2-qz-y-1}) and (\ref{qz-sl-2-qz-y-3}) for the same
reasons as above; they are used only for Theorems
\ref{q-wk-mod}--\ref{q-characterizationofbackslash}.

Let $\lambda\in (W_{1}\otimes W_{2})^{*}$, $w_{(1)}\in W_{1}$ and
$w_{(2)}\in W_{2}$. Using (\ref{LQ'(j)}), (\ref{Y'qdef}), the
commutator formulas for $L(j)$ and $Y_{1}(v, x_{0})$ for $j=-1, 0, 1$
and $v\in V$ (recall Definition \ref{moduleMobius}), and the
commutator formulas for $L(j)$ and $Y_{2}^{o}(v, x)$ for $j=-1, 0, 1$
and $v\in V$ (recall Lemma \ref{sl2opposite}), we obtain
\begin{eqnarray*}
\lefteqn{([L(0), Y'_{Q(z)}(v, x)]\lambda)(w_{(1)}\otimes w_{(2)})}\nn
&&=(Y'_{Q(z)}(v, x)\lambda)(L(0)w_{(1)}\otimes w_{(2)})\nn
&&\quad -z(Y'_{Q(z)}(v, x)\lambda)(L(1)w_{(1)}\otimes w_{(2)})\nn
&&\quad -
(Y'_{Q(z)}(v, x)\lambda)(w_{(1)}\otimes L(0)w_{(2)})\nn
&&\quad +z
(Y'_{Q(z)}(v, x)\lambda)(w_{(1)}\otimes L(-1)w_{(2)})\nn
&&\quad -(L(0)\lambda)(Y^o_1(v, x + z)w_{(1)}
\otimes w_{(2)})\nn
&&\quad +\res_{x_1}x^{-1} \delta \left(\frac{z-x_1}{-x}\right)
(L(0)\lambda)(w_{(1)}
\otimes Y_2(v,x_1)w_{(2)})\nn
&&=\lambda(Y^o_1(v, x + z)L(0)w_{(1)}\otimes w_{(2)})\nn
&&\quad-\res_{x_1}x^{-1} \delta \left(\frac{z-x_1}{-x}\right)
\lambda(L(0)w_{(1)}\otimes Y_2(v,x_1)w_{(2)})\nn
&&\quad -z\lambda(Y^o_1(v, x + z)L(1)w_{(1)}\otimes w_{(2)})\nn
&&\quad +z\res_{x_1}x^{-1} \delta \left(\frac{z-x_1}{-x}\right)
\lambda(L(1)w_{(1)}\otimes Y_2(v,x_1)w_{(2)})\nn
&&\quad -
\lambda(Y^o_1(v, x + z)w_{(1)}\otimes L(0)w_{(2)})\nn
&&\quad +\res_{x_1}x^{-1} \delta \left(\frac{z-x_1}{-x}\right)
\lambda(w_{(1)}\otimes Y_2(v,x_1)L(0)w_{(2)})\nn
&&\quad +z
\lambda(Y^o_1(v, x + z)w_{(1)}\otimes L(-1)w_{(2)})\nn
&&\quad -z\res_{x_1}x^{-1} \delta \left(\frac{z-x_1}{-x}\right)
\lambda(w_{(1)}\otimes Y_2(v,x_1)L(-1)w_{(2)})\nn
&&\quad -\lambda(L(0)Y^o_1(v, x + z)w_{(1)}
\otimes w_{(2)})\nn
&&\quad +z\lambda(L(1)Y^o_1(v, x + z)w_{(1)}
\otimes w_{(2)})\nn
&&\quad +\lambda(Y^o_1(v, x + z)w_{(1)}
\otimes L(0)w_{(2)})\nn
&&\quad -z\lambda(Y^o_1(v, x + z)w_{(1)}
\otimes L(-1)w_{(2)})\nn
&&\quad +\res_{x_1}x^{-1} \delta \left(\frac{z-x_1}{-x}\right)
\lambda(L(0)w_{(1)}
\otimes Y_2(v,x_1)w_{(2)})\nn
&&\quad -z\res_{x_1}x^{-1} \delta \left(\frac{z-x_1}{-x}\right)
\lambda(L(1)w_{(1)}
\otimes Y_2(v,x_1)w_{(2)})\nn
&&\quad -\res_{x_1}x^{-1} \delta \left(\frac{z-x_1}{-x}\right)
\lambda(w_{(1)}
\otimes L(0)Y_2(v,x_1)w_{(2)})\nn
&&\quad +z\res_{x_1}x^{-1} \delta \left(\frac{z-x_1}{-x}\right)
\lambda(w_{(1)}
\otimes L(-1)Y_2(v,x_1)w_{(2)})\nn
&&=\lambda([Y^o_1(v, x + z), L(0)]w_{(1)}\otimes w_{(2)})\nn
&&\quad -z\lambda([Y^o_1(v, x + z), L(1)]w_{(1)}\otimes w_{(2)})\nn
&&\quad +z\res_{x_1}x^{-1} \delta \left(\frac{z-x_1}{-x}\right)
\lambda(w_{(1)}\otimes [L(-1), Y_2(v,x_1)]w_{(2)})\nn
&&\quad -\res_{x_1}x^{-1} \delta \left(\frac{z-x_1}{-x}\right)
\lambda(w_{(1)}
\otimes [L(0), Y_2(v,x_1)]w_{(2)})\nn
&&=\lambda(Y^o_1((L(0)+(x+z)L(-1)v, x + z)w_{(1)}\otimes w_{(2)})\nn
&&\quad -z\lambda(Y^o_1(L(-1)v, x + z)w_{(1)}\otimes w_{(2)})\nn
&&\quad +z\res_{x_1}x^{-1} \delta \left(\frac{z-x_1}{-x}\right)
\lambda(w_{(1)}\otimes Y_2(L(-1)v,x_1)w_{(2)})\nn
&&\quad -\res_{x_1}x^{-1} \delta \left(\frac{z-x_1}{-x}\right)
\lambda(w_{(1)}
\otimes Y_2((L(0)+x_{1}L(-1))v,x_1)w_{(2)})\nn
&&=\lambda(Y^o_1((L(0)+xL(-1)v, x + z)w_{(1)}\otimes w_{(2)})\nn
&&\quad -\res_{x_1}x^{-1} \delta \left(\frac{z-x_1}{-x}\right)
\lambda(w_{(1)}
\otimes Y_2((L(0)+(x-z)L(-1))v,x_1)w_{(2)})\nn
&&=\lambda(Y^o_1((L(0)+xL(-1)v, x + z)w_{(1)}\otimes w_{(2)})\nn
&&\quad -\res_{x_1}x^{-1} \delta \left(\frac{z-x_1}{-x}\right)
\lambda(w_{(1)}
\otimes Y_2((L(0)+xL(-1))v,x_1)w_{(2)})\nn
&&=(Y'_{Q(z)}((L(0)+xL(-1))v, x)\lambda)(w_{(1)}\otimes w_{(2)})\nn
&&=(Y'_{Q(z)}(L(0)v, x)\lambda)(w_{(1)}\otimes w_{(2)})
+(xY'_{Q(z)}(L(-1)v, x)\lambda)(w_{(1)}\otimes w_{(2)}),
\end{eqnarray*}
proving (\ref{qz-sl-2-qz-y-2}).
\epfv

Let $W_3$ also be an object of ${\cal C}$. Note that $V\otimes
\iota_{+}{\mathbb C}[t, t^{-1}, (z+t)^{-1}]$ acts on $W'_3$ in the
natural way.  The following result provides further motivation for the
definition of our action (\ref{5.2}) on $(W_1\otimes W_2)^{*}$; recall
the discussion preceding Proposition \ref{pz}:

\begin{propo}\label{qz}
Let $W_{1}$, $W_{2}$ and $W_{3}$ be generalized $V$-modules.
Under the natural isomorphism described in Remark \ref{alternateformoflemma}
between the space of $\tilde{A}$-compatible linear maps
\[
I:W_{1}\otimes W_{2} \rightarrow \overline{W_{3}}
\]
and the space of $\tilde{A}$-compatible linear maps
\[
J:W'_{3} \rightarrow (W_{1}\otimes W_{2})^{*}
\]
determined by (\ref{IcorrespondstoJ}), the $Q(z)$-intertwining maps
$I$ of type ${W_3\choose W_1\, W_2}$ correspond exactly to the
(grading restricted) $\tilde{A}$-compatible maps $J$ that intertwine
the actions of both
\[
V \otimes \iota_{+}{\mathbb C}[t, t^{-1}, (z+t)^{-1}]
\]
and ${\mathfrak s} {\mathfrak l}(2)$ on $W'_{3}$ and on $(W_1\otimes
W_2)^{*}$.  If $W_3$ is lower bounded, we may replace the grading
restrictions by (\ref{QpinI=0}) and (\ref{Jlowerbdd}).
\end{propo}
\pf 
In view of (\ref{IcorrespondstoJalternateform}), Remark
\ref{I-intw-q} asserts that (\ref{imq:def'}), or equivalently,
(\ref{imq:def}), is equivalent to the condition
\begin{equation}\label{q-j-tau}
J\left(\tau_{W'_3}\left(z^{-1}\delta\left(\frac{x_1-x_0}{z}\right)
Y_{t}(v, x_0)\right)w'_{(3)}\right)
=\tau_{Q(z)}\left(z^{-1}\delta\left(\frac{x_1-x_0}{z}\right)
Y_{t}(v, x_0)\right)J(w'_{(3)}),
\end{equation}
that is, the condition that $J$ intertwines the actions of $V \otimes
\iota_{+}{\mathbb C}[t,t^{- 1}, (z+t)^{-1}]$ on $W'_{3}$ and on
$(W_1\otimes W_2)^{*}$ (recall (\ref{3.12})--(\ref{3.13})).
Similarly, Remark \ref{I-q-intw2} asserts that (\ref{imq:Lj}) is
equivalent to the condition
\begin{equation}\label{q-j-lj}
J(L'(j)w'_{(3)}) = L'_{Q(z)}(j)J(w'_{(3)})
\end{equation}
for $j=-1$, $0$, $1$, that is, the condition that $J$ intertwines the
actions of ${\mathfrak s} {\mathfrak l}(2)$ on $W'_{3}$ and on
$(W_1\otimes W_2)^{*}$.
\epfv

\begin{nota}\label{qscriptN}
{\rm Given generalized $V$-modules $W_1$, $W_2$ and $W_3$, we shall
write ${\cal N}[Q(z)]_{W'_3}^{(W_1 \otimes W_2)^{*}}$ for the
space of (grading restricted) $\tilde{A}$-compatible linear maps
\[
J:W'_{3} \rightarrow (W_{1}\otimes W_{2})^{*}
\]
that intertwine the actions of both
\[
V \otimes \iota_{+}{\mathbb C}[t,t^{- 1}, (z+t)^{-1}]
\]
and ${\mathfrak s} {\mathfrak l}(2)$ on $W'_{3}$ and on $(W_1\otimes
W_2)^{*}$.
Note that Proposition \ref{qz} gives a natural linear isomorphism
\begin{eqnarray*}
{\cal M}[Q(z)]^{W_3}_{W_1 W_2} &
\stackrel{\sim}{\longrightarrow} & {\cal N}[Q(z)]_{W'_3}^{(W_1 \otimes
W_2)^{*}}\nno\\ I & \mapsto & J
\end{eqnarray*}
(recall {}from Definition \ref{im:qimdef} the notation for the space
of $Q(z)$-intertwining maps), and if $W_3$ is lower bounded, the
spaces satisfy the stronger grading restrictions (\ref{QpinI=0}) and
(\ref{Jlowerbdd}).  As in Notation \ref{scriptN}, we still use the
symbol ``prime'' to denote this isomorphism in both directions:
\begin{eqnarray*}
{\cal M}[Q(z)]^{W_3}_{W_1 W_2} & \stackrel{\sim}{\longrightarrow} & {\cal
N}[Q(z)]_{W'_3}^{(W_1 \otimes W_2)^{*}}\nno\\
I & \mapsto & I'\nno\\
J' & \leftarrow\!\!\!{\scriptstyle |} & J,
\end{eqnarray*}
so that in particular,
\[
I'' = I \;\;\mbox{ and }\;\; J'' = J
\]
for $I \in {\cal M}[Q(z)]^{W_3}_{W_1 W_2}$ and 
$J \in {\cal N}[Q(z)]_{W'_3}^{(W_1
\otimes W_2)^{*}}$, and the relation between $I$ and $I'$ is
determined by
\[
\langle w'_{(3)}, I(w_{(1)}\otimes w_{(2)})\rangle
=I'(w'_{(3)})(w_{(1)}\otimes w_{(2)})
\]
for $w_{(1)}\in W_{1}$, $w_{(2)}\in W_{2}$ and $w'_{(3)}\in W'_{3}$,
or equivalently,
\[
w'_{(3)}\circ I = I'(w'_{(3)}).
\]
}
\end{nota}

\begin{rema}
{\rm Combining Proposition \ref{qz} with Proposition \ref{Q-cor}, we see
that for any integer $p$, we also have a natural linear
isomorphism
\[
{\cal N}[Q(z)]_{W'_3}^{(W_1 \otimes W_2)^{*}} \stackrel{\sim}{\longrightarrow}
{\cal V}^{W'_1}_{W'_3 W_2}
\]
{}from ${\cal N}[Q(z)]_{W'_3}^{(W_1 \otimes W_2)^{*}}$ to the space of
logarithmic intertwining operators of type ${W'_1\choose W'_3\,W_2}$.
In particular, given any such logarithmic intertwining
operator ${\cal Y}$ and integer $p$, the map 
\[
(I^{Q(z)}_{{\cal Y}, p})':
W'_3\to (W_1\otimes W_2)^{*}
\]
defined by
\[
(I^{Q(z)}_{{\cal Y}, p})'(w'_{(3)})(w_{(1)}\otimes w_{(2)})=\bra w_{(1)},
{\cal Y}(w'_{(3)}, e^{l_{p}(z)})w_{(2)}\ket_{W'_1}
\]
is $\tilde{A}$-compatible and intertwines both actions on both spaces.
If the modules involved are lower bounded, we have the stronger
grading restrictions (cf. Remark \ref{NisotoV}).}
\end{rema}

We have formulated the notions of $Q(z)$-product and $Q(z)$-tensor
product using $Q(z)$-intertwining maps (Definitions \ref{qz-product}
and \ref{qz-tp}).  Now that we know that $Q(z)$-intertwining maps can
be interpreted as in Proposition \ref{qz} (and Notation
\ref{qscriptN}), we can reformulate the notions of $Q(z)$-product and
$Q(z)$-tensor product correspondingly (the proof of the next result is
the same as that of Proposition \ref{productusingI'}):

\begin{propo}\label{q-productusingI'}
Let ${\cal C}_1$ be either of the categories ${\cal M}_{sg}$ or ${\cal
GM}_{sg}$, as in Definition \ref{qz-product}.  For $W_1, W_2\in
\ob{\cal C}_1$, a $Q(z)$-product $(W_3;I_3)$ of $W_1$ and $W_2$
(recall Definition \ref{qz-product}) amounts to an object $(W_3,Y_3)$
of ${\cal C}_1$ equipped with a map $I'_3 \in {\cal N}[Q(z)]_{W'_3}^{(W_1
\otimes W_2)^{*}}$, that is, equipped with an $\tilde{A}$-compatible
map
\[
I'_3:W'_{3} \rightarrow (W_{1}\otimes W_{2})^{*}
\]
that intertwines the two actions of $V \otimes \iota_{+}{\mathbb
C}[t,t^{- 1}, (z+t)^{-1}]$ and of ${\mathfrak s} {\mathfrak
l}(2)$.  The map $I'_3$ corresponds to the $Q(z)$-intertwining map
\[
I_3:W_{1}\otimes W_{2} \rightarrow \overline{W_{3}}
\]
as above:
\[
I'_3(w'_{(3)}) = w'_{(3)}\circ I_3
\]
for $w'_{(3)} \in W'_{3}$ (recall
\ref{IcorrespondstoJalternateform})).  Denoting this structure by
$(W_3,Y_3;I'_3)$ or simply by $(W_3;I'_3)$, let $(W_4;I'_4)$ be
another such structure.  Then a morphism of $Q(z)$-products {}from $W_3$
to $W_4$ amounts to a module map $\eta: W_3 \to W_4$ such that the
diagram
\begin{center}
\begin{picture}(100,60)
\put(-2,0){$W'_4$}
\put(13,4){\vector(1,0){104}}
\put(119,0){$W'_3$}
\put(38,50){$(W_1\otimes W_2)^*$}
\put(13,12){\vector(3,2){50}}
\put(118,12){\vector(-3,2){50}}
\put(65,8){$\eta'$}
\put(23,27){$I'_4$}
\put(98,27){$I'_3$}
\end{picture}
\end{center}
commutes, where $\eta'$ is the natural map given by (\ref{fprime}).
\epf
\end{propo}

\begin{corol}\label{q-tensorproductusingI'}
Let ${\cal C}$ be a full subcategory of either ${\cal M}_{sg}$ or
${\cal GM}_{sg}$, as in Definition \ref{qz-tp}.  For $W_1, W_2\in
\ob{\cal C}$, a $Q(z)$-tensor product $(W_0; I_0)$ of $W_1$ and $W_2$
in ${\cal C}$, if it exists, amounts to an object $W_0 =
W_1\boxtimes_{Q(z)} W_2$ of ${\cal C}$ and a structure $(W_0 =
W_1\boxtimes_{Q(z)} W_2; I'_0)$ as in Proposition
\ref{q-productusingI'}, with
\[
I'_0: (W_1\boxtimes_{Q(z)} W_2)' \longrightarrow (W_1\otimes W_2)^*
\]
in ${\cal N}[Q(z)]_{(W_1\boxtimes_{Q(z)} W_2)'}^{(W_1 \otimes W_2)^{*}}$,
such that for any such pair $(W; I')$ $(W\in \ob \mathcal{C})$, with
\[
I': W' \longrightarrow (W_1\otimes W_2)^*
\]
in ${\cal N}[Q(z)]_{W'}^{(W_1 \otimes W_2)^{*}}$, there is a unique module
map
\[
\chi: W' \longrightarrow (W_1\boxtimes_{Q(z)} W_2)'
\]
such that the diagram
\begin{center}
\begin{picture}(100,60)
\put(-2,0){$W'$}
\put(13,4){\vector(1,0){104}}
\put(119,0){$(W_1\boxtimes_{Q(z)} W_2)'$}
\put(38,50){$(W_1\otimes W_2)^*$}
\put(13,12){\vector(3,2){50}}
\put(118,12){\vector(-3,2){50}}
\put(65,8){$\chi$}
\put(23,27){$I'$}
\put(98,27){$I'_0$}
\end{picture}
\end{center}
commutes.  Here $\chi = \eta'$, where $\eta$ is a correspondingly
unique module map
\[
\eta: W_1\boxtimes_{Q(z)} W_2 \longrightarrow W.
\]
Also, the map $I_0'$, which is $\tilde{A}$-compatible and which
intertwines the two actions of $V \otimes \iota_{+}{\mathbb C}[t,t^{-
1}, (z+t)^{-1}]$ and of ${\mathfrak s} {\mathfrak l}(2)$, is
related to the $Q(z)$-intertwining map
\[
I_0 = \boxtimes_{Q(z)}: W_1\otimes W_2 \longrightarrow 
\overline{W_1\boxtimes_{Q(z)} W_2}
\]
by
\[
I_0'(w') = w' \circ \boxtimes_{Q(z)}
\]
for $w' \in (W_1\boxtimes_{Q(z)} W_2)'$, that is,
\[
I_0'(w')(w_{(1)}\otimes w_{(2)}) = \langle w',w_{(1)}\boxtimes_{Q(z)}
w_{(2)} \rangle
\]
for $w_{(1)}\in W_{1}$ and $w_{(2)}\in W_{2}$, using the notation 
(\ref{q-boxtensorofelements}).
\epf
\end{corol}

\begin{defi}
{\rm For $W_{1}, W_{2}\in \ob \mathcal{C}$, define the subset 
\[
W_{1}\hboxtr_{Q(z)}W_{2}\subset (W_{1}\otimes W_{2})^{*}
\]
of $(W_{1}\otimes W_{2})^{*}$ to be the union of the images
\[
I'(W')\subset (W_{1}\otimes W_{2})^{*}
\]
as $(W; I)$ ranges through all the $Q(z)$-products of $W_{1}$ and $W_{2}$ with 
$W\in \ob
\mathcal{C}$. Equivalently, $W_{1}\hboxtr_{P(z)}W_{2}$ is the union 
of the images 
$I'(W')$ as $W$ (or $W'$) ranges through $\ob
\mathcal{C}$ and $I'$ ranges through 
$\mathcal{N}[Q(z)]_{W'}^{(W_{1}\otimes W_{2})^{*}}$---the space of 
$\tilde{A}$-compatible linear maps
\[
W'\to (W_{1}\otimes W_{2})^{*}
\]
intertwining the actions of both 
\[
V\otimes \iota_{+}\C[t, t^{-1}, (z+t)^{-1}]
\]
and $\mathfrak{s}\mathfrak{l}(2)$ on both spaces.}
\end{defi}

\begin{rema}
{\rm Since $\mathcal{C}$ is closed under finite direct sums (Assumption 
\ref{assum-c}), it is clear that $W_{1}\hboxtr_{Q(z)}W_{2}$ is
in fact a linear subspace of $(W_{1}\otimes W_{2})^{*}$, and in particular,
it can be defined alternatively as the sum of all the images $I'(W')$:
\begin{equation}\label{q-hboxtr-sum}
W_{1}\hboxtr_{Q(z)}W_{2}=\sum I'(W') = \bigcup I'(W')\subset 
(W_{1}\otimes W_{2})^{*},
\end{equation}
where the sum and union both range over $W\in \ob
\mathcal{C}$, $I\in \mathcal{M}[Q(z)]_{W_{1}W_{2}}^{W}$.}
\end{rema}

For any generalized $V$-modules $W_{1}$ and $W_{2}$,
using the operator $L'_{Q(z)}(0)$ (recall (\ref{LQ'(j)}))
on $(W_{1}\otimes W_{2})^{*}$ we define
the generalized $L'_{Q(z)}(0)$-eigenspaces 
$((W_{1}\otimes W_{2})^{*})_{[n]; Q(z)}$ for $n\in \C$ in the usual way:
\begin{equation}
((W_{1}\otimes W_{2})^{*})_{[n]; Q(z)}=\{w\in (W_{1}\otimes W_{2})^{*}\;|\;
(L'_{Q(z)}(0)-n)^{m}w=0 \;{\rm for}\; m\in \N \;
\mbox{\rm sufficiently large}\}.
\end{equation}
Then we have the (proper) subspace 
\begin{equation}
\coprod_{n\in \C}((W_{1}\otimes W_{2})^{*})_{[n]; Q(z)}\subset 
(W_{1}\otimes W_{2})^{*}.
\end{equation}
We also define the ordinary $L'_{Q(z)}(0)$-eigenspaces
$((W_{1}\otimes W_{2})^{*})_{(n); Q(z)}$ in the usual way:
\begin{equation}
((W_{1}\otimes W_{2})^{*})_{(n); Q(z)}=\{w\in (W_{1}\otimes W_{2})^{*}\;|\;
L'_{P(z)}(0)w=nw \}.
\end{equation}
Then we have the (proper) subspace
\begin{equation}
\coprod_{n\in \C}((W_{1}\otimes W_{2})^{*})_{(n); Q(z)}\subset 
(W_{1}\otimes W_{2})^{*}.
\end{equation}

Just as in Proposition \ref{im:abc}, we have:

\begin{propo}\label{im-q:abc}
Let $W_{1}, W_{2}\in \ob \mathcal{C}$. 

(a)
The elements of $W_1\hboxtr_{Q(z)} W_2$ are exactly the linear functionals
on $W_{1}\otimes W_{2}$ of the form $w'\circ
I(\cdot\otimes \cdot)$ for some $Q(z)$-intertwining map $I$ of type
${W\choose W_1\,W_2}$ and some $w'\in W'$, $W\in\ob{\cal C}$.

(b) Let $(W; I)$ be any $Q(z)$-product of $W_{1}$ and $W_{2}$, with 
$W$ any generalized $V$-module. Then for $n\in \C$,
\[
I'(W'_{[n]}) \subset  ((W_{1}\otimes W_{2})^{*})_{[n]; Q(z)}
\]
and 
\[
I'(W'_{(n)})\subset ((W_{1}\otimes W_{2})^{*})_{(n); Q(z)}.
\]

(c) The structure $(W_1\hboxtr_{Q(z)} W_2,Y'_{Q(z)})$ (recall
(\ref{y'-q-z})) satisfies all the axioms in the definition of
(strongly $\tilde{A}$-graded) generalized $V$-module except perhaps
for the two grading conditions (\ref{set:dmltc}) and
(\ref{set:dmfin}).

(d) Suppose that the objects of the category $\mathcal{C}$ consist
only of (strongly $\tilde{A}$-graded) {\em ordinary}, as opposed to
{\em generalized}, $V$-modules.  Then the structure
$(W_1\hboxtr_{Q(z)} W_2,Y'_{Q(z)})$ satisfies all the axioms in the
definition of (strongly $\tilde{A}$-graded ordinary) $V$-module except
perhaps for (\ref{set:dmltc}) and (\ref{set:dmfin}).
\end{propo}
\pf 
Part (a) is clear {}from the definition of $W_1\hboxtr_{Q(z)} W_2$, and 
(b) follows {}from (\ref{q-j-lj}) with $j=0$.

To prove (c), let $(W; I)$ be any any $Q(z)$-product of $W_{1}$ and
$W_{2}$, with $W$ any generalized $V$-module. Then $(I'(W'),
Y'_{Q(z)})$ satisfies all the conditions in the definition of
(strongly $\tilde{A}$-graded) generalized $V$-module since $I'$ is
$\tilde{A}$-compatible and intertwines the actions of $V\otimes
{\mathbb C}[t,t^{-1}]$ and of ${\mathfrak s}{\mathfrak l}(2)$; the
$\C$-grading follows {}from Part (b).  Note that
\begin{equation}
I':W' \rightarrow I'(W')
\end{equation}
is a map of generalized $V$-modules.  Since $W_1\hboxtr_{Q(z)} W_2$ is
the sum of these structures $I'(W')$ over $W\in \ob \mathcal{C}$
(recall (\ref{hboxtr-sum})), $(W_1\hboxtr_{Q(z)} W_2, Y'_{Q(z)})$
satisfies all the conditions in the definition of generalized module
except perhaps for (\ref{set:dmltc}) and (\ref{set:dmfin}).

Part (d) is proved by the same argument as for (c): For $(W; I)$ any
$Q(z)$-product of possibly generalized $V$-modules $W_{1}$ and
$W_{2}$, with $W$ any ordinary $V$-module, $(I'(W'), Y'_{Q(z)})$
satisfies all the conditions in the definition of (strongly
$\tilde{A}$-graded) ordinary $V$-module; the $\C$-grading (by ordinary
$L'_{Q(z)}(0)$-eigenspaces) again follows {}from Part (b).  \epfv

Just as in Proposition \ref{backslash=union}, we have:

\begin{propo}\label{Qbackslash=union}
Suppose that $\mathcal{C}$ is closed under images (as well as under
contragredients and finite direct sums (Assumption \ref{assum-c})).
Let $W_1,W_2 \in \ob{\cal C}$.  Then the subspace $W_1\hboxtr_{Q(z)}
W_2$ of $(W_{1}\otimes W_{2})^{*}$ is equal to the union and also to
the sum of the objects of $\mathcal{C}$ lying in $(W_{1}\otimes
W_{2})^{*}$:
\[
W_{1}\hboxtr_{Q(z)}W_{2}=\bigcup W = \sum W \subset 
(W_{1}\otimes W_{2})^{*},
\]
where in the union and in the sum, $W$ ranges through the subspaces of
$(W_{1}\otimes W_{2})^{*}$ that are objects of $\mathcal{C}$ when
equipped with the action $Y'_{Q(z)}(\cdot,x)$ of $V$ and the
corresponding action of ${\mathfrak s}{\mathfrak l}(2)$ on
$(W_{1}\otimes W_{2})^{*}$.  In particular, every object of
$\mathcal{C}$ lying in $(W_{1}\otimes W_{2})^{*}$ is a subspace of
$W_{1}\hboxtr_{Q(z)}W_{2}$ (and for this assertion, the assumption
that $\mathcal{C}$ is closed under images is not needed).\epf
\end{propo}

We also have the following generalization of Proposition 5.8 in
\cite{tensor1}, characterizing $W_{1}\boxtimes_{Q(z)}W_{2}$, including its
existence, in terms of $W_1\hboxtr_{Q(z)} W_2$; 
the proof is the same as that of Proposition 
\ref{tensor1-13.7}:

\begin{propo}\label{tensor1-5.7}
Let $W_{1}, W_{2}\in \ob \mathcal{C}$.  If $(W_1\hboxtr_{Q(z)} W_2,
Y'_{Q(z)})$ is an object of ${\cal C}$, denote by
$(W_1\boxtimes_{Q(z)} W_2, Y_{Q(z)})$ its contragredient module:
\[
W_1\boxtimes_{Q(z)} W_2 = (W_1\hboxtr_{Q(z)} W_2)'.
\]
Then the $Q(z)$-tensor product of $W_{1}$ and $W_{2}$ in ${\cal C}$ exists
and is
\[
(W_1\boxtimes_{Q(z)} W_2, Y_{Q(z)}; i'),
\]
where $i$ is the
natural inclusion {}from $W_1\hboxtr_{Q(z)} W_2$ to $(W_1\otimes
W_2)^*$ (recall Notation \ref{qscriptN}).  Conversely, let us assume
that $\mathcal{C}$ is closed under images (recall Definition
\ref{closedunderimages}).  If the $Q(z)$-tensor product of $W_1$ and
$W_2$ in ${\cal C}$ exists, then $(W_1\hboxtr_{Q(z)} W_2, Y'_{Q(z)})$
is an object of ${\cal C}$.\epf
\end{propo}

\begin{rema}
{\rm Suppose that $W_1\hboxtr_{Q(z)} W_2$ is an object of
$\mathcal{C}$.  {}From Corollary \ref{q-tensorproductusingI'} and
Proposition \ref{tensor1-5.7} we see that
\begin{equation}\label{boxpair-q}
\langle\lambda, w_{(1)}\boxtimes_{Q(z)}w_{(2)}\rangle
_{W_1\boxtimes_{Q(z)} W_2}=
\lambda(w_{(1)}\otimes w_{(2)})
\end{equation}
for $\lambda\in W_1\hboxtr_{Q(z)} W_2\subset (W_1\otimes W_2)^*$,
$w_{(1)}\in W_1$ and $w_{(2)}\in W_2$.}
\end{rema}

As in the $P(z)$-case, our next goal is 
to present an alternative description of the
subspace $W_1\hboxtr_{Q(z)} W_2$ of $(W_1\otimes W_2)^*$. The main
ingredient of this description will be the ``$Q(z)$-compatibility
condition,'' as was the case in \cite{tensor1}--\cite{tensor2}.

Take $W_{1}$ and $W_{2}$ to be arbitrary generalized $V$-modules.  Let
$(W,I)$ ($W$ a generalized $V$-module) be a $Q(z)$-product of $W_{1}$
and $W_{2}$ and let $w'\in W'$. Then {}from (\ref{q-j-tau}),
Proposition \ref{q-productusingI'}, (\ref{tau-w-comp}), (\ref{3.7})
and (\ref{y'-q-z}), we have, for all $v\in V$,
\begin{eqnarray}\label{5.18}
\lefteqn{\tau_{Q(z)}\left(z^{-1}\delta\left(\frac{x_1-x_0}{z}\right)
Y_{t}(v, x_0)\right)I'(w')}\nno\\
&&=I'\left(\tau_{W'}\left(z^{-1}\delta\left(\frac{x_1-x_0}{z}\right)
Y_{t}(v, x_0)\right)w'\right)\nno\\
&&=I'\left(z^{-1}\delta\left(\frac{x_1-x_0}{z}\right)
Y_{W'}(v, x_0)w'\right)\nno\\
&&=z^{-1}\delta\left(\frac{x_1-x_0}{z}\right)I'(Y_{W'}(v,
x_0)w')\nno\\
&&=z^{-1}\delta\left(\frac{x_1-x_0}{z}\right)I'(\tau_{W'}(Y_{t}(v,
x_0))w')\nno\\
&&=z^{-1}\delta\left(\frac{x_1-x_0}{z}\right)\tau_{Q(z)}(Y_{t}(v,
x_0))I'(w')\nn
&&=z^{-1}\delta\left(\frac{x_1-x_0}{z}\right)
Y'_{Q(z)}(v, x_0)I'(w').
\end{eqnarray}
That is, 
$I'(w')$ satisfies the following nontrivial and subtle condition on
\[
\lambda \in (W_1\otimes W_2)^{*}:
\]

\begin{description}
\item{\bf The $Q(z)$-compatibility condition}

(a) The {\em $Q(z)$-lower truncation condition}: For all $v\in V$, the formal
Laurent series $Y'_{Q(z)}(v, x)\lambda$ involves only finitely many
negative powers of $x$.

(b) The following formula holds:
\begin{eqnarray}\label{cpb-q}
\lefteqn{\tau_{Q(z)}\left(z^{-1}\delta\left(\frac{x_1-x_0}{z}\right)
Y_{t}(v, x_0)\right)\lambda}\nno\\
&&=z^{-1}\delta\left(\frac{x_1-x_0}{z}\right)
Y'_{Q(z)}(v, x_0)\lambda  \;\;\mbox{ for all }\;v\in V.
\end{eqnarray}
(Note that the two sides of (\ref{cpb-q}) are not {\it a priori} equal
for general $\lambda\in (W_1\otimes W_2)^{*}$. Note also that Condition 
(a) insures that the right-hand side in Condition (b) is 
well defined.)
\end{description}

\begin{nota}
{\rm Note that the set of elements of $(W_1\otimes W_2)^*$ satisfying
either  the full $Q(z)$-compatibility
condition or Part (a) of this condition forms a subspace. 
We shall denote the space of elements of $(W_1\otimes W_2)^*$ satisfying
the $Q(z)$-compatibility
condition by
\[
\comp_{Q(z)}((W_1\otimes W_2)^*).
\]}
\end{nota}

We know that each space $((W_1\otimes W_2)^*)^{(\beta)}$ is 
$L'_{Q(z)}(0)$-stable
(recall Proposition \ref{tau-q-a-comp} and Remark
\ref{L'qjpreservesbetaspace}), so that we may consider the subspaces
\[
\coprod_{n\in \C}((W_1\otimes W_2)^*)_{[n];Q(z)}^{(\beta)} \subset
((W_1\otimes W_2)^*)^{(\beta)}
\]
and 
\[
\coprod_{n\in \C}((W_1\otimes W_2)^*)_{(n);Q(z)}^{(\beta)} \subset
((W_1\otimes W_2)^*)^{(\beta)}
\]
(recall Remark \ref{generalizedeigenspacedecomp}).  We define the
two subspaces
\begin{equation}\label{W1W2_[C];q^Atilde}
((W_1\otimes W_2)^*)_{[{\mathbb C}]; Q(z)}^{( \tilde A )}=
\coprod_{n\in
\C}\coprod_{\beta\in \tilde{A}}((W_1\otimes W_2)^*)_{[n];Q(z)}^{(\beta)}
\subset (W_1\otimes W_2)^*
\end{equation}
and 
\begin{equation}\label{W1W2_(C);q^Atilde}
((W_1\otimes W_2)^*)_{({\mathbb C});Q(z)}^{( \tilde A )}=
\coprod_{n\in
\C}\coprod_{\beta\in \tilde{A}}((W_1\otimes W_2)^*)_{(n);Q(z)}^{(\beta)}
\subset (W_1\otimes W_2)^*.
\end{equation}

\begin{rema}\label{q-singleanddoublegraded}
{\rm Any $L'_{Q(z)}(0)$-stable subspace of $((W_1\otimes
W_2)^*)_{[{\mathbb C}];Q(z)}^{( \tilde A )}$ is graded by generalized
eigenspaces (again recall Remark \ref{generalizedeigenspacedecomp}),
and if such a subspace is also $\tilde A$-graded, then it is doubly
graded; similarly for subspaces of $((W_1\otimes W_2)^*)_{({\mathbb
C});Q(z)}^{( \tilde A )}$.}
\end{rema}

We have:

\begin{lemma}\label{q-a-tilde-comp}
Suppose that $\lambda\in ((W_1\otimes W_2)^*)_{[{\mathbb C}];Q(z)}^{(\tilde
A )}$ satisfies the $Q(z)$-compatibility condition. Then every
$\tilde{A}$-homogeneous component of $\lambda$ also satisfies this
condition.
\end{lemma}
\pf
When $v\in V$ is $\tilde{A}$-homogeneous, 
\[
\tau_{Q(z)}\bigg(z^{-1}\delta\bigg(\frac{x_1-x_{0}}{z} \bigg)
Y_{t}(v, x_0)\bigg)\;\;\mbox{ and }\;\;
z^{-1}\delta\bigg(\frac{x_1-x_{0}}{z} \bigg)Y'_{Q(z)}(v, x_0)
\]
are both $\tilde{A}$-homogeneous as operators.  By comparing the
$\tilde{A}$-homogeneous components of both sides of (\ref{cpb-q}), we
see that the $\tilde{A}$-homogeneous components of $\lambda$ also
satisfy the $Q(z)$-compatibility condition.  \epfv

\begin{rema}\label{q-stableundercomponentops}
{\rm Just as in Remark \ref{stableundercomponentops}, note that both
the spaces $((W_1\otimes W_2)^*)_{[{\mathbb C}];Q(z)}^{( \tilde A )}$
and $((W_1\otimes W_2)^*)_{({\mathbb C});Q(z)}^{( \tilde A )}$ are
stable under the component operators $\tau_{Q(z)}(v\otimes t^m)$ of
the operators $Y'_{Q(z)}(v,x)$ for $v\in V$, $m\in {\mathbb Z}$, and
under the operators $L'_{Q(z)}(-1)$, $L'_{Q(z)}(0)$ and
$L'_{Q(z)}(1)$; this uses Proposition \ref{tau-q-a-comp}, Remark
\ref{L'qjpreservesbetaspace}, Propositions \ref{5.1} and
\ref{qz-comm}, and (\ref{qz-sl-2-qz-y-2}).}
\end{rema}

Again let $(W; I)$ ($W$ a generalized $V$-module) be a $Q(z)$-product
of $W_{1}$ and $W_{2}$ and let $w'\in W'$.  Since $I'$ in particular
intertwines the actions of $V\otimes{\mathbb C}[t, t^{-1}]$ and of
$\mathfrak{s}\mathfrak{l}(2)$, and is $\tilde{A}$-compatible, $I'(W')$
is a generalized $V$-module (recall the proof of Proposition
\ref{im-q:abc}).  Thus for every $w'\in W'$, $I'(w')$ also
satisfies the following condition on
\[
\lambda \in (W_1\otimes W_2)^*:
\]
\begin{description}
\item{\bf The $Q(z)$-local grading restriction condition}

(a) The {\em $Q(z)$-grading condition}: $\lambda$ is a (finite) sum of
generalized eigenvectors  for the operator
$L'_{Q(z)}(0)$ on $(W_1\otimes W_2)^*$ that 
are also homogeneous with respect to $\tilde A$, that is, 
\[
\lambda\in ((W_1\otimes W_2)^*)_{[{\mathbb C}]; Q(z)}^{( \tilde A )}.
\]
\label{q-homo}

(b) Let $W_{\lambda;Q(z)}$ be the smallest doubly graded (or
equivalently, $\tilde A$-graded; recall Remark
\ref{q-singleanddoublegraded}) subspace of $((W_1\otimes
W_2)^*)_{[{\mathbb C}];Q(z)}^{(\tilde A )}$ containing $\lambda$ and
stable under the component operators $\tau_{Q(z)}(v\otimes t^m)$ of
the operators $Y'_{Q(z)}(v,x)$ for $v\in V$, $m\in {\mathbb Z}$, and
under the operators $L'_{Q(z)}(-1)$, $L'_{Q(z)}(0)$ and
$L'_{Q(z)}(1)$.  (In view of Remark \ref{q-stableundercomponentops},
$W_{\lambda;Q(z)}$ indeed exists.)  Then $W_{\lambda; Q(z)}$ has the
properties
\begin{eqnarray}
&\dim(W_{\lambda})^{(\beta)}_{[n];Q(z)}<\infty,&\label{q-lgrc1}\\
&(W_{\lambda})^{(\beta)}_{[n+k];Q(z)}=0\;\;\mbox{ for }\;k\in {\mathbb Z}
\;\mbox{ sufficiently negative},&\label{q-lgrc2}
\end{eqnarray}
for any $n\in {\mathbb C}$ and $\beta\in \tilde A$, where as usual the
subscripts denote the ${\mathbb C}$-grading and the superscripts
denote the $\tilde A$-grading.
\end{description}

In the case that $W$ is an (ordinary) $V$-module and $w'\in W'$,
$I'(w')$ also satisfies the following $L(0)$-semisimple version of
this condition on $\lambda \in (W_1\otimes W_2)^*$:

\begin{description}
\item{\bf The $L(0)$-semisimple $Q(z)$-local grading restriction condition}

(a) The {\em $L(0)$-semisimple 
$Q(z)$-grading condition}: $\lambda$ is a (finite) sum of
eigenvectors  for the operator
$L'_{Q(z)}(0)$ on $(W_1\otimes W_2)^*$ that 
are also homogeneous with respect to $\tilde A$, that is,
\[
\lambda\in ((W_1\otimes W_2)^*)_{({\mathbb C});Q(z)}^{( \tilde A )}.
\]
\label{q-semi-homo}

(b) Consider $W_{\lambda;Q(z)}$ as above, which in this case is in fact the
smallest doubly graded (or equivalently, $\tilde A$-graded) subspace
of $((W_1\otimes W_2)^*)_{({\mathbb C});Q(z)}^{( \tilde A )}$ containing
$\lambda$ and stable under the component operators
$\tau_{Q(z)}(v\otimes t^m)$ of the operators $Y'_{Q(z)}(v,x)$ for
$v\in V$, $m\in {\mathbb Z}$, and under the operators $L'_{Q(z)}(-1)$,
$L'_{Q(z)}(0)$ and $L'_{Q(z)}(1)$.  Then $W_{\lambda;Q(z)}$ has the
properties
\begin{eqnarray}
&\dim(W_{\lambda;Q(z)})^{(\beta)}_{(n);Q(z)}<\infty,&\label{q-semi-lgrc1}\\
&(W_\lambda)^{(\beta)}_{(n+k);Q(z)}=0\;\;\mbox{ for }\;k\in {\mathbb Z}
\;\mbox{ sufficiently negative},&\label{q-semi-lgrc2}
\end{eqnarray}
for any $n\in {\mathbb C}$ and $\beta\in \tilde A$, where  the
subscripts denote the ${\mathbb C}$-grading and the superscripts
denote the $\tilde A$-grading.
\end{description}

\begin{nota}
{\rm Note that the set of elements of $(W_1\otimes W_2)^*$ satisfying
either of these two $Q(z)$-local grading restriction conditions, or
either of the Part (a)'s in these conditions, forms a subspace.  We
shall denote the space of elements of $(W_1\otimes W_2)^*$ satisfying
the $Q(z)$-local grading restriction condition and the
$L(0)$-semisimple $Q(z)$-local grading restriction condition by
\[
\lgr_{[\C]; Q(z)}((W_1\otimes W_2)^*)
\]
and 
\[
\lgr_{(\C); Q(z)}((W_1\otimes W_2)^*),
\]
respectively.}
\end{nota}

We have the following important theorems generalizing 
the corresponding results stated in \cite{tensor1} and 
proved in \cite{tensor2}. The proofs of these theorems
will be given in
the next section.

\begin{theo}\label{6.1}
Let $\lambda$ be an element of $(W_1\otimes W_2)^{*}$ satisfying the
$Q(z)$-compatibility condition. Then when acting on $\lambda$, the Jacobi
identity for $Y'_{Q(z)}$ holds, that is,
\begin{eqnarray}
\lefteqn{x_0^{-1}\delta
\left({\displaystyle\frac{x_1-x_2}{x_0}}\right)Y'_{Q(z)}(u, x_1)
Y'_{Q(z)}(v, x_2)\lambda}\nno\\
&&\hspace{2ex}-x_0^{-1} \delta
\left({\displaystyle\frac{x_2-x_1}{-x_0}}\right)Y'_{Q(z)}(v, x_2)
Y'_{Q(z)}(u, x_1)\lambda\nonumber \\
&&=x_2^{-1} \delta
\left({\displaystyle\frac{x_1-x_0}{x_2}}\right)Y'_{Q(z)}(Y(u, x_0)v,
x_2)\lambda
\end{eqnarray}
for $u, v\in V$.
\end{theo}

\begin{theo}\label{6.2}
The subspace $\comp_{Q(z)}((W_1\otimes W_2)^*)$ of $(W_1\otimes W_2)^{*}$
is stable under the operators
$\tau_{Q(z)}(v\otimes t^{n})$ for $v\in V$ and $n\in {\mathbb Z}$, 
and in the M\"obius case,
also under the operators $L'_{Q(z)}(-1)$, $L'_{Q(z)}(0)$ and
$L'_{Q(z)}(1)$;
similarly for the  subspaces $\lgr_{[\C];
Q(z)}((W_1\otimes W_2)^*$ and $\lgr_{(\C); Q(z)}((W_1\otimes W_2)^*$.
\end{theo}

We have:

\begin{theo}\label{q-wk-mod}
The space $\comp_{Q(z)}((W_1\otimes W_2)^*)$, equipped with the vertex
operator map $Y'_{Q(z)}$ and, in case $V$ is M\"{o}bius, also equipped
with the operators $L_{Q(z)}'(-1)$, $L_{Q(z)}'(0)$ and $L_{Q(z)}'(1)$,
is a weak $V$-module; similarly for the spaces
\[
(\comp_{Q(z)}((W_1\otimes W_2)^*))\cap 
(\lgr_{[\C]; Q(z)}((W_1\otimes W_2)^*))
\]
and 
\[
(\comp_{Q(z)}((W_1\otimes W_2)^*))\cap 
(\lgr_{(\C); Q(z)}((W_1\otimes W_2)^*)).
\]
\end{theo}
\pf By Theorem \ref{6.2}, $Y'_{Q(z)}$ is a map {}from the tensor product
of $V$ with any of these three subspaces to the space of formal
Laurent series with elements of the subspace as coefficients.  By
Proposition \ref{5.1} and Theorem \ref{6.1} and, in the case that $V$
is M\"{o}bius, also by Propositions \ref{q-sl-2} and
\ref{qz-l-y-comm}, we see that all the axioms for weak $V$-module are
satisfied.  \epfv

Moreover, we have the following consequence of Theorem \ref{q-wk-mod}
and Lemma \ref{q-a-tilde-comp}, just as in Theorem \ref{generation}:

\begin{theo}\label{q-generation}
Let 
\[
\lambda\in 
\comp_{Q(z)}((W_1\otimes W_2)^*)\cap \lgr_{[\C]; Q(z)}((W_1\otimes W_2)^*).
\] 
Then $W_{\lambda;Q(z)}$ (recall Part (b) of the $Q(z)$-local grading
restriction condition) equipped with the vertex operator map
$Y_{Q(z)}'$ and, in case $V$ is M\"{o}bius, also equipped with the
operators $L'_{Q(z)}(-1)$, $L'_{Q(z)}(0)$ and $L'_{Q(z)}(1)$, is a
(strongly-graded) generalized $V$-module. If in addition
\[
\lambda\in 
\comp_{Q(z)}((W_1\otimes W_2)^*)\cap \lgr_{(\C); Q(z)}((W_1\otimes W_2)^*),
\] 
that is, $\lambda$ is a sum of eigenvectors of $L_{Q(z)}'(0)$, then
$W_{\lambda;Q(z)}$ ($\subset ((W_{1}\otimes
W_{2})^{*})_{(\C);Q(z)}^{(\tilde{A})}$) is a (strongly-graded) $V$-module.
\epfv
\end{theo}

Finally, as in Theorem \ref{characterizationofbackslash}, we can give
an alternative description of $W_1\hboxtr_{Q(z)} W_2$ by
characterizing the elements of $W_1\hboxtr_{Q(z)} W_2$ using the
$Q(z)$-compatibility condition and the $Q(z)$-local grading
restriction conditions, generalizing Theorem 6.3 in
\cite{tensor1}. The proof of the following theorem is the same as that
of Theorem \ref{characterizationofbackslash}.

\begin{theo}\label{q-characterizationofbackslash}
Suppose that for every element 
\[
\lambda
\in \comp_{Q(z)}((W_1\otimes W_2)^*)\cap \lgr_{[\C]; Q(z)}((W_1\otimes W_2)^*)
\]
the (strongly-graded) generalized module $W_{\lambda;Q(z)}$ given in
Theorem \ref{q-generation} is a generalized submodule of some object
of ${\cal C}$ included in $(W_1\otimes W_2)^*$ (this of course holds
in particular if ${\cal C} = {\cal GM}_{sg}$). Then
\[
W_1\hboxtr_{Q(z)}W_2=\comp_{Q(z)}((W_1\otimes W_2)^*)
\cap \lgr_{[\C]; Q(z)}((W_1\otimes W_2)^*).
\]
Suppose that ${\cal C}$ is a category of strongly-graded $V$-modules
(that is, ${\cal C}\subset {\cal M}_{sg}$) and that for every element
\[
\lambda
\in \comp_{Q(z)}((W_1\otimes W_2)^*)\cap \lgr_{(\C); Q(z)}((W_1\otimes W_2)^*)
\]
the (strongly-graded) $V$-module $W_{\lambda;Q(z)}$ given in Theorem
\ref{q-generation} is a submodule of some object of ${\cal C}$
included in $(W_1\otimes W_2)^*$ (which of course holds in particular
if ${\cal C} = {\cal M}_{sg}$).  Then
\[
W_1\hboxtr_{Q(z)}W_2=\comp_{Q(z)}((W_1\otimes W_2)^*)
\cap \lgr_{(\C); Q(z)}((W_1\otimes W_2)^*). \quad\quad\quad\quad
 \square
\]
\end{theo}

\newpage

\setcounter{equation}{0}
\setcounter{rema}{0}

\section{Proof of the theorems used in the constructions}

The primary goal of this section is to prove Theorems \ref{comp=>jcb},
\ref{stable}, \ref{6.1} and \ref{6.2}.  In Section 6.1 we prove
Theorems \ref{comp=>jcb} and \ref{stable}, and in Section 6.2,
Theorems \ref{6.1} and \ref{6.2}.  The proofs in Section 6.1 are
new, even for the category of (ordinary) modules for a vertex operator
algebra satisfying the finiteness and reductivity conditions treated
in \cite{tensor1}--\cite{tensor3}. In \cite{tensor1}--\cite{tensor3},
for a vertex operator algebra satisfying these conditions, Theorems
\ref{6.1} and \ref{6.2}, in the $Q(z)$ case, were proved first, and
then Theorems \ref{comp=>jcb} and \ref{stable}, in the $P(z)$ case,
were proved using results {}from the $Q(z^{-1})$ case and relations
between $P(z)$-tensor products and $Q(z^{-1})$-tensor products. In
Section 6.1, we prove Theorems \ref{comp=>jcb} and \ref{stable}
directly, without using any results {}from the $Q(z)$ case.  As usual,
the reader should observe the justifiability of each step in the
arguments (the well-definedness of the formal series, etc.); again as
usual, this is sometimes quite subtle.

We continue to work in the setting Section 5.  In particular, we have
Assumptions \ref{assum} and \ref{assum-c}, and $z \in \C^{\times}$.

\subsection{Proofs of Theorems \ref{comp=>jcb} and
\ref{stable}}

We first prove a formula for vertex operators that will be needed in
the proofs of both Theorem \ref{comp=>jcb} and Theorem \ref{stable}.

\begin{lemma}
For $u, v \in V$, we have 
\begin{eqnarray}\label{comp=>jcb-9}
&{\displaystyle x_{2}^{-1}\delta\left(\frac{x_{1}-x_{0}}{x_{2}}\right)
Y(e^{x_{1}L(1)}
(-x_{1}^{-2})^{L(0)}u, -x_{0}x_{1}^{-1}x_{2}^{-1})e^{x_{2}L(1)}
(-x_{2}^{-2})^{L(0)}}v &\nn
&{\displaystyle =x_{2}^{-1}\delta\left(\frac{x_{1}-x_{0}}{x_{2}}\right)
e^{x_{2}L(1)}(-x_{2}^{-2})^{L(0)}
Y(u, x_{0})}v. &
\end{eqnarray}
\end{lemma}
\pf Using (\ref{log:p2}), (\ref{log:p3}), (\ref{xe^Lx}) and
(\ref{deltafunctionsubstitutionformula}), we have
\begin{eqnarray*}
\lefteqn{x_{2}^{-1}\delta\left(\frac{x_{1}-x_{0}}{x_{2}}\right)
Y(e^{x_{1}L(1)}
(-x_{1}^{-2})^{L(0)}u, -x_{0}x_{1}^{-1}x_{2}^{-1})e^{x_{2}L(1)}
(-x_{2}^{-2})^{L(0)}}\nn
&&=x_{2}^{-1}\delta\left(\frac{x_{1}-x_{0}}{x_{2}}\right)e^{x_{2}L(1)}
Y(e^{-x_{2}(1-x_{0}x_{1}^{-1})L(1)}(1-x_{0}x_{1}^{-1})^{-2L(0)}\cdot\nn
&&\quad\quad\quad\quad\quad \cdot
e^{x_{1}L(1)}
(-x_{1}^{-2})^{L(0)}u, -x_{0}x_{1}^{-1}x_{2}^{-1}(1-x_{0} x_{1}^{-1})^{-1})
(-x_{2}^{-2})^{L(0)}\nn
&&=x_{2}^{-1}\delta\left(\frac{x_{2}^{-1}-x_{0}x^{-1}_{1}
x_{2}^{-1}}{x^{-1}_{1}}\right)
e^{x_{2}L(1)}(-x_{2}^{-2})^{L(0)}\cdot\nn
&&\quad\quad\cdot Y((-x_{2}^{2})^{L(0)}
e^{-x_{2}(1-x_{0}x_{1}^{-1})L(1)}(1-x_{0}x_{1}^{-1})^{-2L(0)}\cdot\nn
&&\quad\quad\quad\quad\quad\quad\quad\quad\quad\quad \cdot
e^{x_{1}L(1)}
(-x_{1}^{-2})^{L(0)}u, x_{0}x_{1}^{-1}
(x_{2}^{-1}-x_{0} x_{1}^{-1}x_{2}^{-1})^{-1})\nn
&&=x_{2}^{-1}\delta\left(\frac{x_{2}^{-1}-x_{0}x^{-1}_{1}
x_{2}^{-1}}{x^{-1}_{1}}\right)e^{x_{2}L(1)}(-x_{2}^{-2})^{L(0)}\cdot\nn
&&\quad\quad\cdot
Y(e^{x_{2}^{-1}(1-x_{0}x_{1}^{-1})L(1)}
(-x_{2}^{2})^{L(0)}
(1-x_{0}x_{1}^{-1})^{-2L(0)}
e^{x_{1}L(1)}
(-x_{1}^{-2})^{L(0)}u, x_{0})\nn
&&=x_{2}^{-1}\delta\left(\frac{x_{1}-x_{0}}{x_{2}}\right)
e^{x_{2}L(1)}(-x_{2}^{-2})^{L(0)}\cdot\nn
&&\quad\quad \cdot
Y(e^{x_{2}^{-1}(1-x_{0}x_{1}^{-1})L(1)}
(-(x_{2}^{-1}(1-x_{0}x_{1}^{-1}))^{-2})^{L(0)}
e^{x_{1}L(1)}
(-x_{1}^{-2})^{L(0)}u, x_{0})\nn
&&=x_{2}^{-1}\delta\left(\frac{x_{1}-x_{0}}{x_{2}}\right)
e^{x_{2}L(1)}(-x_{2}^{-2})^{L(0)}\cdot\nn
&&\quad\quad \cdot
Y(e^{x_{2}^{-1}(1-x_{0}x_{1}^{-1})L(1)}
e^{-x_{1}x_{2}^{-2}(1-x_{0}x_{1}^{-1})^{2}L(1)}
(-(x_{2}^{-1}(1-x_{0}x_{1}^{-1}))^{-2})^{L(0)}
(-x_{1}^{-2})^{L(0)}u, x_{0})\nn
&&=x_{2}^{-1}\delta\left(\frac{x_{1}-x_{0}}{x_{2}}\right)
e^{x_{2}L(1)}(-x_{2}^{-2})^{L(0)}\cdot\nn
&&\quad\quad \cdot
Y(e^{x_{2}^{-1}(1-x_{0}x_{1}^{-1})
(1-x_{2}^{-1}(x_{1}-x_{0}))L(1)}
(x_{2}^{-1}(x_{1}-x_{0}))^{-2L(0)}
u, x_{0})\nn
&&=x_{2}^{-1}\delta\left(\frac{x_{1}-x_{0}}{x_{2}}\right)
e^{x_{2}L(1)}(-x_{2}^{-2})^{L(0)}
Y(u, x_{0}).
\end{eqnarray*}
\epfv

\noindent {\it Proof of Theorem \ref{comp=>jcb}} Let $\lambda$ be an
element of $(W_{1}\otimes W_{2})^{*}$ satisfying the
$P(z)$-compatibility condition, that is, satisfying (a) the
$P(z)$-lower truncation condition---for all $v\in V$, the formal
Laurent series $Y'_{P(z)}(v, x)\lambda$ involves only finitely many
negative powers of $x$, and (b) formula (\ref{cpb}) for all $v\in V$.

For $u, v\in V$, $w_{(1)}\in W_{1}$ and $w_{(2)}\in W_{2}$, by 
definition, 
\begin{eqnarray}\label{comp=>jcb-1}
\lefteqn{\left(x_{0}^{-1}\delta\left(\frac{x_{1}-x_{2}}{x_{0}}\right)
Y'_{P(z)}(u, x_{1})Y'_{P(z)}(v, x_{2})\lambda\right)(w_{(1)}
\otimes w_{(2)})}\nn
&&=x_{0}^{-1}\delta\left(\frac{x_{1}-x_{2}}{x_{0}}\right)
\Biggl((Y'_{P(z)}(v, x_{2})\lambda)(w_{(1)}\otimes 
Y_{2}^{o}(u, x_{1})w_{(2)})\nn
&&\quad\quad +\res_{y_{1}}z^{-1}
\delta\left(\frac{x_{1}^{-1}-y_{1}}{z}\right)(Y'_{P(z)}(v, x_{2})\lambda)
(Y_{1}(e^{x_{1}L(1)}(-x_{1}^{-2})^{L(0)}u, y_{1})w_{(1)}\otimes w_{(2)})
\Biggr)\nn
&&=x_{0}^{-1}\delta\left(\frac{x_{1}-x_{2}}{x_{0}}\right)
\Biggl(\lambda(w_{(1)}\otimes 
Y_{2}^{o}(v, x_{2})Y_{2}^{o}(u, x_{1})w_{(2)})\nn
&&\quad\quad +\res_{y_{2}}z^{-1}
\delta\left(\frac{x_{2}^{-1}-y_{2}}{z}\right)\lambda(
Y_{1}(e^{x_{2}L(1)}(-x_{2}^{-2})^{L(0)}v, y_{2})w_{(1)}\otimes 
Y_{2}^{o}(u, x_{1})w_{(2)})\nn
&&\quad\quad +\res_{y_{1}}z^{-1}
\delta\left(\frac{x_{1}^{-1}-y_{1}}{z}\right)(Y'_{P(z)}(v, x_{2})\lambda)
(Y_{1}(e^{x_{1}L(1)}(-x_{1}^{-2})^{L(0)}u, y_{1})w_{(1)}\otimes w_{(2)})
\Biggr)\nn
&&=x_{0}^{-1}\delta\left(\frac{x_{1}-x_{2}}{x_{0}}\right)
\lambda(w_{(1)}\otimes 
Y_{2}^{o}(v, x_{2})Y_{2}^{o}(u, x_{1})w_{(2)})\nn
&&\quad\quad +x_{0}^{-1}\delta\left(\frac{x_{1}-x_{2}}{x_{0}}\right)
\res_{y_{2}}z^{-1}
\delta\left(\frac{x_{2}^{-1}-y_{2}}{z}\right)\cdot\nn
&&\quad\quad\quad\quad\quad \cdot\lambda(
Y_{1}(e^{x_{2}L(1)}(-x_{2}^{-2})^{L(0)}v, y_{2})w_{(1)}\otimes 
Y_{2}^{o}(u, x_{1})w_{(2)})\nn
&&\quad\quad +x_{0}^{-1}\delta\left(\frac{x_{1}-x_{2}}{x_{0}}\right)
\res_{y_{1}}z^{-1}
\delta\left(\frac{x_{1}^{-1}-y_{1}}{z}\right)\cdot\nn
&&\quad\quad\quad \quad\quad\cdot(Y'_{P(z)}(v, x_{2})\lambda)
(Y_{1}(e^{x_{1}L(1)}(-x_{1}^{-2})^{L(0)}u, y_{1})w_{(1)}\otimes w_{(2)}).
\end{eqnarray}

Using (\ref{2termdeltarelation}) and (\ref{cpb}), we see that the
third term on the right-hand side of (\ref{comp=>jcb-1}) is equal to
\begin{eqnarray}\label{comp=>jcb-2}
\lefteqn{x_{0}^{-1}\delta\left(\frac{x^{-1}_{2}-x^{-1}_{1}}
{x_{0}x_{1}^{-1}x_{2}^{-1}}\right)
\res_{y_{1}}z^{-1}
\delta\left(\frac{x_{1}^{-1}-y_{1}}{z}\right)\cdot}
\nn
&&\quad\quad\quad\cdot (Y'_{P(z)}(v, x_{2})\lambda)
(Y_{1}(e^{x_{1}L(1)}(-x_{1}^{-2})^{L(0)}u, y_{1})w_{(1)}\otimes w_{(2)})\nn
&&=\res_{y_{1}}x_{1}^{-1}x_{2}^{-1}(x_{0}x_{1}^{-1}x_{2}^{-1})^{-1}
\delta\left(\frac{x^{-1}_{2}-y_{1}-z}
{x_{0}x_{1}^{-1}x_{2}^{-1}}\right)
z^{-1}
\delta\left(\frac{x_{1}^{-1}-y_{1}}{z}\right)\cdot
\nn
&&\quad\quad\quad\cdot (Y'_{P(z)}(v, x_{2})\lambda)
(Y_{1}(e^{x_{1}L(1)}(-x_{1}^{-2})^{L(0)}u, y_{1})w_{(1)}\otimes w_{(2)})\nn
&&=\res_{y_{1}}x_{1}^{-1}x_{2}^{-1}(x_{0}x_{1}^{-1}x_{2}^{-1}+y_{1})^{-1}
\delta\left(\frac{x^{-1}_{2}-z}
{x_{0}x_{1}^{-1}x_{2}^{-1}+y_{1}}\right)z^{-1}
\delta\left(\frac{x_{1}^{-1}-y_{1}}{z}\right)\cdot
\nn
&&\quad\quad\quad\cdot (Y'_{P(z)}(v, x_{2})\lambda)
(Y_{1}(e^{x_{1}L(1)}(-x_{1}^{-2})^{L(0)}u, y_{1})w_{(1)}\otimes w_{(2)})\nn
&&=\res_{y_{1}}x_{1}^{-1}x_{2}^{-1}z^{-1}
\delta\left(\frac{x_{1}^{-1}-y_{1}}{z}\right)\cdot
\nn
&&\quad\quad\quad\cdot 
\left(\tau_{P(z)}\left((x_{0}x_{1}^{-1}x_{2}^{-1}+y_{1})^{-1}
\delta\left(\frac{x^{-1}_{2}-z}
{x_{0}x_{1}^{-1}x_{2}^{-1}+y_{1}}\right)Y_{t}(v, x_{2})\right)
\lambda\right)\nn
&&\quad\quad\quad\quad\quad\quad(Y_{1}(e^{x_{1}L(1)}
(-x_{1}^{-2})^{L(0)}u, y_{1})w_{(1)}\otimes w_{(2)})\nn
&&=\res_{y_{1}}x_{1}^{-1}x_{2}^{-1}
z^{-1}
\delta\left(\frac{x_{1}^{-1}-y_{1}}{z}\right)\cdot
\nn
&&\quad\quad\cdot \Biggl(z^{-1}
\delta\left(\frac{x_{2}^{-1}-x_{0}x_{1}^{-1}x_{2}^{-1}-y_{1}}{z}\right)
\cdot\nn
&&\quad\quad\quad\; \cdot\lambda(Y_{1}(e^{x_{2}L(1)}
(-x_{2}^{-2})^{L(0)}v, x_{0}x_{1}^{-1}x_{2}^{-1}+y_{1})Y_{1}(e^{x_{1}L(1)}
(-x_{1}^{-2})^{L(0)}u, y_{1})w_{(1)}\otimes w_{(2)})\nn
&&\quad\quad\;\; +
(x_{0}x_{1}^{-1}x_{2}^{-1}+y_{1})^{-1}
\delta\left(\frac{z-x_{2}^{-1}}{-x_{0}x_{1}^{-1}x_{2}^{-1}-y_{1}}\right)
\cdot\nn
&&\quad\quad\quad \;\cdot\lambda(Y_{1}(e^{x_{1}L(1)}
(-x_{1}^{-2})^{L(0)}u, y_{1})w_{(1)}\otimes Y_{2}^{o}(v, x_{2})
w_{(2)})\Biggr)\nn
&&=\res_{y_{1}}x_{2}^{-1}
\delta\left(\frac{z+y_{1}}{x_{1}^{-1}}\right) (z+y_{1})^{-1}
\delta\left(\frac{x_{2}^{-1}-x_{0}x_{1}^{-1}x_{2}^{-1}}{z+y_{1}}\right)
\cdot\nn
&&\quad\quad\quad\; \cdot\lambda(Y_{1}(e^{x_{2}L(1)}
(-x_{2}^{-2})^{L(0)}v, x_{0}x_{1}^{-1}x_{2}^{-1}+y_{1})Y_{1}(e^{x_{1}L(1)}
(-x_{1}^{-2})^{L(0)}u, y_{1})w_{(1)}\otimes w_{(2)})\nn
&&\quad\quad +\res_{y_{1}}x_{2}^{-1}
\delta\left(\frac{z+y_{1}}{x_{1}^{-1}}\right)
(x_{0}x_{1}^{-1}x_{2}^{-1})^{-1}
\delta\left(\frac{z+y_{1}-x_{2}^{-1}}{-x_{0}x_{1}^{-1}x_{2}^{-1}}\right)
\cdot\nn
&&\quad\quad\quad\; \cdot\lambda(Y_{1}(e^{x_{1}L(1)}
(-x_{1}^{-2})^{L(0)}u, y_{1})w_{(1)}\otimes Y_{2}^{o}(v, x_{2})
w_{(2)}).
\end{eqnarray}
By (\ref{deltafunctionsubstitutionformula}) and
(\ref{2termdeltarelation}), the right-hand side of (\ref{comp=>jcb-2})
is equal to
\begin{eqnarray}\label{comp=>jcb-3}
\lefteqn{\res_{y_{1}}x_{2}^{-1}
\delta\left(\frac{z+y_{1}}{x_{1}^{-1}}\right)x_{1}
\delta\left(\frac{x_{2}^{-1}-x_{0}x_{1}^{-1}x_{2}^{-1}}{x_{1}^{-1}}\right)
\cdot}\nn
&&\quad\quad \cdot\lambda(Y_{1}(e^{x_{2}L(1)}
(-x_{2}^{-2})^{L(0)}v, x_{0}x_{1}^{-1}x_{2}^{-1}+y_{1})Y_{1}(e^{x_{1}L(1)}
(-x_{1}^{-2})^{L(0)}u, y_{1})w_{(1)}\otimes w_{(2)})\nn
&&\quad +\res_{y_{1}}x_{2}^{-1}
\delta\left(\frac{z+y_{1}}{x_{1}^{-1}}\right)
(x_{0}x_{1}^{-1}x_{2}^{-1})^{-1}
\delta\left(\frac{x_{1}^{-1}-x_{2}^{-1}}{-x_{0}x_{1}^{-1}x_{2}^{-1}}\right)
\cdot\nn
&&\quad\quad \cdot\lambda(Y_{1}(e^{x_{1}L(1)}
(-x_{1}^{-2})^{L(0)}u, y_{1})w_{(1)}\otimes Y_{2}^{o}(v, x_{2})
w_{(2)})\nn
&&=\res_{y_{1}}z^{-1}
\delta\left(\frac{x_{1}^{-1}-y_{1}}{z}\right)x_{2}^{-1}
\delta\left(\frac{x_{1}-x_{0}}{x_{2}}\right)
\cdot\nn
&&\quad\quad \cdot\lambda(Y_{1}(e^{x_{2}L(1)}
(-x_{2}^{-2})^{L(0)}v, x_{0}x_{1}^{-1}x_{2}^{-1}+y_{1})Y_{1}(e^{x_{1}L(1)}
(-x_{1}^{-2})^{L(0)}u, y_{1})w_{(1)}\otimes w_{(2)})\nn
&&\quad +\res_{y_{1}}z^{-1}
\delta\left(\frac{x_{1}^{-1}-y_{1}}{z}\right)
x_{0}^{-1}
\delta\left(\frac{x_{2}-x_{1}}{-x_{0}}\right)
\cdot\nn
&&\quad\quad \cdot\lambda(Y_{1}(e^{x_{1}L(1)}
(-x_{1}^{-2})^{L(0)}u, y_{1})w_{(1)}\otimes Y_{2}^{o}(v, x_{2})
w_{(2)}).
\end{eqnarray}

Since 
\[
\res_{y_{2}}y_{2}^{-1}
\delta\left(\frac{x_{0}x_{1}^{-1}x_{2}^{-1}+y_{1}}{y_{2}}\right)=1,
\]
the first term on the right-hand side of (\ref{comp=>jcb-3}) 
can be written as
\begin{eqnarray}\label{comp=>jcb-4}
\lefteqn{\res_{y_{1}}z^{-1}
\delta\left(\frac{x_{1}^{-1}-y_{1}}{z}\right)
x_{2}^{-1}
\delta\left(\frac{x_{1}-x_{0}}{x_{2}}\right)
\res_{y_{2}}y_{2}^{-1}
\delta\left(\frac{x_{0}x_{1}^{-1}x_{2}^{-1}+y_{1}}{y_{2}}\right)
\cdot }\nn
&&\quad\quad\quad \cdot\lambda(Y_{1}(e^{x_{2}L(1)}
(-x_{2}^{-2})^{L(0)}v, x_{0}x_{1}^{-1}x_{2}^{-1}+y_{1})Y_{1}(e^{x_{1}L(1)}
(-x_{1}^{-2})^{L(0)}u, y_{1})w_{(1)}\otimes w_{(2)})\nn
&&=x_{2}^{-1}
\delta\left(\frac{x_{1}-x_{0}}{x_{2}}\right)\res_{y_{1}}
\res_{y_{2}}z^{-1}
\delta\left(\frac{x_{1}^{-1}-y_{1}}{z}\right)
(x_{0}x_{1}^{-1}x_{2}^{-1})^{-1}
\delta\left(\frac{y_{2}-y_{1}}{x_{0}x_{1}^{-1}x_{2}^{-1}}\right)
\cdot \nn
&&\quad\quad\quad \cdot\lambda(Y_{1}(e^{x_{2}L(1)}
(-x_{2}^{-2})^{L(0)}v, y_{2})Y_{1}(e^{x_{1}L(1)}
(-x_{1}^{-2})^{L(0)}u, y_{1})w_{(1)}\otimes w_{(2)}),
\end{eqnarray}
where we have also used (\ref{deltafunctionsubstitutionformula}) and
(\ref{2termdeltarelation}).  Again using (\ref{2termdeltarelation})
and (\ref{deltafunctionsubstitutionformula}), we see that the
right-hand side of (\ref{comp=>jcb-4}) is also equal to
\begin{eqnarray}\label{comp=>jcb-5}
\lefteqn{x_{2}^{-1}
\delta\left(\frac{x_{2}^{-1}-x_{0}x_{1}^{-1}x_{2}^{-1}}
{x_{1}^{-1}}\right)\res_{y_{1}}
\res_{y_{2}}z^{-1}
\delta\left(\frac{x_{1}^{-1}-y_{1}}{z}\right)
y_{2}^{-1}
\delta\left(\frac{x_{0}x_{1}^{-1}x_{2}^{-1}+y_{1}}{y_{2}}\right)
\cdot} \nn
&&\quad\quad\quad\quad \cdot\lambda(Y_{1}(e^{x_{2}L(1)}
(-x_{2}^{-2})^{L(0)}v, y_{2})Y_{1}(e^{x_{1}L(1)}
(-x_{1}^{-2})^{L(0)}u, y_{1})w_{(1)}\otimes w_{(2)})\nn
&&=x_{2}^{-1}
\delta\left(\frac{x_{2}^{-1}-x_{0}x_{1}^{-1}x_{2}^{-1}}
{x_{1}^{-1}}\right)\res_{y_{1}}
\res_{y_{2}}z^{-1}
\delta\left(\frac{x_{2}^{-1}-x_{0}x_{1}^{-1}x_{2}^{-1}-y_{1}}{z}\right)
\cdot\nn
&&\quad\quad\quad\quad \cdot y_{2}^{-1}
\delta\left(\frac{x_{0}x_{1}^{-1}x_{2}^{-1}+y_{1}}{y_{2}}\right)
\cdot \nn
&&\quad\quad\quad\quad \cdot\lambda(Y_{1}(e^{x_{2}L(1)}
(-x_{2}^{-2})^{L(0)}v, y_{2})Y_{1}(e^{x_{1}L(1)}
(-x_{1}^{-2})^{L(0)}u, y_{1})w_{(1)}\otimes w_{(2)})\nn
&&=x_{2}^{-1}
\delta\left(\frac{x_{1}-x_{0}}
{x_{2}}\right)\res_{y_{1}}
\res_{y_{2}}z^{-1}
\delta\left(\frac{x_{2}^{-1}-y_{2}}{z}\right)
y_{2}^{-1}
\delta\left(\frac{x_{0}x_{1}^{-1}x_{2}^{-1}+y_{1}}{y_{2}}\right)
\cdot \nn
&&\quad\quad\quad\quad \cdot\lambda(Y_{1}(e^{x_{2}L(1)}
(-x_{2}^{-2})^{L(0)}v, y_{2})Y_{1}(e^{x_{1}L(1)}
(-x_{1}^{-2})^{L(0)}u, y_{1})w_{(1)}\otimes w_{(2)})\nn
&&=x_{2}^{-1}
\delta\left(\frac{x_{1}-x_{0}}
{x_{2}}\right)\res_{y_{1}}
\res_{y_{2}}z^{-1}
\delta\left(\frac{x_{2}^{-1}-y_{2}}{z}\right)
(x_{0}x_{1}^{-1}x_{2}^{-1})^{-1}
\delta\left(\frac{y_{2}-y_{1}}{x_{0}x_{1}^{-1}x_{2}^{-1}}\right)
\cdot \nn
&&\quad\quad\quad\quad \cdot\lambda(Y_{1}(e^{x_{2}L(1)}
(-x_{2}^{-2})^{L(0)}v, y_{2})Y_{1}(e^{x_{1}L(1)}
(-x_{1}^{-2})^{L(0)}u, y_{1})w_{(1)}\otimes w_{(2)}).
\end{eqnarray}
That is, in the middle delta-function expression in the right-hand
side of (\ref{comp=>jcb-4}), we may replace $x_1$ by $x_2$ and $y_1$
by $y_2$.

{}From (\ref{comp=>jcb-1})--(\ref{comp=>jcb-5}), we obtain
\begin{eqnarray}\label{comp=>jcb-6}
\lefteqn{\left(x_{0}^{-1}\delta\left(\frac{x_{1}-x_{2}}{x_{0}}\right)
Y'_{P(z)}(u, x_{1})Y'_{P(z)}(v, x_{2})\lambda\right)(w_{(1)}
\otimes w_{(2)})}\nn
&&=x_{0}^{-1}\delta\left(\frac{x_{1}-x_{2}}{x_{0}}\right)
\lambda(w_{(1)}\otimes 
Y_{2}^{o}(v, x_{2})Y_{2}^{o}(u, x_{1})w_{(2)})\nn
&&\quad +x_{0}^{-1}\delta\left(\frac{x_{1}-x_{2}}{x_{0}}\right)
\res_{y_{2}}z^{-1}
\delta\left(\frac{x_{2}^{-1}-y_{2}}{z}\right)\cdot \nn
&&\quad\quad\quad\quad \cdot\lambda(
Y_{1}(e^{x_{2}L(1)}(-x_{2}^{-2})^{L(0)}v, y_{2})w_{(1)}\otimes 
Y_{2}^{o}(u, x_{1})w_{(2)})\nn
&&\quad +x_{2}^{-1}
\delta\left(\frac{x_{1}-x_{0}}{x_{2}}\right)\res_{y_{1}}
\res_{y_{2}}z^{-1}
\delta\left(\frac{x_{2}^{-1}-y_{2}}{z}\right)
(x_{0}x_{1}^{-1}x_{2}^{-1})^{-1}  
\delta\left(\frac{y_{2}-y_{1}}{x_{0}x_{1}^{-1}x_{2}^{-1}}\right)
\cdot \nn
&&\quad\quad\quad\quad \cdot\lambda(Y_{1}(e^{x_{2}L(1)}
(-x_{2}^{-2})^{L(0)}v, y_{2})Y_{1}(e^{x_{1}L(1)}
(-x_{1}^{-2})^{L(0)}u, y_{1})w_{(1)}\otimes w_{(2)})\nn
&&\quad +x_{0}^{-1}
\delta\left(\frac{x_{2}-x_{1}}{-x_{0}}\right)\res_{y_{1}}
z^{-1}
\delta\left(\frac{x_{1}^{-1}-y_{1}}{z}\right)
\cdot\nn
&&\quad\quad\quad\quad \cdot\lambda(Y_{1}(e^{x_{1}L(1)}
(-x_{1}^{-2})^{L(0)}u, y_{1})w_{(1)}\otimes Y_{2}^{o}(v, x_{2})
w_{(2)})\Biggr).
\end{eqnarray}
{}From (\ref{comp=>jcb-5}) and (\ref{comp=>jcb-6}), replacing $u, v,
x_{1}, x_{2}, x_{0}$ by $v, u, x_{2}, x_{1}, -x_{0}$, respectively,
and also using (\ref{2termdeltarelation}), we find that
\begin{eqnarray}\label{comp=>jcb-7}
\lefteqn{\left(-x_{0}^{-1}\delta\left(\frac{x_{2}-x_{1}}{-x_{0}}\right)
Y'_{P(z)}(v, x_{2})Y'_{P(z)}(u, x_{1})\lambda\right)(w_{(1)}
\otimes w_{(2)})}\nn
&&=-x_{0}^{-1}\delta\left(\frac{x_{2}-x_{1}}{-x_{0}}\right)
\lambda(w_{(1)}\otimes 
Y_{2}^{o}(u, x_{1})Y_{2}^{o}(v, x_{2})w_{(2)})\nn
&&\quad -x_{0}^{-1}\delta\left(\frac{x_{2}-x_{1}}{-x_{0}}\right)
\res_{y_{1}}z^{-1}
\delta\left(\frac{x_{1}^{-1}-y_{1}}{z}\right)\cdot \nn
&&\quad\quad\quad\quad \cdot\lambda(
Y_{1}(e^{x_{1}L(1)}(-x_{1}^{-2})^{L(0)}u, y_{1})w_{(1)}\otimes 
Y_{2}^{o}(v, x_{2})w_{(2)})\nn
&&\quad -x_{2}^{-1}
\delta\left(\frac{x_{1}-x_{0}}{x_{2}}\right)\res_{y_{1}}
\res_{y_{2}}z^{-1}
\delta\left(\frac{x_{2}^{-1}-y_{2}}{z}\right)
(x_{0}x_{1}^{-1}x_{2}^{-1})^{-1}  
\delta\left(\frac{y_{1}-y_{2}}{-x_{0}x_{1}^{-1}x_{2}^{-1}}\right)
\cdot \nn
&&\quad\quad\quad\quad \cdot\lambda(Y_{1}(e^{x_{1}L(1)}
(-x_{1}^{-2})^{L(0)}u, y_{1})Y_{1}(e^{x_{2}L(1)}
(-x_{2}^{-2})^{L(0)}v, y_{2})w_{(1)}\otimes w_{(2)})\nn
&&\quad -x_{0}^{-1}
\delta\left(\frac{x_{1}-x_{2}}{x_{0}}\right)\res_{y_{2}}
z^{-1}
\delta\left(\frac{x_{2}^{-1}-y_{2}}{z}\right)
\cdot\nn
&&\quad\quad\quad\quad \cdot\lambda(Y_{1}(e^{x_{2}L(1)}
(-x_{2}^{-2})^{L(0)}v, y_{2})w_{(1)}\otimes Y_{2}^{o}(u, x_{1})
w_{(2)})\Biggr).
\end{eqnarray}

Using (\ref{comp=>jcb-6}), (\ref{comp=>jcb-7}), the Jacobi identity,
the opposite Jacobi identity (\ref{op-jac-id}) and
(\ref{2termdeltarelation}), we obtain
\begin{eqnarray}\label{comp=>jcb-8}
\lefteqn{\Biggl(x_{0}^{-1}\delta\left(\frac{x_{1}-x_{2}}{x_{0}}\right)
Y'_{P(z)}(u, x_{1})Y'_{P(z)}(v, x_{2})\lambda}\nn
&&\quad
-x_{0}^{-1}\delta\left(\frac{x_{2}-x_{1}}{-x_{0}}\right)
Y'_{P(z)}(v, x_{2})Y'_{P(z)}(u, x_{1})\lambda\Biggr)(w_{(1)}
\otimes w_{(2)})\nn
&&=\lambda\Biggl(w_{(1)}\otimes 
\Biggl(x_{0}^{-1}\delta\left(\frac{x_{1}-x_{2}}{x_{0}}\right)
Y_{2}^{o}(v, x_{2})Y_{2}^{o}(u, x_{1})\nn
&&\quad\quad\quad\quad \quad\quad \quad
-x_{0}^{-1}\delta\left(\frac{x_{2}-x_{1}}{-x_{0}}\right)
Y_{2}^{o}(u, x_{1})Y_{2}^{o}(v, x_{2})\Biggr) w_{(2)}\Biggr)\nn
&&\quad +x_{2}^{-1}
\delta\left(\frac{x_{1}-x_{0}}{x_{2}}\right)\res_{y_{1}}
\res_{y_{2}}z^{-1}
\delta\left(\frac{x_{2}^{-1}-y_{2}}{z}\right)
\cdot \nn
&&\quad\quad \cdot\lambda\Biggl(
\Biggl((x_{0}x_{1}^{-1}x_{2}^{-1})^{-1}
\delta\left(\frac{y_{2}-y_{1}}{x_{0}x_{1}^{-1}x_{2}^{-1}}\right)\cdot\nn
&&\quad\quad\quad\quad\quad\quad\quad \cdot Y_{1}(e^{x_{2}L(1)}
(-x_{2}^{-2})^{L(0)}v, y_{2})Y_{1}(e^{x_{1}L(1)}
(-x_{1}^{-2})^{L(0)}u, y_{1})\nn
&&\quad\quad\quad\quad\;\;-(x_{0}x_{1}^{-1}x_{2}^{-1})^{-1}
\delta\left(\frac{y_{1}-y_{2}}{-x_{0}x_{1}^{-1}x_{2}^{-1}}\right)\cdot\nn
&&\quad\quad\quad\quad\quad\quad \quad\cdot
Y_{1}(e^{x_{1}L(1)}
(-x_{1}^{-2})^{L(0)}u, y_{1})Y_{1}(e^{x_{2}L(1)}
(-x_{2}^{-2})^{L(0)}v, y_{2})\Biggr)w_{(1)}\otimes w_{(2)}\Biggr)\nn
&&=\lambda\Biggl(w_{(1)}\otimes 
\Biggl(x_{2}^{-1}\delta\left(\frac{x_{1}-x_{0}}{x_{2}}\right)
Y_{2}^{o}(Y(u, x_{0})v, x_{2})\Biggr) w_{(2)}\Biggr)\nn
&&\quad +x_{2}^{-1}
\delta\left(\frac{x_{1}-x_{0}}{x_{2}}\right)\res_{y_{1}}
\res_{y_{2}}z^{-1}
\delta\left(\frac{x_{2}^{-1}-y_{2}}{z}\right)
\cdot \nn
&&\quad\quad \cdot\lambda\Biggl(
\Biggl(y_{2}^{-1}
\delta\left(\frac{y_{1}+x_{0}x_{1}^{-1}x_{2}^{-1}}{y_{2}}\right)\cdot\nn
&&\quad\quad\quad\quad \cdot Y_{1}(Y(e^{x_{1}L(1)}
(-x_{1}^{-2})^{L(0)}u, -x_{0}x_{1}^{-1}x_{2}^{-1})e^{x_{2}L(1)}
(-x_{2}^{-2})^{L(0)}v, y_{2})\Biggr)w_{(1)}\otimes w_{(2)}\Biggr)\nn
&&=x_{2}^{-1}\delta\left(\frac{x_{1}-x_{0}}{x_{2}}\right)
\lambda(w_{(1)}\otimes Y_{2}^{o}(Y(u, x_{0})v, x_{2}) w_{(2)})\nn
&&\quad +x_{2}^{-1}
\delta\left(\frac{x_{1}-x_{0}}{x_{2}}\right)\res_{y_{1}}
\res_{y_{2}}z^{-1}
\delta\left(\frac{x_{2}^{-1}-y_{2}}{z}\right)
\cdot \nn
&&\quad\quad \cdot\lambda\Biggl(
\Biggl(y_{1}^{-1}
\delta\left(\frac{y_{2}-x_{0}x_{1}^{-1}x_{2}^{-1}}{y_{1}}\right)\cdot\nn
&&\quad\quad\quad\quad \cdot Y_{1}(Y(e^{x_{1}L(1)}
(-x_{1}^{-2})^{L(0)}u, -x_{0}x_{1}^{-1}x_{2}^{-1})e^{x_{2}L(1)}
(-x_{2}^{-2})^{L(0)}v, y_{2})\Biggr)w_{(1)}\otimes w_{(2)}\Biggr)\nn
&&=x_{2}^{-1}\delta\left(\frac{x_{1}-x_{0}}{x_{2}}\right)
\lambda(w_{(1)}\otimes Y_{2}^{o}(Y(u, x_{0})v, x_{2}) w_{(2)})\nn
&&\quad +x_{2}^{-1}
\delta\left(\frac{x_{1}-x_{0}}{x_{2}}\right)
\res_{y_{2}}z^{-1}
\delta\left(\frac{x_{2}^{-1}-y_{2}}{z}\right)
\cdot \nn
&&\quad\quad \cdot\lambda(
Y_{1}(Y(e^{x_{1}L(1)}
(-x_{1}^{-2})^{L(0)}u, -x_{0}x_{1}^{-1}x_{2}^{-1})e^{x_{2}L(1)}
(-x_{2}^{-2})^{L(0)}v, y_{2})w_{(1)}\otimes w_{(2)})\nn
\end{eqnarray}

Finally, {}from (\ref{comp=>jcb-9}) we see that the right-hand side of
(\ref{comp=>jcb-8}) becomes
\begin{eqnarray}\label{comp=>jcb-10}
\lefteqn{x_{2}^{-1}\delta\left(\frac{x_{1}-x_{0}}{x_{2}}\right)
\lambda(w_{(1)}\otimes Y_{2}^{o}(Y(u, x_{0})v, x_{2}) w_{(2)})}\nn
&&\quad +x_{2}^{-1}
\delta\left(\frac{x_{1}-x_{0}}{x_{2}}\right)
\res_{y_{2}}z^{-1}
\delta\left(\frac{x_{2}^{-1}-y_{2}}{z}\right)
\cdot \nn
&&\quad\quad \cdot\lambda(
Y_{1}(e^{x_{2}L(1)}(-x_{2}^{-2})^{L(0)}
Y(u, x_{0})
v, y_{2})w_{(1)}\otimes w_{(2)})\nn
&&=\left(x_{2}^{-1}\delta\left(\frac{x_{1}-x_{0}}{x_{2}}\right)
Y'_{P(z)}(Y(u, x_{0})v, x_{2})\lambda\right)(w_{(1)}\otimes w_{(2)}),
\end{eqnarray}
and we have proved the Jacobi identity and hence Theorem
\ref{comp=>jcb}.
\epfv

\noindent {\it Proof of Theorem \ref{stable}} Let $\lambda$ be an
element of $(W_{1}\otimes W_{2})^{*}$ satisfying the
$P(z)$-compatibility condition. We first want to prove that the
coefficient of each power of $x$ in $Y'_{P(z)}(u, x_0)Y'_{P(z)}(v,
x)\lambda$ is a formal Laurent series involving only finitely many
negative powers of $x_0$ and that
\begin{eqnarray}\label{stable-1}
\lefteqn{\tau_{P(z)}\left(x_{0}^{-1}\delta\left(\frac{x_{1}^{-1}-z}
{x_{0}}\right)
Y_t(u, x_1)\right)
Y'_{P(z)}(v, x)\lambda}\nno\\
&&=x_{0}^{-1}\delta\left(\frac{x_{1}^{-1}-z}
{x_{0}}\right)
Y'_{P(z)}(u, x_1)Y'_{P(z)}(v,
x) \lambda
\end{eqnarray}
for all $u, v\in V$.  Using the commutator formula (Proposition
\ref{pz-comm}) for $Y'_{P(z)}$, we have
\begin{eqnarray}\label{stable-2}
\lefteqn{Y'_{P(z)}(u, x_0)Y'_{P(z)}(v, x)\lambda}\nno\\
&&=Y'_{P(z)}(v, x)Y'_{P(z)}(u, x_0)\lambda\nno\\
&&\quad -\res_{y}x^{-1}_0\delta\left(\frac{x-y}{x_0}\right)
Y'_{P(z)}(Y(v, y)u, x_0)\lambda.
\end{eqnarray}
Each coefficient in $x$ of the right-hand side of (\ref{stable-2}) is
a formal Laurent series involving only finitely many negative powers
of $x_0$ since $\lambda$ satisfies the $P(z)$-lower truncation
condition.  Thus the coefficients in $x$ of $Y'_{P(z)}(v, x)\lambda$
satisfy the $P(z)$-lower truncation condition.

By (\ref{taudef}) and (\ref{Y'def}), we have
\begin{eqnarray}\label{stable-3}
\lefteqn{\left(\tau_{P(z)}\left(x_{0}^{-1}\delta\left(\frac{x_{1}^{-1}-z}
{x_{0}}\right)
Y_t(u, x_1)\right)
Y'_{P(z)}(v, x) \lambda\right)(w_{(1)}\otimes w_{(2)})}\nno\\
&&=z^{-1}\delta\left(\frac{x_1^{-1}-x_{0}}{z}\right)(Y'_{P(z)}(v, x)
\lambda)(Y_1(e^{x_{1}L(1)}(-x_{1}^{-2})^{L(0)}u, x_0)w_{(1)}
\otimes w_{(2)})\nno\\
&&\quad +x^{-1}_0\delta\left(\frac{z-x_1^{-1}}{-x_0}\right)(Y'_{P(z)}(v, x)
\lambda)(w_{(1)}\otimes Y_2^{o}(u, x_1)w_{(2)})\nno\\
&&=z^{-1}\delta\left(\frac{x_1^{-1}-x_{0}}{z}\right)
\Biggl(\lambda(Y_1(e^{x_{1}L(1)}(-x_{1}^{-2})^{L(0)}u, x_0)w_{(1)}
\otimes Y_{2}^{o}(v, x)w_{(2)})\nno\\
&&\quad \quad +\res_{x_2}
z^{-1}\delta\left(\frac{x^{-1}-x_{2}}{z}\right)\cdot\nn
&&\quad \quad\quad \quad\cdot \lambda(Y_1(e^{xL(1)}(-x^{-2})^{L(0)}v, x_2)
Y_1(e^{x_{1}L(1)}(-x_{1}^{-2})^{L(0)}u,
x_0)w_{(1)}\otimes w_{(2)})\Biggr)\nno\\
&&\quad +x^{-1}_0\delta\left(\frac{z-x_1^{-1}}{-x_0}\right)
\Biggl(\lambda(w_{(1)}\otimes Y_2^o(v, x)Y_2^{o}(u, x_1)w_{(2)})\nno\\
&&\quad \quad +\res_{x_2} z^{-1}\delta\left(\frac{x^{-1}-x_{2}}{z}\right)
\lambda(Y_1(e^{xL(1)}(-x^{-2})^{L(0)}v, x_2)w_{(1)}\otimes
Y_2^{o}(u, x_1)w_{(2)})\Biggr).\nn
&&
\end{eqnarray}
Now the distributive law applies, giving us four terms. Then using the
opposite commutator formula for $Y^{o}_{2}$ (recall (\ref{op-jac-id}))
and the commutator formula for $Y_{2}$, and (\ref{taudef}), we can
write the right-hand side of (\ref{stable-3}) as
\begin{eqnarray}\label{stable-4}
\lefteqn{z^{-1}\delta\left(\frac{x_1^{-1}-x_{0}}{z}\right)\lambda
(Y_1(e^{x_{1}L(1)}(-x_{1}^{-2})^{L(0)}u, x_0)w_{(1)}
\otimes Y_{2}^{o}(v, x)w_{(2)})}\nno\\
&&\quad  +z^{-1}\delta\left(\frac{x_1^{-1}-x_{0}}{z}\right)\res_{x_2}
z^{-1}\delta\left(\frac{x^{-1}-x_{2}}{z}\right)\cdot\nn
&&\quad \quad\cdot \lambda(Y_1(e^{x_{1}L(1)}(-x_{1}^{-2})^{L(0)}u,
x_0)Y_1(e^{xL(1)}(-x^{-2})^{L(0)}v, x_2)w_{(1)}\otimes w_{(2)})\nno\\
&&\quad  +z^{-1}\delta\left(\frac{x_1^{-1}-x_{0}}{z}\right)\res_{x_2}
z^{-1}\delta\left(\frac{x^{-1}-x_{2}}{z}\right)\res_{x_{3}}
x_{0}^{-1}\delta\left(\frac{x_{2}-x_{3}}{x_{0}}\right)\cdot\nn
&&\quad \quad \cdot 
\lambda(Y_1(Y(e^{xL(1)}(-x^{-2})^{L(0)}v, x_3)
e^{x_{1}L(1)}(-x_{1}^{-2})^{L(0)}u,
x_0)w_{(1)}\otimes w_{(2)})\nn
&&\quad +x^{-1}_0\delta\left(\frac{z-x_1^{-1}}{-x_0}\right)
\lambda(w_{(1)}\otimes Y_2^{o}(u, x_1)Y_2^o(v, x)w_{(2)})\nn
&&\quad -x^{-1}_0\delta\left(\frac{z-x_1^{-1}}{-x_0}\right)
\res_{x_{3}}x^{-1}_1\delta\left(\frac{x-x_3}{x_1}\right)
\lambda(w_{(1)}\otimes Y_2^{o}(Y(v, x_{3})u, x_{1})w_{(2)})\nn
&&\quad  +x^{-1}_0\delta\left(\frac{z-x_1^{-1}}{-x_0}\right)
\res_{x_2} z^{-1}\delta\left(\frac{x^{-1}-x_{2}}{z}\right)\cdot\nn
&&\quad \quad \cdot 
\lambda(Y_1(e^{xL(1)}(-x^{-2})^{L(0)}v, x_2)w_{(1)}\otimes
Y_2^{o}(u, x_{1})w_{(2)})\nn
&&=\left(\tau_{P(z)}\left(x_{0}^{-1}\delta\left(\frac{x_{1}^{-1}-z}
{x_{0}}\right)
Y_t(u, x_1)\right)\lambda\right)(w_{(1)}\otimes
Y_2^{o}(v, x)w_{(2)})\nn
&&\quad +\res_{x_2}
z^{-1}\delta\left(\frac{x^{-1}-x_{2}}{z}\right)\cdot\nn
&&\quad \quad\cdot \left(\tau_{P(z)}
\left(x_{0}^{-1}\delta\left(\frac{x_{1}^{-1}-z}
{x_{0}}\right)
Y_t(u, x_1)\right)\lambda\right)
(Y_1(e^{xL(1)}(-x^{-2})^{L(0)}v, x_2)w_{(1)}\otimes w_{(2)})\nn
&&\quad  +\res_{x_2}
z^{-1}\delta\left(\frac{x^{-1}-x_{2}}{z}\right)\res_{x_{3}}
x_{0}^{-1}\delta\left(\frac{x_{2}-x_{3}}{x_{0}}\right)
z^{-1}\delta\left(\frac{x_1^{-1}-x_{0}}{z}\right)\cdot\nn
&&\quad \quad \cdot 
\lambda(Y_1(Y(e^{xL(1)}(-x^{-2})^{L(0)}v, x_3)
e^{x_{1}L(1)}(-x_{1}^{-2})^{L(0)}u,
x_0)w_{(1)}\otimes w_{(2)})\nn
&&\quad -\res_{x_{3}}x^{-1}_1\delta\left(\frac{x-x_3}{x_1}\right)
x^{-1}_0\delta\left(\frac{z-x_1^{-1}}{-x_0}\right)
\lambda(w_{(1)}\otimes Y_2^{o}(Y(v, x_{3})u, x_{1})w_{(2)}).
\end{eqnarray}

Since $\lambda$ satisfies the $P(z)$-compatibility condition
(\ref{cpb}), by (\ref{Y'def}) the sum of the first two terms of
(\ref{stable-4}) is equal to
\begin{eqnarray}\label{stable-5}
\lefteqn{\left(Y'_{P(z)}(v, x)
\tau_{P(z)}\left(x_{0}^{-1}\delta\left(\frac{x_{1}^{-1}-z}
{x_{0}}\right)
Y_t(u, x_1)\right)\lambda\right)(w_{(1)}\otimes
w_{(2)})}\nn
&&=x_{0}^{-1}\delta\left(\frac{x_{1}^{-1}-z}
{x_{0}}\right)\left(Y'_{P(z)}(v, x)
Y'_{P(z)}(u, x_1)\lambda\right)(w_{(1)}\otimes
w_{(2)}).
\end{eqnarray}

Changing the dummy variable $x_{3}$ to $-x_{3}x^{-1}x_{1}^{-1}$ where
we use $x_{3}$ to denote the new dummy variable, using
(\ref{deltafunctionsubstitutionformula}), (\ref{2termdeltarelation})
and (\ref{comp=>jcb-9}), and then evaluating $\res_{x_{2}}$, we see
that the third term of (\ref{stable-4}) is equal to
\begin{eqnarray}\label{stable-6}
\lefteqn{-\res_{x_2}
\res_{x_{3}}
z^{-1}\delta\left(\frac{x^{-1}-x_{2}}{z}\right)x^{-1}x_{1}^{-1}
x_{0}^{-1}\delta\left(\frac{x_{2}+x_{3}x^{-1}x_{1}^{-1}}
{x_{0}}\right)z^{-1}\delta\left(\frac{x_1^{-1}-x_{0}}{z}\right)\cdot}\nn
&&\quad \quad \cdot 
\lambda(Y_1(Y(e^{xL(1)}(-x^{-2})^{L(0)}v, -x_3x^{-1}x_{1}^{-1})
e^{x_{1}L(1)}(-x_{1}^{-2})^{L(0)}u,
x_0)w_{(1)}\otimes w_{(2)})\nn
&&=-\res_{x_2}
z^{-1}\delta\left(\frac{x^{-1}-x_{2}}{z}\right)\res_{x_{3}}x^{-1}x_{1}^{-1}
x_{2}^{-1}\delta\left(\frac{x_{0}-x_{3}x^{-1}x_{1}^{-1}}
{x_{2}}\right)z^{-1}\delta\left(\frac{x_1^{-1}-x_{0}}{z}\right)\cdot\nn
&&\quad \quad \cdot 
\lambda(Y_1(Y(e^{xL(1)}(-x^{-2})^{L(0)}v, -x_3x^{-1}x_{1}^{-1})
e^{x_{1}L(1)}(-x_{1}^{-2})^{L(0)}u,
x_0)w_{(1)}\otimes w_{(2)})\nn
&&=-\res_{x_{3}}\res_{x_2}
z^{-1}\delta\left(\frac{x^{-1}-x_{0}+x_{3}x^{-1}x_{1}^{-1}}{z}\right)
\cdot\nn
&&\quad\quad
 \cdot x^{-1}x_{1}^{-1}
x_{2}^{-1}\delta\left(\frac{x_{0}-x_{3}x^{-1}x_{1}^{-1}}
{x_{2}}\right)z^{-1}\delta\left(\frac{x_1^{-1}-x_{0}}{z}\right)\cdot\nn
&&\quad \quad \cdot 
\lambda(Y_1(Y(e^{xL(1)}(-x^{-2})^{L(0)}v, -x_3x^{-1}x_{1}^{-1})
e^{x_{1}L(1)}(-x_{1}^{-2})^{L(0)}u,
x_0)w_{(1)}\otimes w_{(2)})\nn
&&=-\res_{x_{3}}\res_{x_2}
(z+x_{0})^{-1}\delta\left(\frac{x^{-1}+x_{3}x^{-1}x_{1}^{-1}}{z+x_{0}}\right)
\cdot\nn
&&\quad\quad\cdot x^{-1}x_{1}^{-1}
x_{2}^{-1}\delta\left(\frac{x_{0}-x_{3}x^{-1}x_{1}^{-1}}
{x_{2}}\right)x_{1}\delta\left(\frac{z+x_{0}}{x_1^{-1}}\right)\cdot\nn
&&\quad \quad \cdot 
\lambda(Y_1(Y(e^{xL(1)}(-x^{-2})^{L(0)}v, -x_3x^{-1}x_{1}^{-1})
e^{x_{1}L(1)}(-x_{1}^{-2})^{L(0)}u,
x_0)w_{(1)}\otimes w_{(2)})\nn
&&=-\res_{x_{3}}\res_{x_2}
x_{1}\delta\left(\frac{x^{-1}+x_{3}x^{-1}x_{1}^{-1}}
{x_{1}^{-1}}\right)
\cdot\nn
&&\quad\quad\cdot x^{-1}x_{1}^{-1}
x_{2}^{-1}\delta\left(\frac{x_{0}-x_{3}x^{-1}x_{1}^{-1}}
{x_{2}}\right)x_{1}\delta\left(\frac{z+x_{0}}{x_1^{-1}}\right)\cdot\nn
&&\quad \quad \cdot 
\lambda(Y_1(Y(e^{xL(1)}(-x^{-2})^{L(0)}v, -x_3x^{-1}x_{1}^{-1})
e^{x_{1}L(1)}(-x_{1}^{-2})^{L(0)}u,
x_0)w_{(1)}\otimes w_{(2)})\nn
&&=-\res_{x_{3}}\res_{x_2}
x_{1}^{-1}\delta\left(\frac{x-x_{3}}
{x_{1}}\right)
\cdot\nn
&&\quad\quad\cdot 
x_{2}^{-1}\delta\left(\frac{x_{0}-x_{3}x^{-1}x_{1}^{-1}}
{x_{2}}\right)z^{-1}\delta\left(\frac{x_1^{-1}-x_{0}}{z}\right)\cdot\nn
&&\quad \quad \cdot 
\lambda(Y_1(Y(e^{xL(1)}(-x^{-2})^{L(0)}v, -x_3x^{-1}x_{1}^{-1})
e^{x_{1}L(1)}(-x_{1}^{-2})^{L(0)}u,
x_0)w_{(1)}\otimes w_{(2)})\nn
&&=-\res_{x_{3}}\res_{x_2}
x_{1}^{-1}\delta\left(\frac{x-x_{3}}
{x_{1}}\right)
x_{2}^{-1}\delta\left(\frac{x_{0}-x_{3}x^{-1}x_{1}^{-1}}
{x_{2}}\right)z^{-1}\delta\left(\frac{x_1^{-1}-x_{0}}{z}\right)\cdot\nn
&&\quad \quad \cdot 
\lambda(Y_1(e^{x_{1}L(1)}(-x_{1}^{-2})^{L(0)}
Y(v, x_{3})u,
x_0)w_{(1)}\otimes w_{(2)})\nn
&&=-\res_{x_{3}}
x_{1}^{-1}\delta\left(\frac{x-x_{3}}
{x_{1}}\right)z^{-1}\delta\left(\frac{x_1^{-1}-x_{0}}{z}\right)\cdot\nn
&&\quad \quad \cdot 
\lambda(Y_1(e^{x_{1}L(1)}(-x_{1}^{-2})^{L(0)}
Y(v, x_{3})u,
x_0)w_{(1)}\otimes w_{(2)}).
\end{eqnarray}

{}From (\ref{stable-6}), (\ref{taudef}) and (\ref{cpb}), the sum of
the last two terms of (\ref{stable-4}) becomes
\begin{eqnarray}\label{stable-7}
\lefteqn{-\res_{x_{3}}
x_{1}^{-1}\delta\left(\frac{x-x_{3}}
{x_{1}}\right)z^{-1}\delta\left(\frac{x_1^{-1}-x_{0}}{z}\right)\cdot}\nn
&&\quad \quad \cdot
\lambda(Y_1(e^{x_{1}L(1)}(-x_{1}^{-2})^{L(0)}
Y(v, x_{3})u,
x_0)w_{(1)}\otimes w_{(2)})\nn
&&\quad -\res_{x_{3}}x^{-1}_1\delta\left(\frac{x-x_3}{x_1}\right)
x^{-1}_0\delta\left(\frac{z-x_1^{-1}}{-x_0}\right)
\lambda(w_{(1)}\otimes Y_2^{o}(Y(v, x_{3})u, x_{1})w_{(2)})\nn
&&=-\res_{x_{3}}
x_{1}^{-1}\delta\left(\frac{x-x_{3}}{x_{1}}\right)\left(\tau_{P(z)}\left(
x_{0}^{-1}\delta\left(\frac{x_1^{-1}-z}{x_{0}}\right)
Y_{t}(Y(v, x_{3})u, x_{1})\right)\lambda\right)(w_{(1)}\otimes w_{(2)})\nn
&&=-x_{0}^{-1}\delta\left(\frac{x_1^{-1}-z}{x_{0}}\right)\res_{x_{3}}
x_{1}^{-1}\delta\left(\frac{x-x_{3}}{x_{1}}\right)
\left(Y'_{P(z)}(Y(v, x_{3})u, x_{1})
\lambda\right)(w_{(1)}\otimes w_{(2)}).
\end{eqnarray}
Using (\ref{stable-5}), (\ref{stable-7}) and the commutator formula
for $Y'_{P(z)}$, we now see that the right-hand side of
(\ref{stable-4}) is equal to
\begin{eqnarray}\label{stable-8}
\lefteqn{x_{0}^{-1}\delta\left(\frac{x_{1}^{-1}-z}
{x_{0}}\right)\left(Y'_{P(z)}(v, x)
Y'_{P(z)}(u, x_1)\lambda\right)(w_{(1)}\otimes
w_{(2)})}\nn
&&-x_{0}^{-1}\delta\left(\frac{x_1^{-1}-z}{x_{0}}\right)\res_{x_{3}}
x_{1}^{-1}\delta\left(\frac{x-x_{3}}{x_{1}}\right)
\left(Y'_{P(z)}(Y(v, x_{3})u, x_{1})
\lambda\right)(w_{(1)}\otimes w_{(2)})\nn
&&=x_{0}^{-1}\delta\left(\frac{x_{1}^{-1}-z}
{x_{0}}\right)\left(
Y'_{P(z)}(u, x_1)Y'_{P(z)}(v, x)\lambda\right)(w_{(1)}\otimes
w_{(2)}).
\end{eqnarray}
The formulas (\ref{stable-3}), (\ref{stable-4}) and (\ref{stable-8})
together prove (\ref{stable-1}), as desired.  For the M\"obius case,
the corresponding verification for $L'_{P(z)}(-1)$, $L'_{P(z)}(0)$ and
$L'_{P(z)}(1)$ is straightforward, as usual, and we omit this
verification.  The first half of Theorem \ref{stable} holds.

For the second half of Theorem \ref{stable}, suppose that $\lambda\in
(W_{1}\otimes W_{2})^{*}$ satisfies either the $P(z)$-local grading
restriction condition or the $L(0)$-semisimple condition.  Assume
without loss of generality that $\lambda$ is doubly homogeneous.  {}From
Remark \ref{stableundercomponentops}, we see that for $v \in V$ doubly
homogeneous, $m \in \Z$ and $j=-1,0,1$, the elements
$\tau_{P(z)}(v\otimes t^{m})\lambda$ and $L'_{P(z)}(j)\lambda$ are
also doubly homogeneous.  Each such element $\mu$ lies in
$W_{\lambda}$, and so $W_{\mu} \subset W_{\lambda}$.  Thus $\mu$
satisfies the $P(z)$-local grading restriction condition (or the
$L(0)$-semisimple condition), and the second half of Theorem
\ref{stable} follows.  \epfv

\subsection{Proofs of Theorems \ref{6.1} and \ref{6.2}}

In this section, we follow \cite{tensor2}; the arguments given
there carry over to our more general situation with very little
change.  We first prove certain formulas that will be useful later.

Let $\lambda\in (W_{1}\otimes W_{2})^{*}$, $w_{(1)}\in W_{1}$ and
$w_{(2)}\in W_{2}$.  {}From (\ref{LQ'(j)}) we have
\begin{eqnarray}\label{9.1}
\lefteqn{(L'_{Q(z)}(0)\lambda)(w_{(1)} \otimes w_{(2)})}\nno\\
&&= \lambda((L(0) -
zL(1))w_{(1)} \otimes w_{(2)})- \lambda(w_{(1)} \otimes (L(0) -
zL(-1))w_{(2)}),\;\;\;\;\;
\end{eqnarray}
where (as usual) we have used the same notations $L(0), L(-1), L(1)$
to denote operators on both $W_1$ and $W_2$.  For convenience we write
$L(-1)=L'_{Q(z)}(-1)$, $L(0) = L'_{Q(z)}(0)$ and $L(1)=L'_{Q(z)}(1)$
in the rest of this section.  There will be no confusion since the
operators act on different spaces.

\begin{lemma}\label{1-y1zL(0)}
For $\lambda\in (W_{1}\otimes W_{2})^{*}$, $w_{(1)}\in W_{1}$
and $w_{(2)}\in W_{2}$,
we have
\begin{eqnarray}\label{9.2}
\lefteqn{\biggl(\biggl(1-\frac{y_1}{z}\biggr)^{L(0)} \lambda\biggr)(w_{(1)}
\otimes w_{(2)})}\nno\\
&&= \lambda\biggl(\biggl(1-\frac{y_1}{z}\biggr)^{L(0)-zL(1)}w_{(1)}
\otimes \biggl(1-\frac{y_1}{z}\biggr)^{-(L(0)-z L(-
1))}w_{(2)}\biggr).
\end{eqnarray}
\end{lemma}
\pf
{}From (\ref{9.1}),
\begin{eqnarray}
\lefteqn{\biggl(\biggl(1-\frac{y_1}{z}\biggr)^{L(0)} \lambda)(w_{(1)}
\otimes w_{(2)})}\nno\\
&&= (e^{L(0)\log(1-\frac{y_1}{z})} \lambda)(w_{(1)}
\otimes w_{(2)})\nno\\
&&= \lambda(e^{(L(0)-z L(1))\log(1-\frac{y_1}{z})}w_{(1)}
\otimes e^{-(L(0)-z L(-1))\log(1-\frac{y_1}{z})}w_{(2)})\nno\\
&&= \lambda\biggl(\biggl(1-\frac{y_1}{z}\biggr)^{L(0)-zL(1)}w_{(1)}
\otimes \biggr(1-\frac{y_1}{z}\biggr)^{-(L(0)-z L(-
1))}w_{(2)}\biggr). \hspace{1.5em}\square
\end{eqnarray}

\begin{lemma}\label{Y'Q(z)L(0)}
For $v\in V$,
\begin{equation}\label{9.4}
Y'_{Q(z)}(v,x) = \biggl(1-
\frac{y_1}{z}\biggr)^{L(0)}Y'_{Q(z)}
\biggl(\biggl(1-\frac{y_1}{z}\biggr)^{-L(0)}
v,\frac{x}{1-y_1/z}\biggr)\biggr(1-\frac{y_1}{z}\biggr)^{-
L(0)}.
\end{equation}
This formula also holds for the vertex operators associated with any
generalized $V$-module.
\end{lemma}
\pf
The identity (\ref{9.4}) will follow {}from the formula
\begin{equation}\label{9.5}
e^{yL(0)}Y'_{Q(z)}(v, x)e^{-yL(0)}=Y'_{Q(z)}(e^{yL(0)}
v, e^{y}x).
\end{equation}
To prove this, assume without loss of generality that $\mbox{wt}\
v=h\in {\mathbb Z}$, and use the $L(-1)$-derivative property
(\ref{QL-1}) and the commutator formulas (\ref{commu-q-z}) and
(\ref{qz-sl-2-qz-y-2}) to get
\begin{equation}\label{9.6}
[L(0), Y'_{Q(z)}(v, x)]=\left(x\frac{d}{dx}+h\right)Y'_{Q(z)}(v, x).
\end{equation}
Formula (\ref{9.5}) now follows {}from (an easier version of) the proof
of (\ref{710}).
\epf

\begin{lemma}\label{L(0)L(-1)formula}
Let $L(-1)$ and $L(0)$ be any operators satisfying the commutator
relation
\begin{equation}
[L(0), L(-1)]=L(-1).
\end{equation}
 Then
\begin{equation}\label{9.8}
\biggr(1-\frac{y_1}{x}\biggr)^{L(0)-xL(-1)}=e^{y_1L(-1)}
\biggr(1-\frac{y_1}{x}\biggr)^{L(0)}.
\end{equation}
\end{lemma}
\pf
We first prove that the derivative with respect to $y$ of
$$(1-y)^{L(0)-xL(-1)}
(1-y)^{-L(0)}
e^{-xyL(-1)}$$ is $0$.  Write $A=(1-y)^{L(0)-xL(-1)}$,
$B=(1-y)^{-L(0)}$,
$C=e^{-xyL(-1)}$.  Then
\begin{eqnarray}\label{9.9}
\frac{d}{dy}(ABC)
&=&-A(1-y)^{-1}(L(0)-xL(-1)))BC\nno\\
&&+A(1-y)^{-1}L(0)BC\nno\\
&&-xABL(-1)C.
\end{eqnarray}
Using (\ref{log:SL2-2}) we have
\begin{eqnarray}\label{9.10}
BL(-1)&=&(1-y)^{-L(0)}L(-1)\nno\\
&=&e^{(-\log(1-y))L(0)}L(-1)\nno\\
&=&L(-1)e^{(-\log(1-y))L(0)}e^{-\log(1-y)}\nno\\
&=&L(-1)(1-y)^{-L(0)}
(1-y)^{-1}\nno\\
&=&(1-y)^{-1}L(-1)B,
\end{eqnarray}
and substituting (\ref{9.10}) into (\ref{9.9}) gives
\begin{eqnarray*}
\frac{d}{dy}(ABC)&=&-A(1-y)^{-1}L(0)BC
+xA(1-y)^{-1}L(-1)BC
\nno\\
&&+A(1-y)^{-1}L(0)BC
-xA(1-y)^{-1}L(-1)BC\nno\\
&=&0.
\end{eqnarray*}
Thus $ABC$ is constant in $y$, and since $ABC\lbar_{y=0}=1$, we have
$ABC=1$, which is equivalent to (\ref{9.8}).\epfv

\noindent {\it Proof of Theorem \ref{6.1}} As always, the reader
should again observe the justifiability of each formal step in the
argument.

Let $\lambda$ be an element of  $(W_{1}\otimes W_{2})^{*}$
satisfying the $Q(z)$-compatibility condition, that is,
(a) the $Q(z)$-lower truncation condition---for all $v\in V$,
$Y'_{Q(z)}(v,x)\lambda = \tau_{Q(z)}(Y_t(v,x))\lambda$
involves only finitely many negative powers of $x$,
and (b)
\begin{eqnarray}\label{10.3}
\lefteqn{
\tau_{Q(z)}\left(z^{-1}\delta\left(\frac{x_{1}-x_0}{z}\right)Y_t(v,x_0)
\right)\lambda}\nno\\
&&=z^{-1}\delta\left(\frac{x_{1}-x_0}{z}\right)
Y'_{Q(z)}(v,x_0)\lambda\;\;\;
\mbox{\rm for all}\ v\in V.
\end{eqnarray}
By (\ref{5.2}) and (\ref{Y'qdef}), (\ref{10.3}) is equivalent to
\begin{eqnarray}
\lefteqn{x^{-1}_0 \delta\left(\frac{x_1-
z}{x_0}\right)\lambda(Y^o_1(v,x_1)w_{(1)} \otimes w_{(2)})}\nno\\
&&\quad- x^{-1}_0 \delta\left(\frac{z-x_1}{-x_0}\right)\lambda(w_{(1)}
\otimes
Y_2(v,x_1)w_{(2)})\nno\\
&&= z^{-1}\delta\left(\frac{x_{1}-x_0}{z}\right)\biggl(\res_{y_{1}}x^{-1}_0
\delta\left(\frac{y_1-
z}{x_0}\right)\lambda(Y^o_1(v,y_1)w_{(1)} \otimes
w_{(2)})\nno\\
&&\quad- \res_{y_1}x^{-1}_0 \delta\left(\frac{z-y_1}{-
x_0}\right)\lambda(w_{(1)} \otimes Y_2(v,y_1)w_{(2)})\biggr)
\end{eqnarray}
for all $v\in V$, $w_{(1)}\in W_{1}$ and $w_{(2)}\in W_{2}$.  It is
important to note that on the right-hand side the distributive law is
not valid since the two individual products are not defined. One
critical feature of the argument that follows is that we must rewrite
expressions to allow the application of distributivity.

By (\ref{Y'qdef}), we have
\begin{eqnarray}\label{10.5}
\lefteqn{\left(x^{-
1}_0\delta\left(\frac{x_1-x_2}{x_0}\right)
Y'_{Q(z)}(v_1,x_1)Y'_{Q(z)}(v_2,x_2)
\lambda\right)(w_{(1)} \otimes w_{(2)})}\nno\\
&&= x^{-1}_0 \delta\left(\frac{x_1-
x_2}{x_0}\right)(Y'_{Q(z)}(v_1,x_1)Y'_{Q(z)}(v_2,x_2) \lambda)(w_{(1)} \otimes
w_{(2)})\nno\\
&&= x^{-1}_0 \delta\left(\frac{x_1-x_2}{x_0}\right)\biggl(\res_{y_1}x^{-
1}_1 \delta\left(\frac{y_1-z}{x_1}\right)\cdot\nno\\
&&\hspace{6em}\cdot (Y'_{Q(z)}(v_2,x_2)
\lambda)(Y^o_1(v_1,y_1)w_{(1)} \otimes w_{(2)})\nno\\
 &&\quad -\res_{y_1} x^{-
1}_1 \delta \left(\frac{z - y_1}{-x_1}\right)(Y'_{Q(z)}(v_2,x_2)
\lambda)(w_{(1)} \otimes Y_2(v_1,y_1)w_{(2)})\biggr)\nno\\
&&= x^{-1}_0
\delta\left(\frac{x_1-x_2}{x_0}\right)\biggl(\res_{y_1}x^{-1}_1
\delta\left(\frac{y_1-
z}{x_1}\right)\res_{y_2}x^{-1}_2 \delta \left(\frac{y_2-
z}{x_2}\right)\cdot\nno\\
&&\hspace{6em}\cdot \lambda(Y^o_1(v_2,y_2)Y^o_1(v_1,y_1)w_{(1)} \otimes
w_{(2)})\nno\\
&&\quad -\res_{y_1}x^{-1}_1 \delta\left(\frac{y_1-
z}{x_1}\right)\res_{y_2}x^{-1}_2 \delta\left(\frac{z-y_2}{-
x_2}\right)\cdot \nno\\
&&\hspace{6em}\cdot \lambda(Y^o_1(v_1,y_1)w_{(1)} \otimes Y_2(v_2,y_2)w_{(2)})\nno\\
&&\quad -
\res_{y_1}x^{-1}_1 \delta\left(\frac{z-y_1}{-x_1}\right)(Y'_{Q(z)}(v_2,x_2)
\lambda)(w_{(1)} \otimes Y_2(v_1,y_1)w_{(2)})\biggr).\;\;\;\;\;\;
\end{eqnarray}
{}From the properties of the formal $\delta$-function,
we see that the right-hand side of (\ref{10.5})
is equal to
\begin{eqnarray}\label{10.6}
\lefteqn{x^{-1}_0
\delta\left(\frac{x_1-x_2}{x_0}\right)\biggl(\res_{y_1}y^{-1}_1
\delta\left(\frac{x_1+z}{y_1}\right)\res_{y_2}y^{-1}_2
\delta\left(\frac{x_2+z}{y_2}\right)\cdot} \nno\\
&&\hspace{6em}\cdot \lambda(Y^o_1(v_2,x_2+z)Y^o_1(v_1,x_1+z)w_{(
1)} \otimes w_{(2)})\nno\\
&&\quad -\res_{y_1}y^{-1}_1
\delta\left(\frac{x_1+z}{y_1}\right)\res_{y_2}x^{-1}_2
\delta\left(\frac{z-y_2}{-
x_2}\right)\cdot \nno\\
&&\hspace{6em}\cdot \lambda(Y^o_1(v_1,x_1 +z)w_{(1)} \otimes
Y_2(v_2,y_2)w_{(2)})\nno\\
&&\quad - \res_{y_1}x^{-1}_1 \delta\left(\frac{z-y_1}{-
x_1}\right)(Y'_{Q(z)}(v_2,x_2) \lambda)(w_{(1)} \otimes
Y_2(v_1,y_1)w_{(2)})\biggr)\nno\\
&&= x^{-1}_0 \delta \left(\frac{x_1-
x_2}{x_0}\right)\biggl(\lambda(Y^o_1(v_2,x_2+z)
Y^o_1(v_1,x_1+z)w_{(1)}
\otimes w_{(2)})\nno\\
&&\quad - \res_{y_2}x^{-1}_2 \delta\left(\frac{z-y_2}{-
x_2}\right)\lambda(Y^o_1(v_1,x_1+z)w_{(1)} \otimes
Y_2(v_2,y_2)w_{(2)})\nno\\
&&\quad - \res_{y_1}x^{-1}_1 \delta\left(\frac{z-y_1}{-
x_1}\right)(Y'_{Q(z)}(v_2,x_2) \lambda)(w_{(1)} \otimes
Y_2(v_1,y_1)w_{(2)})\biggr).\;\;\;\;\;\;\;\;\;
\end{eqnarray}
{}From the $L(-1)$-derivative property for $Y_{1}$, the
$L(-1)$-derivative property (\ref{yo-l-1}) for $Y_{1}^{o}$ and the
commutator formulas for $L(-1)$, $Y_{1}(\cdot, x)$ and for $L(1)$,
$Y^{o}_{1}(\cdot, x)$ (recall Lemma \ref{sl2opposite}), we obtain
\begin{eqnarray}\label{10.7}
Y_{1}(v, x+z)&=&Y_{1}(e^{zL(-1)}v, x)\nno\\
&=&e^{zL(-1)}Y_{1}(v, x)e^{-zL(-1)}\nno\\
&=&\sum_{n\ge
0}\frac{z^{n}}{n!}\frac{d^{n}}{dx^{n}}Y_{1}(v, x)
\end{eqnarray}
and
\begin{eqnarray}\label{10.8}
Y^{o}_{1}(v, x+z)&=&Y^{o}_{1}(e^{zL(-1)}v, x)\nno\\
&=&e^{-zL(1)}Y^{o}_{1}(v, x)e^{zL(1)}\nno\\
&=&\sum_{n\ge
0}\frac{z^{n}}{n!}\frac{d^{n}}{dx^{n}}Y^{o}_{1}(v, x).
\end{eqnarray}
(Note that all these expressions are in fact defined.)
Using (\ref{10.8}), we see that the right-hand side of (\ref{10.6})
can be written as
\begin{eqnarray}\label{10.9}
\lefteqn{x^{-1}_0 \delta \left(\frac{x_1-
x_2}{x_0}\right)\biggl(\lambda(Y^o_1(e^{zL(-1)}v_2,x_2)Y^o_1(e^{zL(-1)}v_1,x_1)w_{(1)}
\otimes w_{(2)})}\nno\\
&&\quad - \res_{y_2}x^{-1}_2 \delta\left(\frac{z-y_2}{-
x_2}\right)\lambda(Y^o_1(e^{zL(-1)}v_1,x_1)w_{(1)} \otimes
Y_2(v_2,y_2)w_{(2)})\nno\\
&&\quad - \res_{y_1}x^{-1}_1 \delta\left(\frac{z-y_1}{-
x_1}\right)(Y'_{Q(z)}(v_2,x_2) \lambda)(w_{(1)} \otimes
Y_2(v_1,y_1)w_{(2)})\biggr)\nno\\
&&=x^{-1}_0 \delta \left(\frac{x_1-
x_2}{x_0}\right)\biggl(\lambda(e^{-zL(1)}Y^o_1(v_2,x_2)Y^o_1(v_1,x_1)e^{zL(1)}
w_{(1)}
\otimes w_{(2)})\nno\\
&&\quad - \res_{y_2}x^{-1}_2 \delta\left(\frac{z-y_2}{-
x_2}\right)\lambda(e^{-zL(1)}Y^o_1(v_1,x_1)e^{zL(1)}w_{(1)} \otimes
Y_2(v_2,y_2)w_{(2)})\nno\\
&&\quad - \res_{y_1}x^{-1}_1 \delta\left(\frac{z-y_1}{-
x_1}\right)(Y'_{Q(z)}(v_2,x_2) \lambda)(w_{(1)} \otimes
Y_2(v_1,y_1)w_{(2)})\biggr).
\end{eqnarray}
Note that it is easier to verify the well-definedness of the terms on
the right-hand side of (\ref{10.6}) than that of the terms in (\ref{10.9}),
though
(\ref{10.9}) is sometimes easier to use if it is known that every term is
well defined.  Below we shall write expressions like those on the right-hand
side of (\ref{10.6}) in whichever way suits our needs.
The distributive law applies to the right-hand side of (\ref{10.6})
(or (\ref{10.9})) since
all three of the following expressions are defined:
$$x^{-1}_0
\delta\left(\frac{x_1-x_2}{x_0}\right)\lambda(Y^o_1(v_2,x_2+z)
Y^o_1(v_1,x_1+z)w_{(1)} \otimes w_{(2)}),$$ $$x^{-1}_0
\delta\left(\frac{x_1-x_2}{x_0}\right)\res_{y_2}x^{-1}_2
\delta\left(\frac{z-y_2}{-x_2}\right)\lambda(Y^o_1(v_1,x_1+z)w_{(1)} \otimes
Y_2(v_2,y_2)w_{(2)}),$$
$$x^{-1}_0
\delta\left(\frac{x_1-x_2}{x_0}\right)\res_{y_1}x^{-1}_1 \delta\left(\frac{z-
y_1}{-x_1}\right)(Y'_{Q(z)}(v_2,x_2) \lambda)(w_{(1)} \otimes
Y_2(v_1,y_1)w_{(2)}).$$

Now we examine the last expression in (\ref{10.9}).  Rewriting the
formal $\delta$-functions
$x_{0}^{-1}\delta\left(\frac{x_{1}-x_{2}}{x_{0}}\right)$ and
$x_{1}^{-1}\delta\left(\frac{z-y_1}{-x_1}\right)$, and using Lemma
\ref{Y'Q(z)L(0)} and (\ref{deltafunctionsubstitutionformula}), we
have:
\begin{eqnarray}\label{10.10}
\lefteqn{x^{-1}_0\delta\left(\frac{x_1-x_2}{x_0}\right)\res_{y_{1}}x^{-1}_1
\delta\left(\frac{z-y_1}{-x_1}\right)\cdot}\nno\\
&&\hspace{4em}\cdot (Y'_{Q(z)}(v_2, x_2)\lambda)(w_{(1)}\otimes
 Y_2(v_1,
y_1)w_{(2)})\nno\\
&&=\biggl(\frac{x_1}{z}\biggr)^{-1}\biggl(\frac{x_0}{x_1/z}\biggr)^{-
1}\delta\left(\frac{z+(\frac{x_2}{-x_1/z})}{\frac{x_0}{x_1/z}}\right)
\res_{y_1}x_1^{-1}\delta\left(\frac{1-y_1/z}{-x_1/z}\right)\cdot \nno\\
&&\hspace{4em}\cdot\biggl(\biggl(1-\frac{y_1}{z}\biggr)^{L(0)}Y'_{Q(z)}
\biggl(\biggl(1-\frac{y_1}{z}\biggr)^{-
L(0)}v_2, \frac{x_2}{1-y_{1}/z}\biggr)\cdot \nno\\
&&\hspace{5em}\cdot \biggl(1-\frac{y_1}{z}\biggr)^{-L(0)}
\lambda\biggr)(w_{(1)}\otimes Y_2(v_1, y_1)w_{(2)})\nno\\
&&=\biggl(\frac{x_1}{z}\biggr)^{-1}\biggl(\frac{x_0}{x_1/z}\biggr)^{-
1}\delta\left(\frac{z+(\frac{x_2}{-x_1/z})}{\frac{x_0}{x_1/z}}\right)
\res_{y_1}x_1^{-1}\delta\left(\frac{1-y_1/z}{-x_1/z}\right)\cdot \nno\\
&&\hspace{4em}\cdot\biggl(\biggl(1-\frac{y_1}{z}\biggr)^{L(0)}Y'_{Q(z)}
\biggl(\biggl(1-\frac{y_1}{z}\biggr)^{-
L(0)}v_2, \frac{x_2}{-x_{1}/z}\biggr)\cdot \nno\\
&&\hspace{5em}\cdot \biggl(1-\frac{y_1}{z}\biggr)^{-L(0)}
\lambda\biggr)(w_{(1)}\otimes Y_2(v_1, y_1)w_{(2)}).
\end{eqnarray}
By Lemma \ref{1-y1zL(0)} and (\ref{10.3}), the right-hand side of
 (\ref{10.10}) is equal to
\begin{eqnarray}\label{10.11}
\lefteqn{\biggl(\frac{x_1}{z}\biggr)^{-1}\biggl(\frac{x_0}{x_1/z}\biggr)^{-
1}\delta\left(\frac{z+(\frac{x_2}{-x_1/z})}{\frac{x_0}{x_1/z}}\right)
\res_{y_1}x_1^{-1}\delta\left(\frac{1-y_1/z}{-x_1/z}\right)\cdot} \nno\\
&&\hspace{1em}\cdot \biggl(Y'_{Q(z)}\biggl(\biggl(1-\frac{y_1}{z}\biggr)^{-L(0)}
v_2, \frac{x_2}{-x_1/z}\biggr)\biggl(1-
\frac{y_1}{z}\biggr)^{-L(0)} \lambda\biggr)\cdot \nno\\
&&\hspace{3em}\cdot \biggl(\biggl(1-\frac{y_1}{z}\biggr)^{L(0)-
zL(1)}w_{(1)}\otimes\biggl(1-\frac{y_1}{z}\biggr)^{-
(L(0)-zL(-1))}Y_2(v_1, y_1)w_{(2)}\biggr)\nno\\
&&=\biggl(\frac{x_1}{z}\biggr)^{-1}
\res_{y_1}x_1^{-1}\delta\left(\frac{1-y_1/z}{-x_1/z}\right)\cdot \nno\\
&&\hspace{2em} \cdot\Biggl(\tau_{Q(z)}\Biggl(\biggl(\frac{x_0}{x_1/z}
\biggr)^{-
1}\delta\left(\frac{z+(\frac{x_2}{-x_1/z})}{\frac{x_0}{x_1/z}}\right)
\cdot \nno\\
&&\hspace{3em}\cdot Y_{t}\biggl(\biggl(1-\frac{y_1}{z}\biggr)^{-L(0)}v_2,
\frac{x_2}{-x_1/z}\biggr)\Biggr)\biggl(1-
\frac{y_1}{z}\biggr)^{-L(0)} \lambda\Biggr)\cdot \nno\\
&&\hspace{4em}\cdot \biggl(\biggl(1-\frac{y_1}{z}\biggr)^{L(0)-
zL(1)}w_{(1)}\otimes\biggl(1-\frac{y_1}{z}\biggr)^{-
(L(0)-zL(-1))}Y_2(v_1, y_1)w_{(2)}\biggr).
\end{eqnarray}
By (\ref{5.2}), the right-hand side of (\ref{10.11}) becomes
\begin{eqnarray}\label{10.12}
\lefteqn{\biggl(\frac{x_1}{z}\biggr)^{-1}
\res_{y_1}x_1^{-1}\delta\left(\frac{1-y_1/z}{-x_1/z}\right)
\biggl(\biggl(\frac{x_2}{-x_1/z}\biggr)^{-
1}\delta\left(\frac{\frac{x_0}{x_1/z}-z}{\frac{x_2}{-x_1/z}}\right)\cdot}
\nno\\
&&\hspace{1em}\cdot \biggl(\biggl(1-\frac{y_1}{z}\biggr)^{-L(0)}\lambda\biggr)
\biggl(Y_1^o\biggl(\biggl(1-\frac{y_1}{z}\biggr)^{-L(0)}v_2,\frac{x_0}
{x_{1}/z}\biggr)\cdot \nno\\
&&\hspace{2em}\cdot \biggl(1-\frac{y_1}{z}\biggr)^{L(0)-zL(1)}
w_{(1)}\otimes \biggl(1-\frac{y_1}{z}\biggr)^{-(L(0)-zL(1))}Y_2(v_1,
y_1)w_{(2)}\biggr)\nno\\
&&\quad -\biggl(\frac{x_2}{-x_1/z}\biggr)^{-
1}\delta\left(\frac{z-\frac{x_0}{x_1/z}}{-\frac{x_2}{-x_1/z}}\right)
\biggl(\biggl(1-\frac{y_1}{z}\biggr)^{-L(0)}\lambda\biggr)
\biggl(\biggl(1-\frac{y_1}{z}\biggr)^{L(0)-zL(1)} \cdot\nno\\
&&\hspace{2em}\cdot w_{(1)} \otimes
 Y_2\biggl(\biggl(1-\frac{y_1}{z}\biggr)^{-L(0)}v_2,
\frac{x_0}{x_1/z}\biggr)\cdot\nno\\
&&\hspace{3em}\cdot
\biggl(1-\frac{y_1}{z}\biggr)^{-(L(0)-zL(-1))}Y_2(v_1,
y_1)w_{(2)}\biggr)\biggr).\;\;\;\;\;\;\;
\end{eqnarray}
Using Lemma \ref{1-y1zL(0)} again but with $1-\frac{y_{1}}{z}$
replaced by $(1-\frac{y_{1}}{z})^{-1}$, rewriting formal
$\delta$-functions and then using the distributive law, we see that
(\ref{10.12}) is equal to
\begin{eqnarray}\label{10.13}
\lefteqn{\res_{y_1}x_1^{-1}\delta\left(\frac{1-y_1/z}{-x_1/z}\right)
\biggl(-x_2^{-1}\delta\left(\frac{x_0-x_1}{-x_2}\right)\cdot} \nno\\
&&\hspace{4em}\cdot \lambda\biggl(\biggl(1-\frac{y_1}{z}\biggr)^{-(L(0)-zL(1))}
Y_1^o\biggl(\biggl(1-\frac{y_1}{z}\biggr)^{-L(0)}v_2,\frac{x_0}{-(1-
y_1/z)}\biggr)\cdot \nno\\
&&\hspace{5em}\cdot \biggl(1-\frac{y_1}{z}\biggr)^{L(0)-zL(1)}
w_{(1)}\otimes Y_2(v_1, y_1)w_{(2)}\biggr)\nno\\
&&\quad +x_2^{-1}\delta\left(\frac{x_1-x_0}{x_2}\right)
\lambda\biggl(w_{(1)}\otimes\biggl(1-\frac{y_1}{z}\biggr)^{L(0)-zL(-1)}\cdot \nno\\
&&\hspace{4em}\cdot Y_2\biggl(\biggl(1-\frac{y_1}{z}\biggr)^{-L(0)}v_2,
\frac{x_0}{-(1-
y_1/z)}\biggr)\cdot \nno\\
&&\hspace{5em}\cdot \biggl(1-\frac{y_1}{z}\biggr)^{-
(L(0)-zL(-1))}Y_2(v_1, y_1)w_{(2)})\biggr)\nno\\
&&=-\res_{y_1}x_1^{-1}\delta\left(\frac{1-y_1/z}{-x_1/z}\right)
x_2^{-1}\delta\left(\frac{x_0-x_1}{-x_2}\right)\cdot \nno\\
&&\hspace{4em}\cdot \lambda\biggl(\biggl(1-\frac{y_1}{z}\biggr)^{-(L(0)-zL(1))}
Y_1^o\biggl(\biggl(1-\frac{y_1}{z}\biggr)^{-L(0)}v_2,\frac{x_0}{-(1-
y_1/z)}\biggr)\cdot \nno\\
&&\hspace{5em}\cdot \biggl(1-\frac{y_1}{z}\biggr)^{L(0)-zL(1)}
w_{(1)}\otimes Y_2(v_1, y_1)w_{(2)}\biggr)\nno\\
&&\quad +\res_{y_1}x_1^{-1}\delta\left(\frac{1-y_1/z}{-x_1/z}\right)
x_2^{-1}\delta\left(\frac{x_1-x_0}{x_2}\right)
\lambda\biggl(w_{(1)}\otimes\nno\\
&&\hspace{4em}\otimes\biggl(1-\frac{y_1}{z}\biggr)^{L(0)-zL(-1)}
Y_2\biggl(\biggl(1-\frac{y_1}{z}\biggr)^{-L(0)}v_2,
\frac{x_0}{-(1-
y_1/z)}\biggr)\cdot \nno\\
&&\hspace{5em}\cdot \biggl(1-\frac{y_1}{z}\biggr)^{-
(L(0)-zL(-1))}Y_2(v_1, y_1)w_{(2)}\biggr).
\end{eqnarray}

But by Lemmas \ref{L(0)L(-1)formula} and \ref{Y'Q(z)L(0)},
\begin{eqnarray}\label{10.14}
\lefteqn{\biggl(1-\frac{y_1}{z}\biggr)^{L(0)-zL(-1)} Y_2
\biggl(\biggl(1-\frac{y_1}{z}\biggr)^{-L(0)}v_2, \frac{x_0}{-(1-
y_1/z)}\biggr)\biggl(1-\frac{y_1}{z}\biggr)^{-(L(0)-zL(-
1))}}\nno\\
&&=e^{y_{1}L(-1)}\biggl(1-\frac{y_1}{z}\biggr)^{L(0)}Y_2
\biggl(\biggl(1-\frac{y_1}{z}\biggr)^{-L(0)}v_2, \frac{x_0}{-(1-
y_1/z)}\biggr)\cdot\hspace{4em}\nno\\
&&\hspace{10em}\cdot\biggl(1-\frac{y_1}{z}\biggr)^{-L(0)}e^{-y_{1}L(-1)}\nno\\
&&=e^{y_{1}L(-1)}Y_2(v_2, -x_0) e^{-y_{1}L(-1)}\nno\\
&&=Y_2(v_2, -x_0+ y_1)\nno\\
&&=Y_2(v_2, -x_0-(z-y_1)+z)\nno\\
&&=Y_2(e^{zL(-1)}v_2, -x_0-
(z-y_1)).\hspace{17em}
\end{eqnarray}
We similarly have, using Lemmas \ref{L(0)L(-1)formula} and
\ref{Y'Q(z)L(0)} for $Y'_1(v_2, x)$ and then using (\ref{y'}) and
Theorem \ref{set:W'},
\begin{eqnarray}\label{10.15}
\lefteqn{\biggl(1-\frac{y_1}{z}\biggr)^{-(L(0)-zL(1))}
Y_1^o\biggl(\biggl(1-\frac{y_1}{z}\biggr)^{-L(0)}v_2,\frac{x_0}{-(1-
y_1/z)}\biggr)\biggl(1-\frac{y_1}{z}\biggr)^{L(0)-zL(1)}}\nno\\
&&=e^{-y_{1}L(1)}Y_{1}^{o}(v_{2}, -x_{0})e^{y_{1}L(1)}\nno\\
&&=Y_{1}^{o}(v_{2}, -x_{0}+y_{1})\nno\\
&&=Y_{1}^{o}(v_{2}, -x_{0}-(z-y_{1})+z)\nno\\
&& =Y_1^o(e^{zL(-1)}v_2, -x_0-
(z-y_1)).\hspace{17em}
\end{eqnarray}
Substituting (\ref{10.14}) and (\ref{10.15}) into the
right-hand side of (\ref{10.13}) and then
combining  with  (\ref{10.10})--(\ref{10.13}),
we obtain
\begin{eqnarray}\label{10.16}
\lefteqn{x^{-1}_0\delta\left(\frac{x_1-x_2}{x_0}\right)\res_{y_1}
x^{-1}_1\delta\left(\frac{z-y_1}{-x_1}\right)(Y'_{Q(z)}(v_2, x_2) \lambda)
(w_{(1)}\otimes Y_2(v_1, y_1)w_{(2)})}\nno\\
&&=-\res_{y_1} x^{-1}_1\delta\left(\frac{z-y_1}{-x_1}\right)
x^{-1}_2\delta\left(\frac{x_0-x_1}{-x_2}\right)\cdot \nno\\
&&\hspace{4em}\cdot \lambda(Y_1^o\biggl(v_2,
-x_0-(z-y_1)+z\biggr)
w_{(1)}\otimes Y_2(v_1, y_1)w_{(2)})\nno\\
&&\quad +\res_{y_1} x^{-1}_1\delta\left(\frac{z-y_1}{-x_1}\right)
x^{-1}_2\delta\left(\frac{x_1-x_0}{x_2}\right)\cdot \nno\\
&&\hspace{4em}\cdot \lambda(w_{(1)}\otimes Y_2(v_2, -x_0+y_1)
Y_2(v_1, y_1)w_{(2)}).\hspace{8em}
\end{eqnarray}
(We choose the form of the expression {}from (\ref{10.15})
in anticipation of the
next step.)

By (\ref{deltafunctionsubstitutionformula}) and (\ref{10.7}), the
right-hand side of (\ref{10.16}) is equal to
\begin{eqnarray}\label{10.17}
\lefteqn{-x^{-1}_2\delta\left(\frac{-x_0+x_1}{x_2}\right)\res_{y_1}
x^{-1}_1\delta\left(\frac{z-y_1}{-x_1}\right)\cdot} \nno\\
&&\hspace{4em}\cdot \lambda(Y_1^o(v_2, -x_0+x_1+z)w_{(1)}\otimes Y_2(v_1,
y_1)w_{(2)})\nno\\
&&\quad +x^{-1}_2\delta\left(\frac{x_1-x_0}{x_2}\right)\res_{y_1}
x^{-1}_1\delta\left(\frac{z-y_1}{-x_1}\right)\cdot \nno\\
&&\hspace{4em}\cdot \lambda(w_{(1)}\otimes Y_2(v_2, -x_0+y_1)
Y_2(v_1, y_1)w_{(2)})\nno\\
&&=-x^{-1}_2\delta\left(\frac{-x_0+x_1}{x_2}\right)\res_{y_1}
x^{-1}_1\delta\left(\frac{z-y_1}{-x_1}\right)\cdot \nno\\
&&\hspace{4em}\cdot \lambda(Y_1^o(e^{zL(-1)}v_2, x_{2})w_{(1)}\otimes Y_2(v_1,
y_1)w_{(2)})\nno\\
&&\quad +x^{-1}_2\delta\left(\frac{x_1-x_0}{x_2}\right)\res_{y_1}
x^{-1}_1\delta\left(\frac{z-y_1}{-x_1}\right)\cdot \nno\\
&&\hspace{4em}\cdot \lambda(w_{(1)}\otimes Y_2(v_2, -x_0+y_1)
Y_2(v_1, y_1)w_{(2)}).
\end{eqnarray}
Since
$$\res_{y_{2}}y_{2}^{-1}\delta\left(\frac{-x_{0}+y_{1}}{y_{2}}\right)=1,$$
the right-hand side of (\ref{10.17}) can be written as
\begin{eqnarray}\label{10.18}
&&-x^{-1}_2\delta\left(\frac{-x_0+x_1}{x_2}\right)\res_{y_1}
x^{-1}_1\delta\left(\frac{z-y_1}{-x_1}\right)\cdot \nno\\
&&\hspace{4em}\cdot \lambda(Y_1^o(e^{zL(-1)}v_2, x_{2})w_{(1)}\otimes Y_2(v_1,
y_1)w_{(2)})\nno\\
&&\quad +x^{-1}_2\delta\left(\frac{x_1-x_0}{x_2}\right)\res_{y_1}
x^{-1}_1\delta\left(\frac{z-y_1}{-x_1}\right)\res_{y_2}
y^{-1}_2\delta\left(\frac{-x_0+y_1}{y_2}\right)\cdot \nno\\
&&\hspace{4em}\cdot \lambda(w_{(1)}\otimes Y_2(v_2, -x_0+y_1)
Y_2(v_1, y_1)w_{(2)}).
\end{eqnarray}

By (\ref{2termdeltarelation}) and 
(\ref{deltafunctionsubstitutionformula}), (\ref{10.18}) becomes
\begin{eqnarray}\label{10.19}
\lefteqn{x^{-1}_0\delta\left(\frac{x_2-x_1}{-x_0}\right)
\res_{y_1} x^{-1}_1\delta\left(\frac{z-y_1}{-x_1}\right)\cdot} \nno\\
&&\hspace{4em}\cdot \lambda(Y_1^o(e^{zL(-1)}v_2, x_2)w_{(1)}\otimes Y_2(v_1,
y_1)w_{(2)})\nno\\
&&\quad +x^{-1}_2\delta\left(\frac{x_1-x_0}{x_2}\right)\res_{y_1}
x^{-1}_1\delta\left(\frac{z-y_1}{-x_1}\right)\res_{y_2}
y^{-1}_2\delta\left(\frac{-x_0+y_1}{y_2}\right)\cdot \nno\\
&&\hspace{4em}\cdot \lambda(w_{(1)}\otimes Y_2(v_2, y_2)Y_2(v_1, y_1)w_{(2)})\nno\\
&&=x^{-1}_0\delta\left(\frac{x_2-x_1}{-x_0}\right)\res_{y_1}
x^{-1}_1\delta\left(\frac{z-y_1}{-x_1}\right)\cdot \nno\\
&&\hspace{4em}\cdot \lambda(Y_1^o(e^{zL(-1)}v_2, x_2)w_{(1)}\otimes Y_2(v_1,
y_1)w_{(2)})\nno\\
&&\quad +x^{-1}_2\delta\left(\frac{x_1-x_0}{x_2}\right)\res_{y_1}
\res_{y_2} x^{-1}_1\delta\left(\frac{z-y_1}{-x_1}\right)
y^{-1}_2\delta\left(\frac{-x_0+y_1}{y_2}\right)\cdot \nno\\
&&\hspace{4em}\cdot \lambda(w_{(1)}\otimes Y_2(v_2, y_2)Y_2(v_1, y_1)w_{(2)})\nno\\
&&=x^{-1}_0\delta\left(\frac{x_2-x_1}{-x_0}\right)\res_{y_1}
x^{-1}_1\delta\left(\frac{z-y_1}{-x_1}\right)\cdot \nno\\
&&\hspace{4em}\cdot \lambda(Y_1^o(e^{zL(-1)}v_2, x_2)w_{(1)}\otimes Y_2(v_1,
y_1)w_{(2)})\nno\\
&&\quad -x^{-1}_2\delta\left(\frac{x_1-x_0}{x_2}\right)\res_{y_1} \res_{y_2}
\cdot\nno\\
&&\hspace{3em}\cdot
(x_2+x_0)^{-1}\delta\left(\frac{z-y_1}{-x_2-x_0}\right)
x^{-1}_0\delta\left(\frac{y_2-y_1}{-x_0}\right)\cdot \nno\\
&&\hspace{4em}\cdot \lambda(w_{(1)}\otimes Y_2(v_2, y_2)Y_2(v_1, y_1)w_{(2)})\nno\\
&&=x^{-1}_0\delta\left(\frac{x_2-x_1}{-x_0}\right)\res_{y_1}
x^{-1}_1\delta\left(\frac{z-y_1}{-x_1}\right)\cdot \nno\\
&&\hspace{4em}\cdot \lambda(Y_1^o(e^{zL(-1)}v_2, x_2)w_{(1)}\otimes Y_2(v_1,
y_1)w_{(2)})\nno\\
&&\quad -x^{-1}_2\delta\left(\frac{x_1-x_0}{x_2}\right)\res_{y_1} \res_{y_2}
x^{-1}_2\delta\left(\frac{z-y_2}{-x_2}\right)x^{-1}_0\delta\left(\frac{y_2-
y_1}{-x_0}\right)\cdot \nno\\
&&\hspace{4em}\cdot \lambda(w_{(1)}\otimes Y_2(v_2, y_2)Y_2(v_1, y_1)w_{(2)}).
\end{eqnarray}
Substituting  (\ref{10.16})--(\ref{10.19}) into (\ref{10.9}) we obtain
\begin{eqnarray}\label{10.20}
\lefteqn{\left(x^{-1}_0\delta\left(\frac{x_1-x_2}{x_0}\right)Y'_{Q(z)}(v_1, x_1)
Y'_{Q(z)}(v_2,
x_2) \lambda\right)(w_{(1)}\otimes w_{(2)})}\nno\\
&&=x^{-1}_0\delta\left(\frac{x_1-x_2}{x_0}\right)\lambda(Y_1^o(e^{zL(-
1)}v_2, x_2)Y_1^o(e^{zL(-
1)}v_1, x_1)w_{(1)}\otimes w_{(2)})\nno\\
&&\quad -x^{-1}_0\delta\left(\frac{x_1-x_2}{x_0}\right)
\res_{y_2} x^{-1}_2\delta\left(\frac{z-y_2}{-x_2}\right)\cdot\nno\\
&&\hspace{6em}\cdot\lambda(Y_1^o(e^{zL(-1)}v_1, x_1)w_{(1)}\otimes Y_2(v_2,
x_2)w_{(2)})\nno\\
&&\quad -x^{-1}_0\delta\left(\frac{x_2-x_1}{-x_0}\right)\res_{y_1}
x^{-1}_1\delta\left(\frac{z-y_1}{-x_1}\right)\cdot\nno\\
&&\hspace{6em}\cdot\lambda(Y_1^o(e^{zL(-1)}v_2, x_2)w_{(1)}\otimes Y_2(v_1,
y_1)w_{(2)})\nno\\
&&\quad +x^{-1}_2\delta\left(\frac{x_1-x_0}{x_2}\right)\res_{y_1}
\res_{y_2} x^{-1}_2\delta\left(\frac{z-y_2}{-x_2}\right)
x^{-1}_0\delta\left(\frac{y_2-y_1}{-x_0}\right)\cdot\nno\\
&&\hspace{6em}\cdot\lambda(w_{(1)}\otimes Y_2(v_2, y_2)Y_2(v_1, y_1)w_{(2)}).
\end{eqnarray}

Now consider the result of the calculation {}from (\ref{10.16})
to (\ref{10.19})
except for the last two steps in (\ref{10.19}). Reversing
the subscripts 1 and 2 of the symbols $v$, $x$ and $y$ and replacing 
$x_0$ by $-x_0$ in this result and then using (\ref{2termdeltarelation}), 
we have
\begin{eqnarray}\label{10.21}
\lefteqn{x^{-1}_0\delta\left(\frac{x_2-x_1}{-x_0}\right)\res_{y_2}
x^{-1}_2\delta\left(\frac{z-y_2}{-x_2}\right)\cdot }\nno\\
&&\hspace{6em}\cdot (Y'_{Q(z)}(v_1, x_1) \lambda)
(w_{(1)}\otimes Y_2(v_2, y_2)w_{(2)})\nno\\
&& =x^{-1}_0\delta\left(\frac{x_1-x_2}{x_0}\right)\res_{y_2}
x^{-1}_2\delta\left(\frac{z-y_2}{-x_2}\right)\cdot \nno\\
&& \hspace{6em}\cdot \lambda(Y_1^o(e^{zL(-1)}v_1, x_1)w_{(1)}\otimes Y_2(v_2,
x_2)w_{(2)})\nno\\
&&\quad -x^{-1}_2\delta\left(\frac{x_1-x_0}{x_2}\right)\res_{y_1}\res_{y_2}
x^{-1}_2\delta\left(\frac{z-y_2}{-x_2}\right)
x^{-1}_0\delta\left(\frac{y_1-y_2}{x_0}\right)\cdot \nno\\
&&\hspace{6em}\cdot \lambda(w_{(1)}\otimes Y_2(v_1, y_1)Y_2(v_2, y_2)w_{(2)}).
\end{eqnarray}
{}From (\ref{10.5})--(\ref{10.9}),
again reversing the subscripts 1 and 2 of the symbols
$v$, $x$ and $y$
and replacing $x_0$ by $-x_0$, and (\ref{10.21}), we have
 \begin{eqnarray}\label{10.22}
\lefteqn{\left(-x^{-1}_0\delta\left(\frac{x_2-x_1}{-x_0}\right)Y'_{Q(z)}(v_2,
x_2)Y'_{Q(z)}(v_1, x_1)\cdot \lambda\right) (w_{(1)}\otimes w_{(2)})}\nno\\
&&=-x^{-1}_0\delta\left(\frac{x_2-x_1}{-
x_0}\right)\lambda(Y_1^o(e^{zL(-1)}v_1, x_1)Y_1^o(e^{zL(-
1)}v_2, x_2)w_{(1)}\otimes w_{(2)})\nno\\
&&\quad +x^{-1}_0\delta\left(\frac{x_2-x_1}{-x_0}\right)\res_{y_1}
x^{-1}_1\delta\left(\frac{z-y_1}{-x_1}\right)\cdot \nno\\
&&\hspace{6em}\cdot \lambda(Y_1^o(e^{zL(-1)}v_2,
x_2)w_{(1)}\otimes Y_2(v_1, y_1)w_{(2)})\nno\\
&&\quad +x^{-1}_0\delta\left(\frac{x_1-x_2}{x_0}\right)\res_{y_2}
x^{-1}_2\delta\left(\frac{z-y_2}{-x_2}\right)\cdot \nno\\
&&\hspace{6em}\cdot \lambda(Y_1^o(e^{zL(-1)}v_1, x_1)w_{(1)}\otimes Y_2(v_2,
x_2)w_{(2)})\nno\\
&&-x^{-1}_2\delta\left(\frac{x_1-x_0}{x_2}\right)\res_{y_1}\res_{y_2}
x^{-1}_2\delta\left(\frac{z-y_2}{-x_2}\right) x^{-1}_0\delta\left(\frac{y_1-
y_2}{x_0}\right)\cdot \nno\\
&&\hspace{6em}\cdot \lambda(w_{(1)}\otimes Y_2(v_1, y_1)Y_2(v_2, y_2)w_{(2)}).
\end{eqnarray}

The formulas (\ref{10.20}) and (\ref{10.22}) give:
\begin{eqnarray}\label{10.23}
\lefteqn{\biggl(\biggl(x^{-1}_0\delta\left(\frac{x_1-x_2}{x_0}\right)
Y'_{Q(z)}(v_1, x_1)
Y'_{Q(z)}(v_2,
x_2)}\nno\\
&&\quad  -x^{-1}_0\delta\left(\frac{x_2-x_1}{-x_0}\right)
Y'_{Q(z)}(v_2, x_2)Y'_{Q(z)}(v_1,
x_1)\biggr) \lambda\biggr)(w_{(1)}\otimes w_{(2)})\nno\\
&&=\lambda\biggl(\biggl(x^{-1}_0\delta\left(\frac{x_1-x_2}{x_0}\right)
Y_1^o(e^{zL(-
1)}v_2, x_2)Y_1^o(e^{zL(-1)}v_1, x_1)\nno\\
&&\hspace{2em}-x^{-1}_0\delta\left(\frac{x_2-x_1}{-x_0}\right)
Y_1^o(e^{zL(-1)}v_1,
x_1)Y_1^o(e^{zL(-1)}v_2, x_2)\biggr)w_{(1)}\otimes w_{(2)}\biggr)\nno\\
&&\quad -x^{-1}_2\delta\left(\frac{x_1-x_0}{x_2}\right)\res_{y_1}\res_{y_2}
x^{-1}_2\delta\left(\frac{z-y_2}{-x_2}\right)\cdot \nno\\
&&\hspace{4em}\cdot \lambda\biggl(w_{(1)}\otimes \biggl(x^{-1}_0
\delta\left(\frac{y_1-y_2}{x_0}\right)
Y_2(v_1, y_1)Y_2(v_2, y_2)\nno\\
&&\hspace{2em}-x^{-1}_0\delta\left(\frac{y_2-y_1}{-x_0}\right)
Y_2(v_2, y_2)Y_2(v_1,
y_1)\biggr)w_{(2)}\biggr).
\end{eqnarray}
{}From the Jacobi identities for $Y_{1}^{o}$ and $Y_{2}$ and 
(\ref{log:p1}), the right-hand
side of (\ref{10.23}) is equal to
\begin{eqnarray}\label{10.24}
\lefteqn{x^{-1}_2\delta\left(\frac{x_1-x_0}{x_2}\right)\lambda(Y_1^o(Y(e^{zL(-
1)}v_1, x_0)e^{zL(-1)}v_2, x_2)w_{(1)}\otimes w_{(2)})}\nno\\
&&\quad -x^{-1}_2\delta\left(\frac{x_1-x_0}{x_2}\right)\res_{y_1}\res_{y_2}
x^{-1}_2\delta\left(\frac{z-y_2}{-x_2}\right) y^{-1}_2\delta\left(\frac{y_1-
x_0}{y_2}\right)\cdot \nno\\
&&\hspace{6em}\cdot \lambda(w_{(1)}\otimes Y_2(Y(v_1, x_0)v_2, y_2)w_{(2)})
\nno\\
&&=x^{-1}_2\delta\left(\frac{x_1-x_0}{x_2}\right)\lambda(Y_1^o(e^{zL(-
1)}Y(v_1, x_0)v_2, x_2)w_{(1)}\otimes w_{(2)})\nno\\
&&\quad -x^{-1}_2\delta\left(\frac{x_1-x_0}{x_2}\right)\res_{y_1}\res_{y_2}
x^{-1}_2\delta\left(\frac{z-y_2}{-x_2}\right)
y^{-1}_2\delta\left(\frac{y_1-
x_0}{y_2}\right)\cdot \nno\\
&&\hspace{6em}\cdot \lambda(w_{(1)}\otimes Y_2(Y(v_1, x_0)v_2, y_2)w_{(2)}).
\end{eqnarray}
Using (\ref{10.8}), evaluating
$\res_{y_{1}}$ and then using the definition of $Y'_{Q(z)}$ (recall
(\ref{Y'qdef})),
we finally see that the
right-hand side of (\ref{10.24}) is equal to
\begin{eqnarray}\label{10.25}
\lefteqn{x^{-1}_2\delta\left(\frac{x_1-x_0}{x_2}\right)\lambda(Y_1^o(Y(v_1,
x_0)v_2, x_2+z)w_{(1)}\otimes w_{(2)})}\nno\\
&&\quad -x^{-1}_2\delta\left(\frac{x_1-x_0}{x_2}\right)\res_{y_2}
x^{-1}_2\delta\left(\frac{z-y_2}{-x_2}\right)\cdot \nno\\
&&\hspace{6em}\cdot \lambda(w_{(1)}\otimes Y_2(Y(v_1, x_0)v_2, y_2)w_{(2)})\nno\\
&&=x^{-1}_2\delta\left(\frac{x_1-x_0}{x_2}\right)(Y'_{Q(z)}(Y(v_1,
x_0)v_2, x_2) \lambda)(w_{(1)}\otimes w_{(2)}),
\end{eqnarray}
proving Theorem \ref{6.1}.
\epfv

\noindent {\it Proof of Theorem \ref{6.2}} Let $\lambda$ be an element
of $(W_{1}\otimes W_{2})^{*}$ satisfying the $Q(z)$-compatibility
condition. We first want to prove that each coefficient in $x$ of
$Y'_{Q(z)}(u, x_0)Y'_{Q(z)}(v, x)\lambda$ is a formal Laurent series
involving only finitely many negative powers of $x_0$ and that
\begin{eqnarray}\label{11.1}
\lefteqn{\tau_{Q(z)}\left(z^{-1}\delta\left(\frac{x_{1}-x_0}{z}\right)
Y_t(u, x_0)\right)
Y'_{Q(z)}(v, x)\lambda}\nno\\
&&=z^{-1}\delta\left(\frac{x_{1}-x_0}{z}\right)
Y'_{Q(z)}(u, x_0)Y'_{Q(z)}(v,
x) \lambda
\end{eqnarray}
for all $u, v\in V$.
Using the commutator formula for $Y'_{Q(z)}$, we have
\begin{eqnarray}\label{11.2}
\lefteqn{Y'_{Q(z)}(u, x_0)Y'_{Q(z)}(v, x)\lambda}\nno\\
&&=Y'_{Q(z)}(v, x)Y'_{Q(z)}(u, x_0)\lambda\nno\\
&&\quad -\res_{y}x^{-1}_0\delta\left(\frac{x-y}{x_0}\right)
Y'_{Q(z)}(Y(v, y)u, x_0)\lambda.
\end{eqnarray}
Each coefficient in $x$ of the right-hand side of (\ref{11.2}) is a
formal Laurent series involving only finitely many negative powers of
$x_0$ since $\lambda$ satisfies the $Q(z)$-lower truncation condition.
Thus the coefficients in $x$ of $Y'_{Q(z)}(v, x)\lambda$ satisfy the
$Q(z)$-lower truncation condition.

By (\ref{5.2}) and (\ref{Y'qdef}), we have
\begin{eqnarray}\label{11.3}
\lefteqn{\left(\tau_{Q(z)}\left(z^{-1}\delta\left(\frac{x_{1}-x_0}{z}\right)
Y_t(u, x_0)\right)
Y'_{Q(z)}(v, x) \lambda\right)(w_{(1)}\otimes w_{(2)})}\nno\\
&&=x^{-1}_0\delta\left(\frac{x_1-z}{x_0}\right)(Y'_{Q(z)}(v, x)
\lambda)(Y_1^o(u, x_1)w_{(1)}\otimes w_{(2)})\nno\\
&&\quad -x^{-1}_0\delta\left(\frac{z-x_1}{-x_0}\right)(Y'_{Q(z)}(v, x)
\lambda)(w_{(1)}\otimes Y_2(u, x_1)w_{(2)})\nno\\
&&=x^{-1}_0\delta\left(\frac{x_1-z}{x_0}\right)\biggl(\lambda(Y_1^o(v,
x+z)Y_1^o(u, x_1)w_{(1)}\otimes w_{(2)})\nno\\
&&\quad \quad -\res_{x_2}x^{-1}\delta\left(\frac{z-x_2}{-x}\right)
\lambda(Y_1^o(u,
x_1)w_{(1)}\otimes Y_2(v, x_2)w_{(2)})\biggr)\nno\\
&&\quad -x^{-1}_0\delta\left(\frac{z-x_1}{-x_0}\right)\biggl(\lambda(Y_1^o(e^{zL(-
1)}v, x)w_{(1)}\otimes Y_2(u, x_1)w_{(2)})\nno\\
&&\quad \quad-\res_{x_2} x^{-1}\delta\left(\frac{z-x_2}{-x}\right)
\lambda(w_{(1)}\otimes
Y_2(v, x_2)Y_2(u, x_1)w_{(2)})\biggr).\;\;\;\;\;
\end{eqnarray}
Now the distributive law applies, giving us four terms. Inserting
\[
\res_{x_{4}}x_{4}^{-1}\delta\left(\frac{x+z}{x_{4}}\right)=1
\]
into the first of these terms and correspondingly replacing $x+z$ by $x_{4}$
in $Y^{o}_{1}(v, x+z)$, we can apply the commutator formula for $Y^{o}_{1}$
in the usual way. Also using the commutator formula for $Y_{2}$,
(\ref{2termdeltarelation}) and 
(\ref{deltafunctionsubstitutionformula}), we write the
 right-hand side of (\ref{11.3}) as
\begin{eqnarray}\label{11.4}
\lefteqn{x^{-1}_0\delta\left(\frac{x_1-z}{x_0}\right)\lambda(Y_1^o(u, x_1)Y_1^o(v,
x+z)w_{(1)}\otimes w_{(2)})}\nno\\
&&\quad -x^{-1}_0\delta\left(\frac{x_1-z}{x_0}\right)\res_{x_4}\res_{x_3}
x^{-1}_1\delta\left(\frac{x_{4}-x_3}{x_1}\right)
x^{-1}_4\delta\left(\frac{x+z}{x_4}\right)\cdot \nno\\
&&\hspace{6em}\cdot \lambda(Y_1^o(Y(v, x_3)u, x_1)w_{(1)}\otimes
w_{(2)})\nno\\
&&\quad -x^{-1}_0\delta\left(\frac{x_1-z}{x_0}\right)\res_{x_2}
x^{-1}\delta\left(\frac{z-x_2}{-x}\right)\lambda(Y_1^o(u, x_1)w_{(1)}\otimes
Y_2(v, x_2)w_{(2)})\nno\\
&&\quad -x^{-1}_0\delta\left(\frac{z-x_1}{-x_0}\right)\lambda(Y_1^o(e^{zL(-
1)}v, x)w_{(1)}\otimes Y_2(u, x_1)w_{(2)})\nno\\
&&\quad +x^{-1}_0\delta\left(\frac{z-x_1}{-x_0}\right)\res_{x_2}
x^{-1}\delta\left(\frac{z-x_2}{-x}\right)\lambda(w_{(1)}\otimes
Y_2(u, x_1)Y_2(v, x_2)w_{(2)})\nno\\
&&\quad +x^{-1}_0\delta\left(\frac{z-x_1}{-x_0}\right)\res_{x_2}
x^{-1}\delta\left(\frac{z-x_2}{-x}\right)\res_{x_3}
x^{-1}_1\delta\left(\frac{x_{2}-x_3}{x_1}\right)\cdot \nno\\
&&\hspace{6em}\cdot \lambda(w_{(1)}\otimes Y_2(Y(v, x_3)u, x_1)w_{(2)})\nno\\
&&=\Biggl(x^{-1}_0\delta\left(\frac{x_1-z}{x_0}\right)\lambda(Y_1^o(u,
x_1)Y_1^o(e^{zL(-1)}v, x)w_{(1)}\otimes w_{(2)})\nno\\
&&\quad \quad-x^{-1}_0\delta\left(\frac{z-x_1}{-x_0}\right)\lambda(Y_1^o(e^{zL(-
1)}v, x)w_{(1)}\otimes Y_2(u, x_1)w_{(2)})\Biggr)\nno\\
&&\quad -\res_{x_2} x^{-1}\delta\left(\frac{z-x_2}{-x}\right)
\Biggl(x^{-1}_0\delta\left(\frac{x_1-z}{x_0}\right)\cdot\nno\\
&&\hspace{6em}\cdot\lambda(Y_1^o(u,
x_1)w_{(1)}\otimes Y_2(v, x_2)w_{(2)})\nno\\
&&\quad \quad -x^{-1}_0\delta\left(\frac{z-x_1}{-x_0}\right)
\lambda(w_{(1)}\otimes
Y_2(u, x_1)Y_2(v, x_2)w_{(2)})\Biggr)\nno\\
&&\quad -\Biggl( \res_{x_4}\res_{x_3}(x_0+z)^{-1}\delta\left(\frac{x+z-x_3}{x_0+
z}\right) \cdot\nno\\
&&\hspace{5em}\cdot x^{-1}_4\delta\left(\frac{x+z}{x_4}\right)
x^{-1}_0\delta\left(\frac{x_1-z}{x_0}\right)\cdot \nno\\
&&\hspace{6em}\cdot \lambda(Y_1^o(Y(v, x_3)u, x_1)w_{(1)}\otimes w_{(2)})\nno\\
&&\quad \quad-\res_{x_2}\res_{x_3}x^{-1}\delta\left(\frac{z-(x_1+x_3)}{-
x}\right) x^{-1}_2\delta\left(\frac{x_1+x_3}{x_2}\right)\cdot \nno\\
&&\hspace{6em}\cdot x^{-1}_0
\delta\left(\frac{z-
x_1}{-x_0}\right)\lambda(w_{(1)}\otimes Y_2(Y(v, x_3)u, x_1)w_{(2)})\Biggr).
\end{eqnarray}

Using (\ref{5.2}) and (\ref{2termdeltarelation}) and evaluating
suitable residues, we see that the right-hand side of (\ref{11.4}) is
equal to
\begin{eqnarray}\label{11.5}
\lefteqn{\left(\tau_{Q(z)}\left(z^{-1}\delta\left(\frac{x_{1}-x_0}{z}\right)Y_t(u,
x_0)\right)\lambda\right)(Y_1^o(e^{zL(-1)}v, x)w_{(1)}\otimes w_{(2)})}\nno\\
&&\quad -\res_{x_2} x^{-1}\delta\left(\frac{z-x_2}{-x}\right)\cdot \nno\\
&&\hspace{4em}\cdot
\left(\tau_{Q(z)}\left(z^{-1}\delta\left(\frac{x_{1}-x_0}{z}\right)
Y_t(u, x_0)\right)\lambda\right)(w_{(1)}\otimes Y_2(v, x_2)w_{(2)})\nno\\
&&\quad -\biggl(\res_{x_3}x^{-1}_0
\delta\left(\frac{x-x_3}{x_0}\right)\cdot\nno\\
&&\hspace{4em}\cdot x^{-1}_0\delta\left(\frac{x_1-z}{x_0}\right)
\lambda(Y_1^o(Y(v, x_3)u,
x_1)w_{(1)}\otimes w_{(2)})\nno\\
&&\quad -\res_{x_3}x^{-1}\delta
\left(\frac{x_0+x_{3}}{x}\right)\cdot \nno\\
&&\hspace{4em}\cdot x^{-1}_0\delta\left(\frac{z-x_1}{-x_0}\right)
\lambda(w_{(1)}\otimes Y_2(Y(v, x_3)u, x_1)w_{(2)})\biggr).
\end{eqnarray}
Using (\ref{5.2}) and (\ref{Y'qdef}), we find that the right-hand side of
(\ref{11.5}) becomes
\begin{eqnarray}\label{11.6}
\lefteqn{\left(Y'_{Q(z)}(v, x)\tau_{Q(z)}\left(z^{-1}
\delta\left(\frac{x_{1}-x_0}{z}\right)
Y_t(u, x_0) \right)\lambda\right)(w_{(1)}\otimes w_{(2)})}\nno\\
&&\quad -\biggl( \res_{x_3} x^{-1}_0\delta\left(\frac{x-x_3}{x_0}\right)\cdot \nno\\
&&\hspace{3em}\cdot \tau_{Q(z)}\biggl(z^{-1}\delta\left(\frac{x_{1}-x_0}{z}
\right)
Y_t(Y(v, x_3)u,
x_0)\biggr) \lambda\biggr)(w_{(1)}\otimes w_{(2)}).\;\;\;\;\;
\end{eqnarray}
By the compatibility condition for $\lambda$ and the commutator formula for
$Y'_{Q(z)}$, the right-hand side of (\ref{11.6}) is
equal to
\begin{eqnarray}
\lefteqn{z^{-1}\delta\left(\frac{x_{1}-x_0}{z}\right)
(Y'_{Q(z)}(v, x)Y'_{Q(z)}(u,
x_0) \lambda)(w_{(1)}\otimes w_{(2)})}\nno\\
&&\quad -z^{-1}\delta\left(\frac{x_{1}-x_0}{z}\right)\biggl( \res_{x_3}
x^{-1}_0\delta\left(\frac{x-x_3}{x_0}\right)\cdot \nno\\
&&\hspace{6em}\cdot Y'_{Q(z)}(Y(v, x_3)u, x_0)
\lambda\biggr)(w_{(1)}\otimes w_{(2)})\nno\\
&&= z^{-1}\delta\left(\frac{x_{1}-x_0}{z}\right)
\biggl(\biggl(Y'_{Q(z)}(v, x)Y'_{Q(z)}(u, x_0)\nno\\
&&\quad -\res_{x_3} x^{-1}_0\delta\left(\frac{x-x_3}{x_0}\right)
Y'_{Q(z)}(Y(v, x_3)u,
x_0)\biggr) \lambda\biggr)(w_{(1)}\otimes w_{(2)})\nno\\
&&=z^{-1}\delta\left(\frac{x_{1}-x_0}{z}\right)
(Y'_{Q(z)}(u, x_0)Y'_{Q(z)}(v,
x) \lambda)(w_{(1)}\otimes w_{(2)}).
\end{eqnarray}
This proves (\ref{11.1}).  In the M\"obius case, the three operators
are handled in the usual way.  The first part of Theorem \ref{6.2} is
established.

The proof of the second half of Theorem \ref{6.2} is exactly like that
for Theorem \ref{stable}.
\epfv


\bigskip

\noindent {\small \sc Department of Mathematics, Rutgers University,
Piscataway, NJ 08854 (permanent address)}

\noindent {\it and}

\noindent {\small \sc Beijing International Center for Mathematical Research,
Peking University, Beijing, China}

\noindent {\em E-mail address}: yzhuang@math.rutgers.edu

\vspace{1em}

\noindent {\small \sc Department of Mathematics, Rutgers University,
Piscataway, NJ 08854}

\noindent {\em E-mail address}: lepowsky@math.rutgers.edu

\vspace{1em}

\noindent {\small \sc Department of Mathematics, Rutgers University,
Piscataway, NJ 08854}

\noindent {\em E-mail address}: linzhang@math.rutgers.edu

\end{document}